\documentclass[1p,authoryear]{elsarticle}

\usepackage{amsfonts}
\usepackage{amsmath}
\usepackage{amssymb}
\usepackage[applemac]{inputenc}
\usepackage{calrsfs} 

\usepackage{graphicx}
\usepackage{subfigure}
\usepackage{mathabx}

\usepackage{natbib}
\usepackage{geometry}

\newcommand{\vs}[1]{\vspace{#1}}    
\newcommand{\hs}[1]{\hspace{#1}}    
\newcommand{\vset}[1]{\vspace*{#1}}    
\newcommand{\hset}[1]{\hspace*{#1}}    
\newcommand{\begitem}{\begin{itemize}}    
\newcommand{\finit}{\end{itemize}}    
\newcommand{\begenum}{\begin{enumerate}}    
\newcommand{\finenum}{\end{enumerate}}    
\newcommand{\begar}{\begin{array}}    
\newcommand{\finar}{\end{array}}    
\newcommand{\begeq}[1]{\begin{equation} \label{#1}}    
\newcommand{\fineq}{\end{equation}}    
\newcommand{\begct}{\begin{center}}    
\newcommand{\finct}{\end{center}}    
\newcommand{\ra}{\rightarrow}    
\newcommand{\lra}{\longrightarrow}

\newcommand{\mkb}{\makebox}               
\newcommand{\zun}{\vs{0.1cm}} 
\newcommand{\zundemi}{\vs{0.1cm}}     
\newcommand{\zdeux}{\vs{0.2cm}} 
\newcommand{\ztrois}{\vs{0.3cm}}    
    
\newcommand{\tinf}{\rightarrow\infty}         
\newcommand{\tqdninf}{\stackrel{n\tinf}{\longrightarrow}} 

\newcommand{\bE}{\mathbb{E}}

\newcommand{\bR}{\mathbb{R}}
\newcommand{\bV}{\mathbb{V}}

\newcommand{\bP}{\mathbb{P}}

\newcommand{\bI}{\mathbb{I}}

\newcommand{\Var}{{\bV}ar}

\newcommand{\cause}{{\mathcal C}}
\newcommand{\ksi}{\xi}
\newcommand{\Zi}{Z_{(i)}}
\newcommand{\Zim}{Z_i^-}
\newcommand{\Zimun}{Z_{(i-1)}}
\newcommand{\deltai}{\delta_{(i)}}
\newcommand{\xpn}{{x}^{(k)}_{p_n}}
\newcommand{\xchapeau}{\hat {x}^{(k)}_{p_n,t_n}}

\newcommand{\indksik}{\bI_{\ksi_i=k}}
\newcommand{\indksipk}{\bI_{\ksi_{(i)}=k}}
\newcommand{\indcausek}{\bI_{\cause_i=k}}
\newcommand{\indZitn}{\bI_{Z_i>t_n}}

\newcommand{\Xin}{X_{i,n}}

\newcommand{\gamchapnk}{\widehat\gamma_{n,k}}
\newcommand{\gamtchenk}{\widecheck\gamma_{n,k}}
\newcommand{\gamtildenk}{\widetilde\gamma_{n,k}}
\newcommand{\gammank}{\gamma_{n,k}}

\newcommand{\racinevn}{\sqrt{v_n}}

\newcommand{\Fbar}{\widebar{F}}
\newcommand{\Gbar}{\widebar{G}}
\newcommand{\Hbar}{\widebar{H}}

\newcommand{\Fk}{F^{(k)}}
\newcommand{\Fbark}{\widebar{F}^{(k)}}
\newcommand{\Fnk}{F_n^{(k)}}
\newcommand{\Fbarnk}{\widebar F_n^{(k)}}

\newcommand{\Hbarn}{\widebar H_n}
\newcommand{\Gbarn}{\widebar G_n}
\newcommand{\GbarnZim}{\Gbarn(\Zim)}
\newcommand{\GbarnZimun}{\Gbarn(\Zimun)}

\newcommand{\Fbarnun}{\widebar F_n^{(1)}}

\newcommand{\HO}{H^{(0)}}
\newcommand{\HbarO}{\widebar H^{(0)}}
\newcommand{\HnO}{H_n^{(0)}}
\newcommand{\HbarnO}{\widebar H_n^{(0)}}

\newcommand{\Hunk}{H^{(1,k)}}

\newcommand{\Hnunk}{H_n^{(1,k)}}

\newcommand{\Fbarnktn}{\Fbarnk(t_n)}
\newcommand{\unsurFbarnktn}{\frac 1{\Fbarnk(t_n)}}
\newcommand{\Fbarktn}{\Fbark(t_n)}
\newcommand{\unsurFbarktn}{\frac 1{\Fbark(t_n)}}

\newcommand{\fracFbark}[2]{\frac{\Fbark(#1)}{\Fbark(#2)}}
\newcommand{\fracGbar}[2]{\frac{\Gbar(#1)}{\Gbar(#2)}}

\newcommand{\dfracFk}[1]{\frac{d\Fk(#1)}{\Fbark(t_n)}}
\newcommand{\dfracFbark}[1]{\frac{d\Fbark(#1)}{\Fbark(t_n)}}

\newcommand{\Uin}[1]{U^{(#1)}_{i,n}}
\newcommand{\Vin}[1]{V^{(#1)}_{i,n}}
\newcommand{\Ubarn}[1]{\widebar U_n^{(#1)}}

\newcommand{\Uunn}{U_{1,n}}
\newcommand{\Vunn}{V_{1,n}}
\newcommand{\Uunnb}[1]{U^{(#1)}_{1,n}}
\newcommand{\Vunnb}[1]{V^{(#1)}_{1,n}}
\newcommand{\tildeUin}{\tilde U_{i,n}}
\newcommand{\tildeUunn}{\tilde U_{1,n}}
\newcommand{\tildeVin}{\tilde V_{i,n}}
\newcommand{\tildeVunn}{\tilde V_{1,n}}

\newcommand{\Bin}{B_{i,n}}
\newcommand{\Cin}{C_{i,n}}
\newcommand{\Rn}[1]{R_{n,#1}}
\newcommand{\CC}[1]{C_n^{(#1)}}
\newcommand{\CCtilde}[1]{\widetilde C_n^{(#1)}}
\newcommand{\CCdoubletilde}[1]{\raisebox{0.27cm}{$\approx$}\hset{-0.27cm} C_n^{(#1)}}
\newcommand{\Cchap}[1]{\widehat C_n^{(#1)}}

\newcommand{\het}{\underline{h}}

\newcommand{\calUn}{{\cal U}_n}
\newcommand{\calVn}{{\cal V}_n}

\newcommand{\gdH}{{\cal H}}
\newcommand{\gdHet}{\underline{\cal H}}
\newcommand{\hunpt}{h_{1\bullet}}
\newcommand{\hptun}{h_{{\bullet}1}}
\newcommand{\hetunptpt}{\het_{1\bullet\bullet}}
\newcommand{\hetptunpt}{\het_{{\bullet}1{\bullet}}}
\newcommand{\hetptptun}{\het_{{\bullet\bullet}1}}

\newcommand{\gdHs}{{\gdH}^*}
\newcommand{\gdHss}{{\gdH}^{**}}
\newcommand{\gdHun}{{\gdH}_1}
\newcommand{\gdHsunpt}{\gdHs_{1\bullet}}
\newcommand{\gdHsptun}{\gdHs_{{\bullet}1}}
\newcommand{\gdHssunpt}{\gdHss_{1\bullet}}
\newcommand{\gdHssptun}{\gdHss_{{\bullet}1}}
\newcommand{\gdHunpt}{\gdH_{1\bullet}}
\newcommand{\gdHptun}{\gdH_{{\bullet}1}}
\newcommand{\sumI}{\sum_{I\in{\cal I}}}

\newcommand{\gdHets}{{\gdHet}^*}
\newcommand{\gdHetss}{{\gdHet}^{**}}
\newcommand{\gdHetun}{{\gdHet}_1}
\newcommand{\gdHetsunptpt}{\gdHets_{1\bullet\bullet}}
\newcommand{\gdHetsptunpt}{\gdHets_{{\bullet}1\bullet}}
\newcommand{\gdHetsptptun}{\gdHets_{{\bullet\bullet}1}}

\begin{document}

\newtheorem{theo}{Theorem}
\newtheorem{prop}{Proposition}
\newtheorem{defi}{Definition}
\newtheorem{lem}{Lemma}
\newtheorem{cor}{Corollary}
\newtheorem{rmk}{Remark}
  
\begin{center}
{\Large {\sc Extreme value statistics for censored data with heavy tails under competing risks }} 
\bigskip

\large  Julien Worms (1) \& Rym Worms\footnote{Corresponding author} (2)

\vset{1.cm}

\vset{1.cm}

 (1) Universit\'e Paris-Saclay / Universit\'e de Versailles-Saint-Quentin-En-Yvelines\\
 Laboratoire de Math\'ematiques de Versailles (CNRS UMR 8100), \\
 F-78035 Versailles Cedex, France, \\
 e-mail : {\tt julien.worms@uvsq.fr}
\bigskip
\vset{1.cm}

(2) Universit\'e Paris-Est \\
Laboratoire d'Analyse et de Math\'ematiques Appliqu\'ees \\
(CNRS UMR 8050), \\
 UPEMLV, UPEC, F-94010, Cr\'eteil, France, \\
 e-mail : {\tt rym.worms@u-pec.fr}
\end{center}
\vspace{1.cm}

 \newpage 
 
\begin{center}
{\Large {\sc Extreme value statistics for censored data with heavy tails under competing risks }} 
\bigskip

\end{center}
\vspace{1.cm}

\begin{center}
{\bf Abstract}\vspace{0.3cm}\\
\parbox{10.cm}{
This paper addresses the problem of estimating,  in the presence of random censoring as well as competing risks,  the extreme value index of the (sub)-distribution function associated to one particular cause,  in the heavy-tail case. Asymptotic normality of the proposed estimator (which has the form of an Aalen-Johansen integral, and is the first estimator proposed in this context) is established. A small simulation study exhibits its performances for finite samples. Estimation of extreme quantiles of the cumulative incidence function is also addressed.
} 
\end{center}

\vfill

  \vspace{0.1cm} 

\setlength{\belowdisplayskip}{0.2cm}
\setlength{\abovedisplayskip}{0.2cm}

\newpage

\section{Introduction}
 \label{intro}

The study of duration data (lifetime, failure time, re-employment time...) subject to random censoring is a major topic of the domain of statistics, which finds applications in many areas (in the sequel we will, for convenience, talk about lifetimes to refer to these observed durations, but without restricting our scope to lifetime data analysis). In general, the interest lies in obtaining informations about the central characteristics of the underlying lifetime distribution (mean lifetime or survival probabilities for instance), often with the objective of comparing results between different conditions under which the lifetime data are acquired. In this work, we will address the problem of inferring about the (upper)  tail of the lifetime distribution, for data subject both to random (right) censoring and competing risks.
\zdeux

Suppose indeed that we are interested in the lifetimes of $n$ individuals or items, which are subject to $K$ different causes of death or failure, and to random censorship (from the right) as well. We are particularly interested in one of these  causes (this main cause will be considered as cause number $k$ thereafter, where $k\in\{1,\ldots,K\}$), and we suppose that all causes are exclusive and are likely to be  dependent on the others. The censoring time is assumed to be independent of the different causes of death or failure and of the observed lifetime itself. However, since the other causes (different from the $k$-th cause of interest) generally cannot be considered as independent of the main cause, in no way they can be included in the censoring mechanism. This prevents us from relying on the basic independent censoring statistical framework, and we are thus in the presence of what is called a competing risks framework (see \cite{MoeschbergerKlein95}).\zdeux

For instance, if a patient is suffering from a very serious disease and starts some treatment, then the final outcome of the treatment can be death due to the main disease, or death due to other causes (nosocomial infection for instance). 
And censoring can occur due to loss of follow up or end of the clinical study. Another example, in a reliability experiment, is that the failure of some mechanical system can be due to the failure of a particular subpart, or component, of the system : since separating the different components for studying the reliability of only one of them is generally not possible, accounting for these different competing causes of failure is necessary.  Another field where competing risks often arise are labor economics, for instance in re-employment studies (see \cite{Fermanian03} for practical examples). 
\zdeux

One way of formalising this is to say that we observe a sample of $n$ independent couples $(Z_i,\ksi_i)_{1\leq i\leq n}$ where
\[
 Z_i = \min(X_i,C_i) , \hs{0.3cm}   \delta_i=\bI_{X_i\leq C_i} , \hs{0.3cm} \ksi_i=\left\{\begar{ll} 0 & \mbox{if $\delta_i=0$,} \\ \cause_i & \mbox{if $\delta_i=1$.} \finar\right. 
\]
The i.i.d. samples $(X_i)_{i\leq n}$ and $(C_i)_{i\leq n}$, of respective continuous distribution functions $F$ and $G$, represent the lifetimes and censoring times of the individuals, and are supposed to be independent. For convenience, we will suppose in this work that they are non-negative. The variables   $(\cause_i)_{i\leq n}$ form a discrete sample with values in $\{1,\ldots,K\}$, and represent the causes of failure or death of the $n$ individuals or items. It is important to note that these causes are observed only when the data is uncensored ({\it i.e.} when $\delta_i=1$), therefore we only observe the $\ksi_i$'s, not the complete $\cause_i$'s.  
\zun

One way of considering the failure times $X_i$ is to write 
\[
 X_i = \min(X_{i,1},\ldots,X_{i,K}), 
\]
where the variable $X_{i,k}$ is a (rather artificial) variable representing the imaginary latent lifetime of the $i$-th individual when the latter is only affected by the $k$-th cause (the other causes being absent). This viewpoint may be interesting in its own right, but we will not keep on considering it in the sequel, one reason being that such variables $X_{i,1},\ldots,X_{i,K}$ cannot be realistically considered as independent, and their respective distributions are of no practical use or interpretability (as explained and demonstrated in the competing risks literature, these distributions are in fact not statistically identifiable, see   \cite{Tsiatis75}  for example).  
\ztrois

The object of interest is the probability that a subject dies or fails after some given time $t$, due to the $k$-th cause, for high values of $t$. This quantity, denoted by
\[
\Fbark(t) = \bP\, [ \,  X > t \, , \, \cause =k  \, ] , 
\]
is related to the so-called {\it cumulative incidence function} $\Fk$ defined by
\[
\Fk(t) = \bP\, [ \,  X \leq t \, , \, \cause =k  \, ].
\]
Note that $\Fbark(t)$ is not equal to $1-\Fk(t)$, but to $\bP(\cause=k)-\Fk(t)$, because $\Fk$ is only a sub-distribution function. However we have  $ \Fbark(t)= \int_t^{\infty} d\Fk(u)$. In the sequel, the notation $\bar{S}(.)= S(\infty)- S(.)$ will be used, for any non-decreasing function $S$. 
\zdeux

In this paper, we are interested in investigating the behaviour of $\Fbark(t)$ for large values of $t$. This amounts to statistically study extreme values in a context of censored data under competing risks, and will lead us to consider some extreme value index $\gamma_k$ related to $\Fbark$, which will be defined in a few lines. Equivalently, the object of interest is the high quantile $x^{(k)}_p = (\Fbark)^-(p)=\inf \{ \, x \in \bR \, ; \, \Fbark(x) \geq p \, \}$ when $p$ is close to $0$, which can be interpreted as follows (in the context of lifetimes of individuals or failure times of systems) : in the presence of the other competing causes, a given individual (or item) will die (or fail), due to cause $k$ after such a time $x^{(k)}_p$, only with small probability $p$. A nonparametric inference for quantiles of fixed (and therefore not extreme) order, in the competing risk setting, has been already proposed in \cite{PengFine07}. 
\zdeux

One way of addressing this problem could be through a parametric point of view (see \cite{Crowder01} for further methods in the competing risk setting), however, the non-parametric approach is the most common choice of people faced with data presenting  censorship or competing risks. Of course, the standard Kaplan-Meier method for survival analysis does not yield valid results for a particular risk if failures from other causes are treated as censoring times, because the  other causes cannot always be considered  independent of the particular cause of interest. 
\zdeux

The commonly used  nonparametric estimator of the cumulative incidence function $\Fk$ is the so-called Aalen-Johansen estimator (see \cite{AalenJohansen78}, or \cite{Geffray09} equation $(7)$) defined by
\[
 \Fnk(t) = \sum_{Z_i \leq t} \frac{\delta_i\bI_{\cause_i=k}}{n\Gbarn(Z_i^-)}, 
\]
where $\Gbarn$ denotes the standard Kaplan-Meier estimator of $G$ (and $\Gbarn(t^-)$ denotes $\lim_{s\uparrow t}\Gbarn(s)$), so that we can introduce the following estimator for $\Fbark$ : 
\[
 \Fbarnk(t) = \sum_{Z_i > t} \frac{\delta_i\bI_{\cause_i=k}}{n\Gbarn(Z_i^-)} .
\]
But if the value $t$ considered  is so high that only very few  (if any) observations $Z_i$ (such that $\cause_i=k$) exceed $t$, then this purely nonparametric approach will lead to very unstable estimations $\Fbarnk(t)$ of $\Fbark(t)$. This is why a semiparametric approach is desirable, and the one we will consider here is the one inspired by classical extreme value theory. 
\zdeux

First note that in this paper, we will only consider situations where the underlying distributions $F$ and $G$ of the variables $X$ and $C$ are supposed to present power-like tails (also commonly named heavy tails), and we will focus on the evaluation of the order of this tail. 
Our working hypothesis will be thus that the different  functions $\Fbark$ (for $k=1,\ldots,K$) as well as $\Gbar=1-G$ belong to the Fr\'echet maximum domain of attraction. In other words, we assume that they are (see Definition \ref{defRV} in the Appendix)  regularly varying at infinity,  with respective negative indices $-1/\gamma_1,\ldots,-1/\gamma_K$ and $-1/\gamma_C$
\begeq{Ordre1}
 \forall 1\leq k\leq K, \ \ \forall x>0, \ \lim_{t\ra +\infty} \Fbark(tx)/\Fbark(t) = x^{-1/\gamma_k} 
 \makebox[1.5cm][c]{and} \lim_{t\ra +\infty} \Gbar(tx)/\Gbar(t) = x^{-1/\gamma_C}.
\fineq
Consequently, $\Fbar=1-F= \sum_{k=1}^K \Fbark$ and $\Hbar = \Fbar \Gbar$ (the survival function of  $Z$)  are  regularly varying (at $+\infty$) with  respective indices $-1/\gamma_F$ and $-1/\gamma$, where $\gamma_F=\max(\gamma_1, \ldots,  \gamma_K)$ and $\gamma$ satisfies $\gamma^{-1} ={\gamma_F}^{-1} + {\gamma_C}^{-1} $ (these relations are constantly used in this paper). 
\zdeux

The estimation of $\gamma_F$ has been already  studied in the literature, as it corresponds to the random (right) censoring framework, without competing risks. We can cite \cite{Beirlant07} and \cite{Einmahl08}, where the authors propose to use consistent estimators of $\gamma$ divided by the proportion of non-censored  observations in the tail, or \cite{WormsWorms14}, where  two Hill-type  estimators are proposed for $\gamma_F$, based on  survival analysis techniques.  However, our target here is $\gamma_k$ (for a fixed $k=1,\ldots,K$) and the point is that there seems to be no way to deduce an estimator of  $\gamma_k$ from an estimator of $\gamma_F$. Note that the useful trick used in \cite{Beirlant07} and \cite{Einmahl08} to construct an estimator of $\gamma_F$ does not seem to  be extendable to  this competing risks setting.  To the best of our knowledge, our present  paper  is the first one addressing the problem of estimating the cause-specific extreme value index $\gamma_k$ . 
\zdeux

Considering assumption (\ref{Ordre1}), it is simple to check that, for a given $k$, we have 
\[
  \lim_{t\rightarrow +\infty} \frac{1}{\Fbark(t)} \int_{t}^{+\infty} \log(u/t) \, d\Fk(t) \, = \, \gamma_k.
\]
It is therefore most natural to propose the following (Hill-type)  estimator of $\gamma_k$ , for some given threshold value $t_n$ (assumptions on this threshold are detailed in the next section) : 
 \[
  \gamchapnk \; = \; \int \widehat\phi_n(u) d\Fnk(u)
  \makebox[1.8cm][c]{where} 
  \widehat\phi_n(u) \; = \; \unsurFbarnktn \log \left( \frac u {t_n} \right)  \bI_{u>t_n} ,
 \]
 which can be also written as
\[
 \gamchapnk = \frac 1 {n \Fbarnktn} \sum_{i=1}^n \frac{\log(Z_i/t_n)  }{ \GbarnZim }  \indksik \indZitn
 = \frac 1 {n \Fbarnktn} \sum_{ \Zi > t_n} \frac{\log(\Zi/t_n)}{ \GbarnZimun } \deltai \bI_{\cause_{(i)}=k} , 
\]
where $Z_{(1)} \leq \ldots \leq Z_{(n)}$ are the ordered random variables associated to $Z_1, \ldots, Z_n$, 
and $\delta_{(i)}$  and $\cause_{(i)}$ are the censoring indicator and cause number which correspond to the order statistic $Z_{(i)}$. 
It is clear that this estimator is a generalisation of one of the estimators proposed in \cite{WormsWorms14}, in which the situation $K=1$ (with only one cause of failure/death) was considered. The asymptotic result we prove in the present work is then valid in the  situation studied in the latter, where only consistency was proved and a random threshold was used.  
\zdeux

Our paper is organized as follows: in Section \ref{HypoTLC}, we state the asymptotic normality result of the proposed estimator, and of a corresponding estimator of  an extreme quantile of the cumulative incidence function. Section \ref{sec-proofs} is devoted to the proofs. In Section \ref{simuls}, we  present some simulations in order to illustrate finite sample behaviour of our estimator. Some technical aspects of the proofs are postponed to the Appendix.

\section{Assumptions and  Statement of the results}
\label{HypoTLC}
The central limit theorem  which is going to be proved has the rate $\sqrt{v_n} $ where $v_n= n \Fbarnktn\Gbar(t_n)$ and $t_n$ is a threshold tending to $\infty$ with the following constraint 
\begeq{condvntn}
v_n \tqdninf +\infty  \makebox[2.2cm][c]{  such that  }  n^{-\eta_0} v_n \tqdninf + \infty \mbox{ \ \ for some \ }   \eta_0 >0.
\fineq
If we note  $l_k$ the slowly varying function associated to $\Fbark$ ({\it i.e.} such that  $\Fbark (x)= x^{-1/\gamma_k} l_k(x)$ in condition (\ref{Ordre1})), the second order condition we consider is the classical $SR2$ condition for $l_k$ (see \citet{BinghamGoldieTeugels87}),
\begeq{Ordre2}
\forall x >0 , \ \frac{l_k(tx)}{l_k(t)} -1 \ \stackrel{t\tinf}{\sim} \ h_{\rho_k}(x) \  g(t) \hspace{0.3cm} (\forall x >1),
\fineq
where $g$ is a positive measurable function, slowly varying with index $\rho_k \leq 0$, and $h_{\rho_k}(x) = \frac{x^{\rho_k}-1}{\rho_k}$  when  $\rho_k <0$, or  $h_{\rho_k}(x) =\log x$ when  $\rho_k = 0$.


\begin{theo} 
\label{TLCGammachap}
Under assumptions $(\ref{Ordre1})$,  $(\ref{condvntn})$ and $(\ref{Ordre2})$,  if   there exists $\lambda  \geq 0$  such that  $\racinevn g(t_n) \tqdninf  \lambda $, and if $\gamma_k < \gamma_C$ 
then  we  have 
\[
\racinevn  ( \gamchapnk  - \gamma_k )  \ \stackrel{d}{\longrightarrow} \  {\cal N}(\lambda m,\sigma^2) \hs{0.5cm} \mbox{as $n\tinf$}
\]
where 
\[
m= \left\{ \begar{ll}    \frac{\gamma_k^2}{1-\gamma_k \rho_k} & \mbox{ if } \rho_k <0, \zdeux\\  \gamma_k^2 & \mbox{ if } \rho_k = 0, \finar \right.  \mbox{ and  } \  \sigma^2 =   \frac{\gamma_k^2}{(1-r)^3} \left( (1+r^2) - 2cr \right),  
 \]
with $c=\lim_{x\tinf} \Fbark(x)/\Fbar(x) \in [0,1]$ and $r=\gamma_k /\gamma_C \in ]0,1[$. 
\end{theo}

\begin{rmk}
Note that when $\gamma_k  < \gamma_F$, then $c=0$, and,  when $\gamma_k=\gamma_F$ and $c=1$ (for instance when there is only one cause of failure/death), then $\sigma^2$ reduces to $\gamma_F^2/(1-r)$.
\end{rmk}


\begin{prop}
\label{consistance}
Under assumptions $(\ref{Ordre1})$  and $(\ref{condvntn})$, we have  
\[
\gamchapnk    \ \stackrel{\bP}{\longrightarrow} \  \gamma_k \hs{0.5cm} \mbox{as $n\tinf$}.
\]
\end{prop}

\begin{rmk}
The condition $\gamma_k < \gamma_C$ (weak censoring) is  not necessary for the consistency of $\gamchapnk$.
\end{rmk}

Now, concerning the estimation of an extreme quantile $x^{(k)}_{p_n}$ (of order  $p_n$ tending to $0$) associated to   $\Fbark$, we propose the usual Weissman-type estimator (in this heavy tailed context), associated to the threshold $t_n$ used in the estimation of $\gamma_k$,
\[
 \xchapeau = t_n \left(  \frac{\Fbarnk(t_n) }{ p_n }  \right)^{\gamchapnk} ,
\]
where $p_n$ is assumed to satisfy the constraint $p_n = o\big( \Fbark(t_n) \big)$. Remind that by definition $\Fbark( {x}^{(k)}_{p_n} ) = p_n$, and thus the definition of this estimator is based on the fact that, by the assumed regular variation of $\Fbark$, the ratio $\Fbark( {x}^{(k)}_{p_n} ) / \Fbark(t_n)$ is close to $({x}^{(k)}_{p_n}/t_n)^{-1/\gamma_k}$. 

\begin{cor}\label{coroquantiles}
Under the assumptions of Theorem \ref{TLCGammachap}, if in addition $\rho_k<0$ (in (\ref{Ordre2})) and $d_n=   \Fbark(t_n)/p_n\tinf$ satisfies the condition 
\begeq{conditiondn}
 \sqrt{v_n}\;/\log(d_n) \; \tqdninf \; \infty,
\fineq
then (with $\lambda$, $m$ and $\sigma^2$ being defined in the statement of Theorem \ref{TLCGammachap})
\[
 \frac{\sqrt{v_n}}{\log(d_n)} \left( \frac{ \xchapeau }{ {x}^{(k)}_{p_n} } - 1 \right) 
 \ \stackrel{d}{\longrightarrow} \  {\cal N}(\lambda m,\sigma^2) \hs{0.5cm} \mbox{as $n\tinf$}.
\]
\end{cor}

\section{Simulations}  \label{simuls}

In this section, a small simulation study is conducted in order to illustrate the finite-sample behaviour of our new estimator in some simple cases, and discuss the main issues associated with the competing risks setting. 
\zun

\begin{figure}[hbtp]
\centering
\subfigure[Fr\'echet case, \, $\gamma_1=0.1,\; \gamma_2=0.25,\; \gamma_C=0.3$]{
\label{CadreFrechet1-1}
\includegraphics[height=3.5cm,width=.46\textwidth]{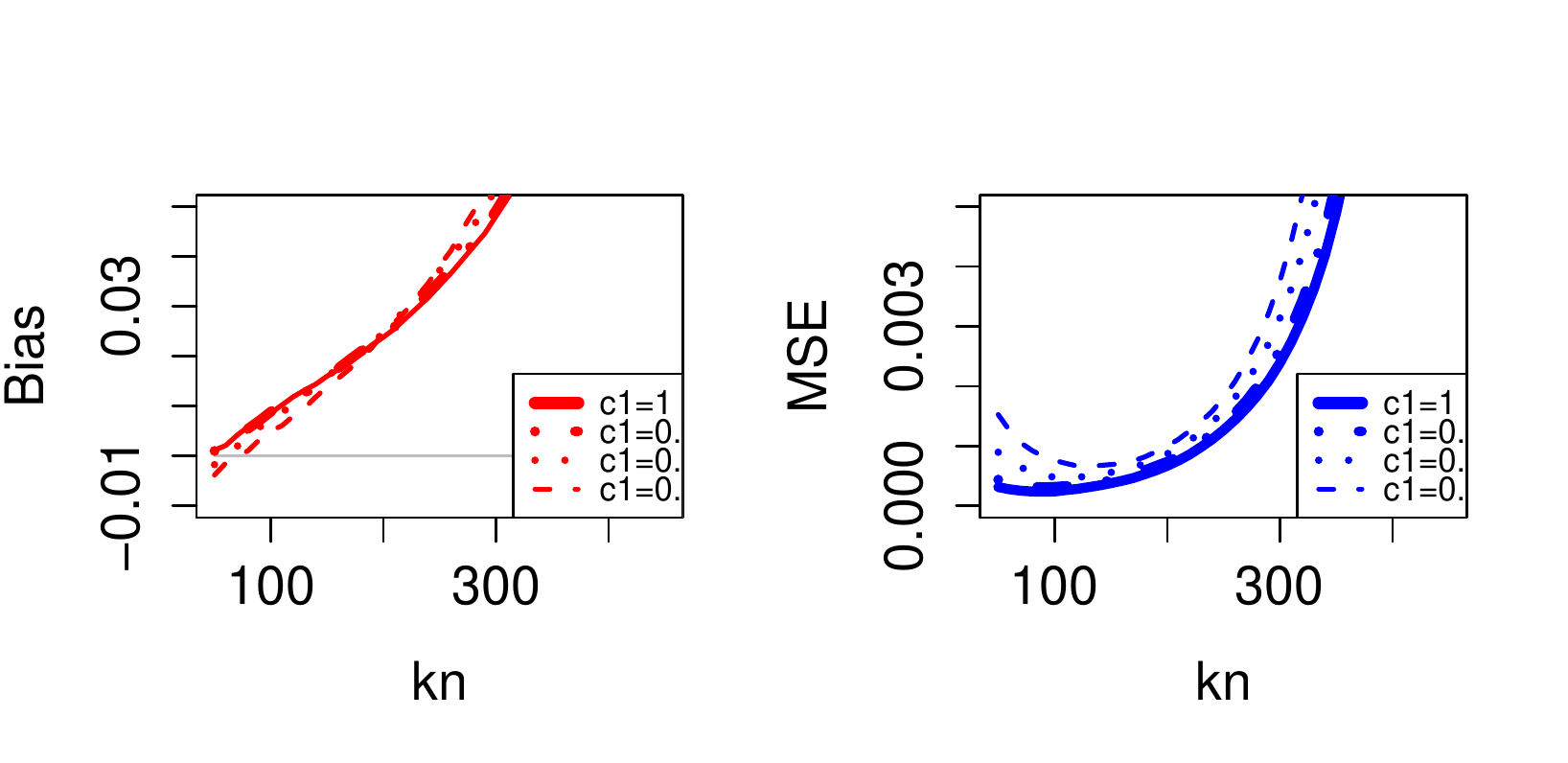}}
\hspace{.1in}
\subfigure[Same case as (a) but for Burr distribution]{
\label{CadreFrechet1-09}
\includegraphics[height=3.5cm,width=.46\textwidth]{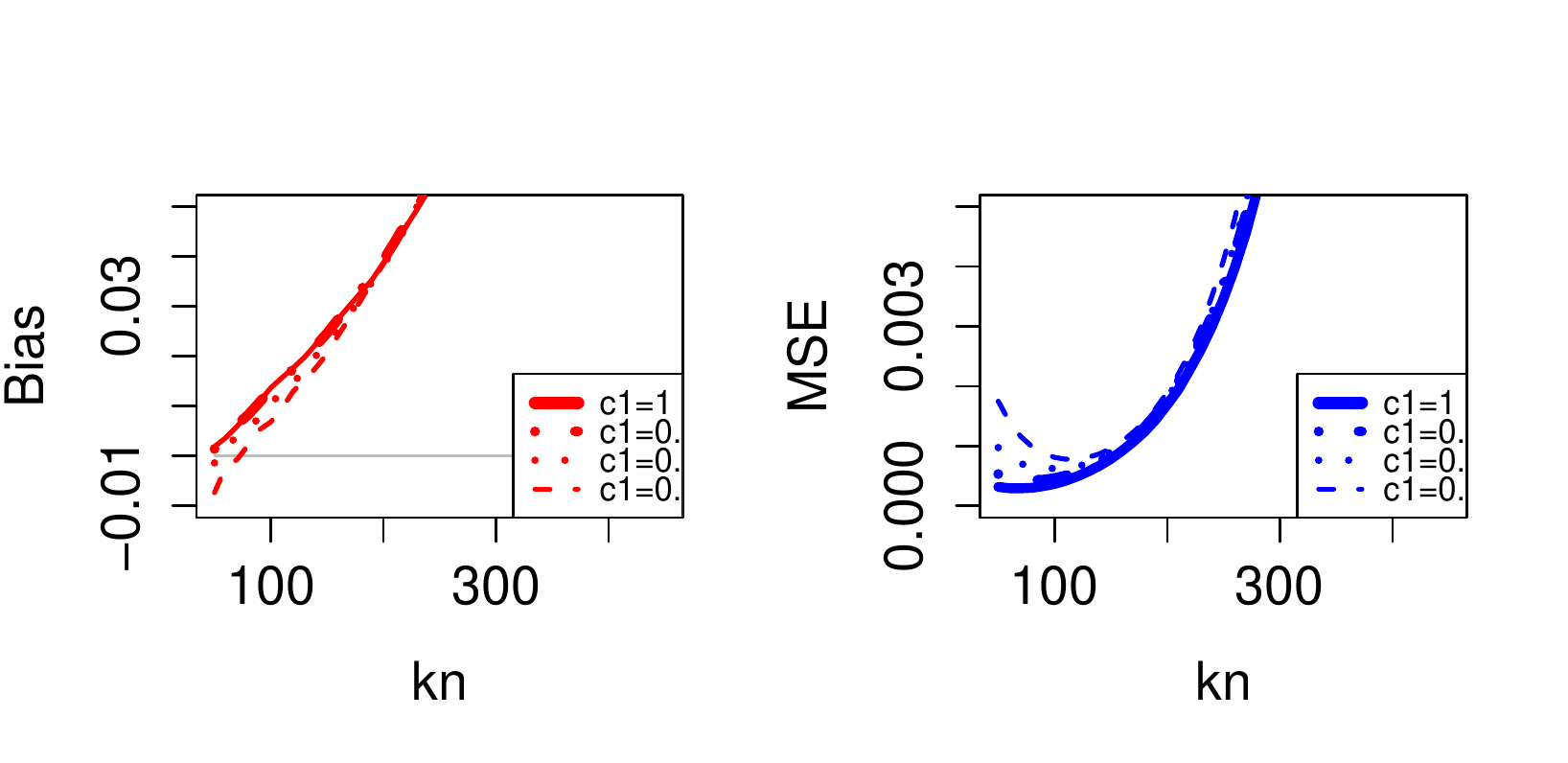}}
\hspace{.1in}
\subfigure[Fr\'echet case, \, $\gamma_1=0.1,\; \gamma_2=0.25,\; \gamma_C=0.2$]{
\label{CadreFrechet2-1}
\includegraphics[height=3.5cm,width=.46\textwidth]{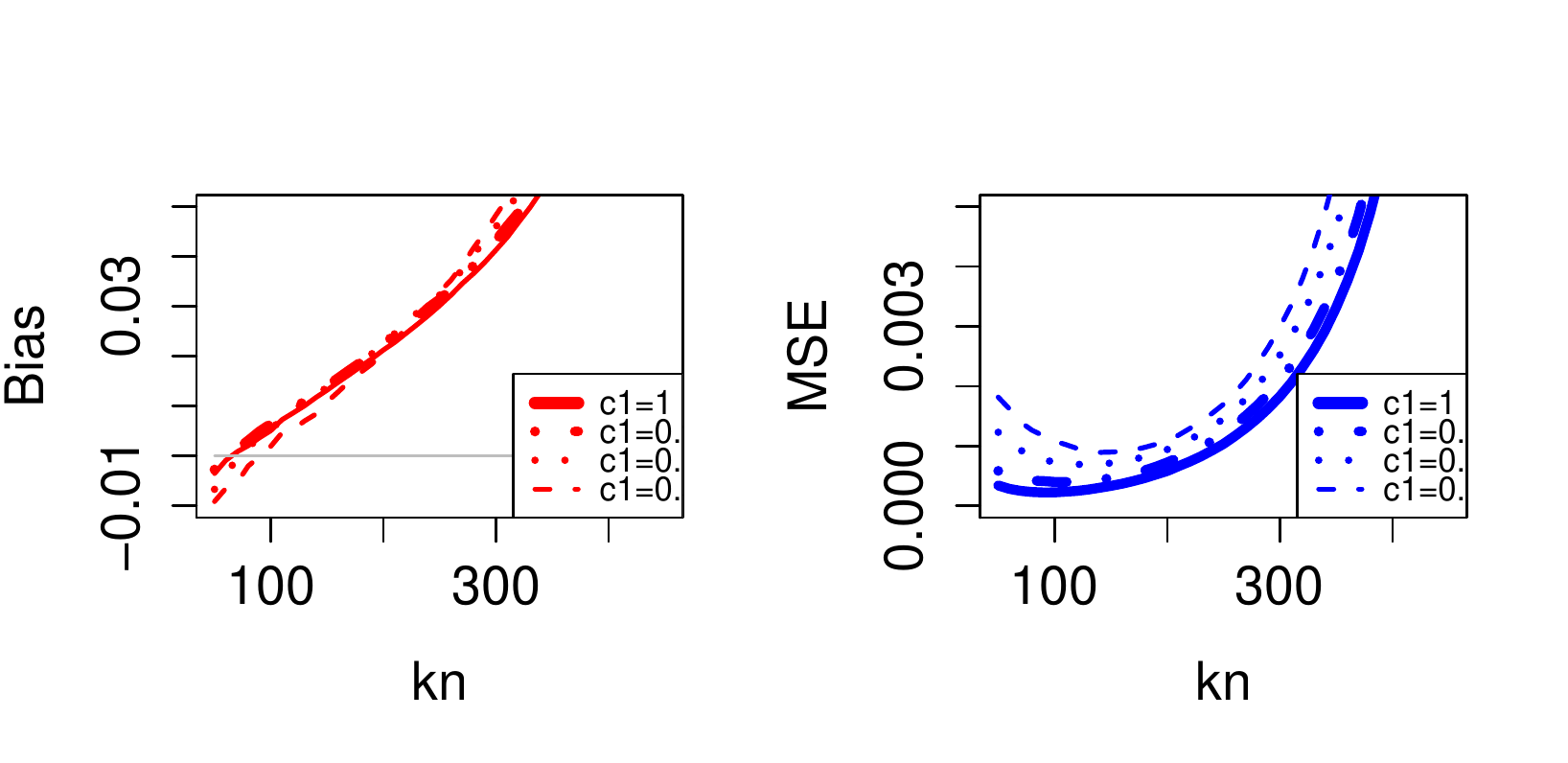}}
\hspace{.1in}
\subfigure[Same case as (c) but for Burr distribution]{
\label{CadreFrechet2-09}
\includegraphics[height=3.5cm,width=.46\textwidth]{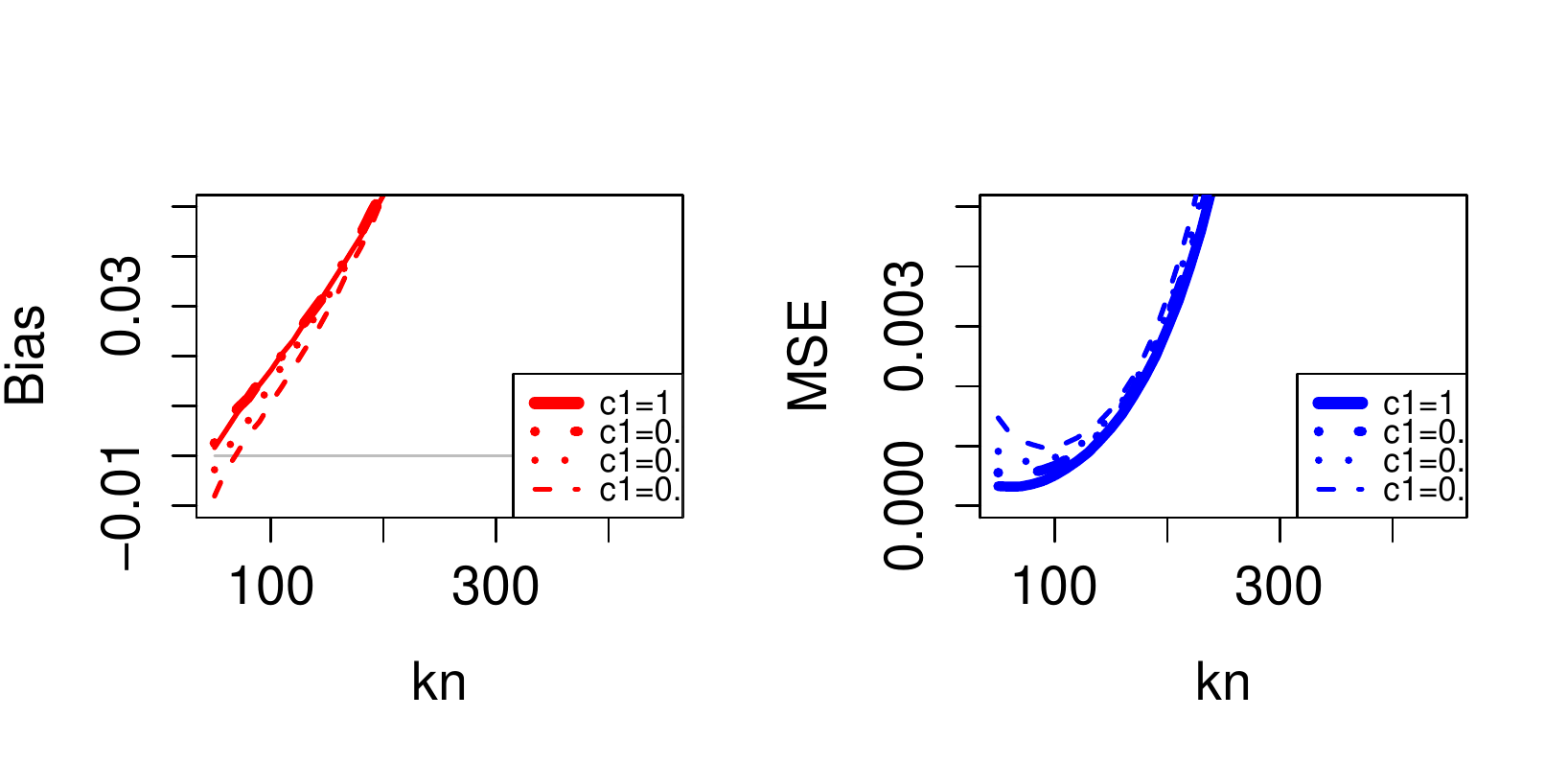}}
\hspace{.1in}
\subfigure[Fr\'echet case, \, $\gamma_1=0.25,\; \gamma_2=0.1,\; \gamma_C=0.45$]{
\label{CadreFrechet3-1}
\includegraphics[height=3.5cm,width=.46\textwidth]{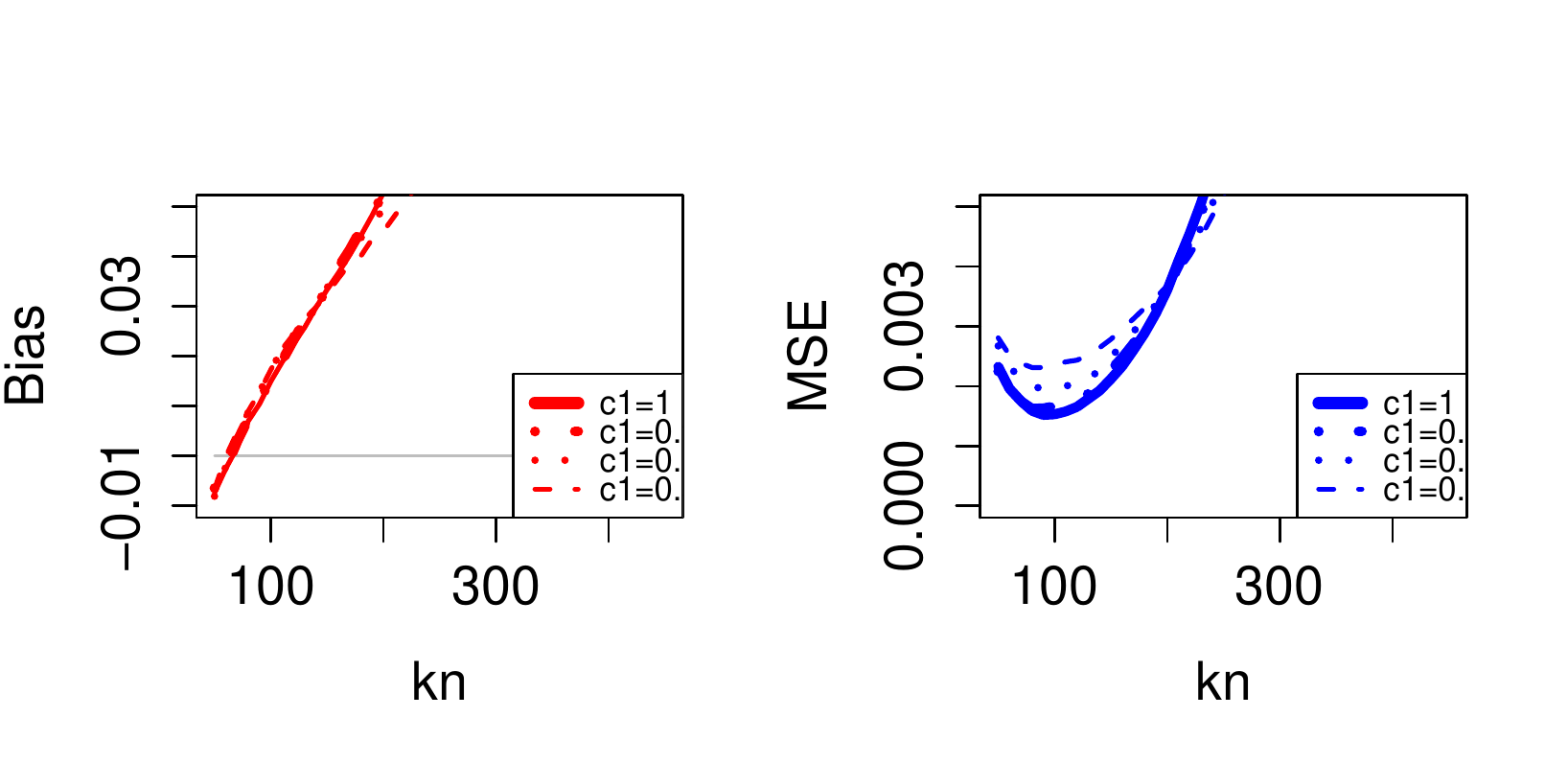}}
\hspace{.1in}
\subfigure[Same case as (e) but for Burr distribution]{
\label{CadreFrechet3-09}
\includegraphics[height=3.5cm,width=.46\textwidth]{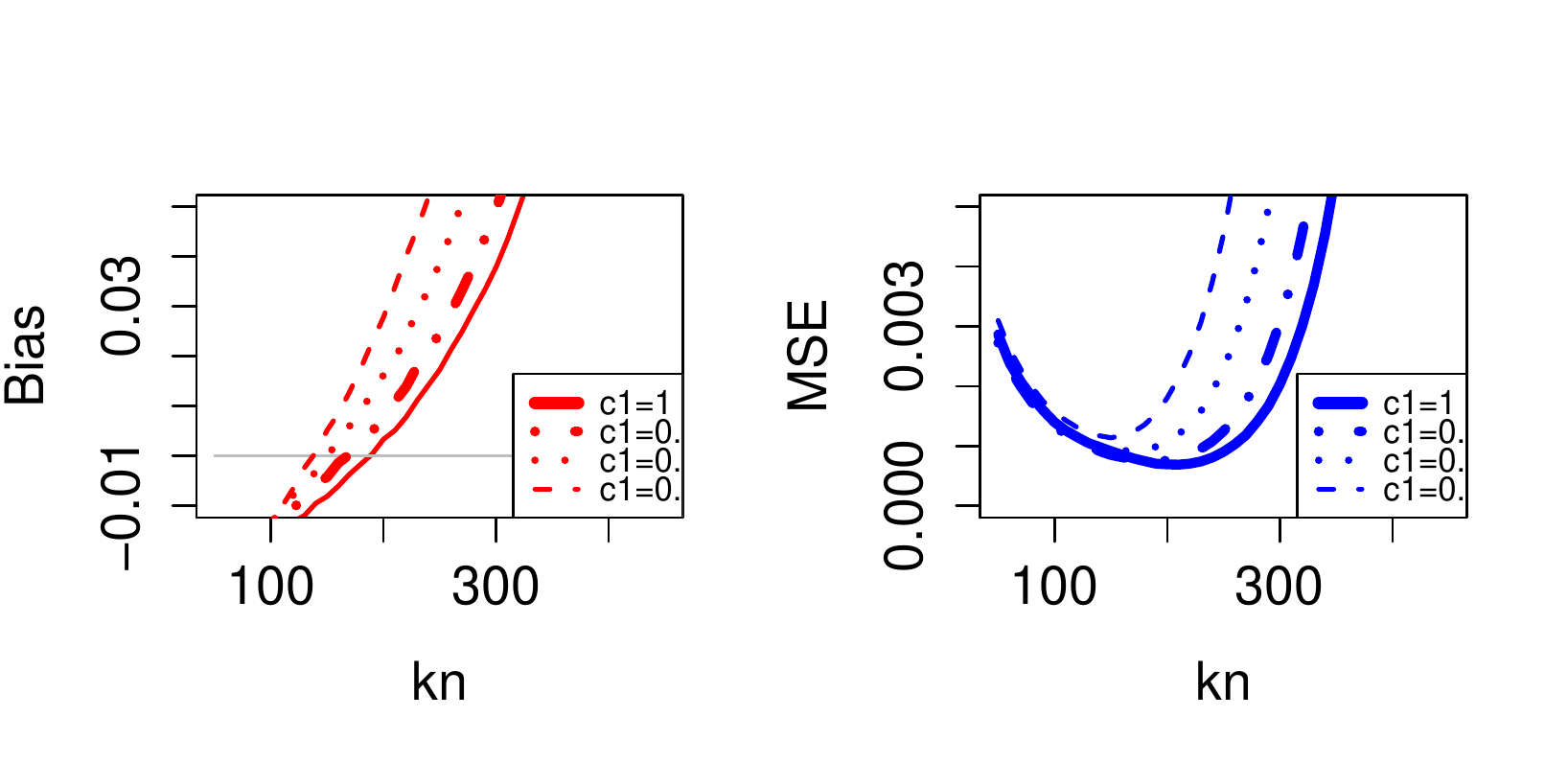}}
\hspace{.1in}
\caption{Comparison of bias  and MSE (respectively left and right in each subfigure) of $\gamchapnk$ for different values of $c_1$ ; in figures $(a)$, $(c)$ and $(e)$, $X$ and $C$ are Fr\'echet distributed, but in figures $(b)$, $(d)$ and $(f)$ they are Burr distributed. }
\label{figinfluencec}
\end{figure}

For simplicity, we focus on the situation with two competing risks ($K=2$), also called causes below, and our aim is the extreme value index $\gamma_1$ associated to the first cause. Data are generated from one of the following two models : for $c_1$, $c_2$ non-negative constants satisfying $c_1+c_2=1$, we consider the following (sub-)distribution for each cause-specific function $\Fbark$ ($k\in\{1,2\}$) : 
\zun\\
\hspace*{0.3cm} $-$ Fr\'echet : $\Fbark(t)= c_k \ \exp(-t^{-1/\gamma_k}) $, for $t \geq 0$ ;  \zun\\
\hspace*{0.3cm} $-$ Burr  :   $\Fbark(t)= c_k  \ (1+t^{\tau_k}/\beta)^{-1/(\gamma_k\tau_k)} $, for $t \geq 1$, where $\tau_k>0$, $\beta>0$.   
\zdeux\\
The lifetime $X$, of survival function $\bar F=\bar F^{(1)}+\bar F^{(2)}$,  is generated by the inversion method (with numerical computation of $\Fbar^{-1}$). Censoring times are then generated  from a Fr\'echet or a Burr distribution : 
\[
\Gbar(t)=  \exp(-t^{-1/\gamma_C})  \; (t\geq 0) 
\makebox[1.4cm][c]{ or } \Gbar(t)= (1+t^{\tau_C}/\beta)^{-1/(\gamma_C\tau_C)}  \; (t\geq 1) . 
\zun
\]

\noindent In this section,  we consider (as it is often done in simulation studies) that the threshold $t_n$ used in the definition of our new estimator $\hat{\gamma}_{n,1} $ is taken equal to $Z_{(n-k_n)}$ (i.e. we consider it as random).   One aim of this section is to show how our  estimator  (with random threshold)
\[
\hat{\gamma}_1  =  \frac 1 {n\Fbarnun(Z_{(n-k_n)})} \sum_{i=1}^{k_n} \frac{\log(Z_{(n-i+1)}/Z_{(n-k_n)})}{\Gbarn(Z_{(n-i,n)})}  \delta_{(n-i+1)}\bI_{\cause_{(n-i+1)}=1}  
\]
 of $\gamma_1$ behaves when the proportion $c_1$ of  cause $1$ events varies : we consider $c_1\in\{1,0.9,0.7,0.5\}$, the case $c_1=1$ corresponding to the simple censoring framework, without competing risk. 
\zun

\begin{figure}[hbt]
\centering
\subfigure[Fr\'echet case,\; $\gamma_1=0.1,\, \gamma_2=0.25,\, \gamma_C=0.3$, and $c_1=1$]{
\label{CadreFrechet1-1}
\includegraphics[height=3.5cm,width=.46\textwidth]{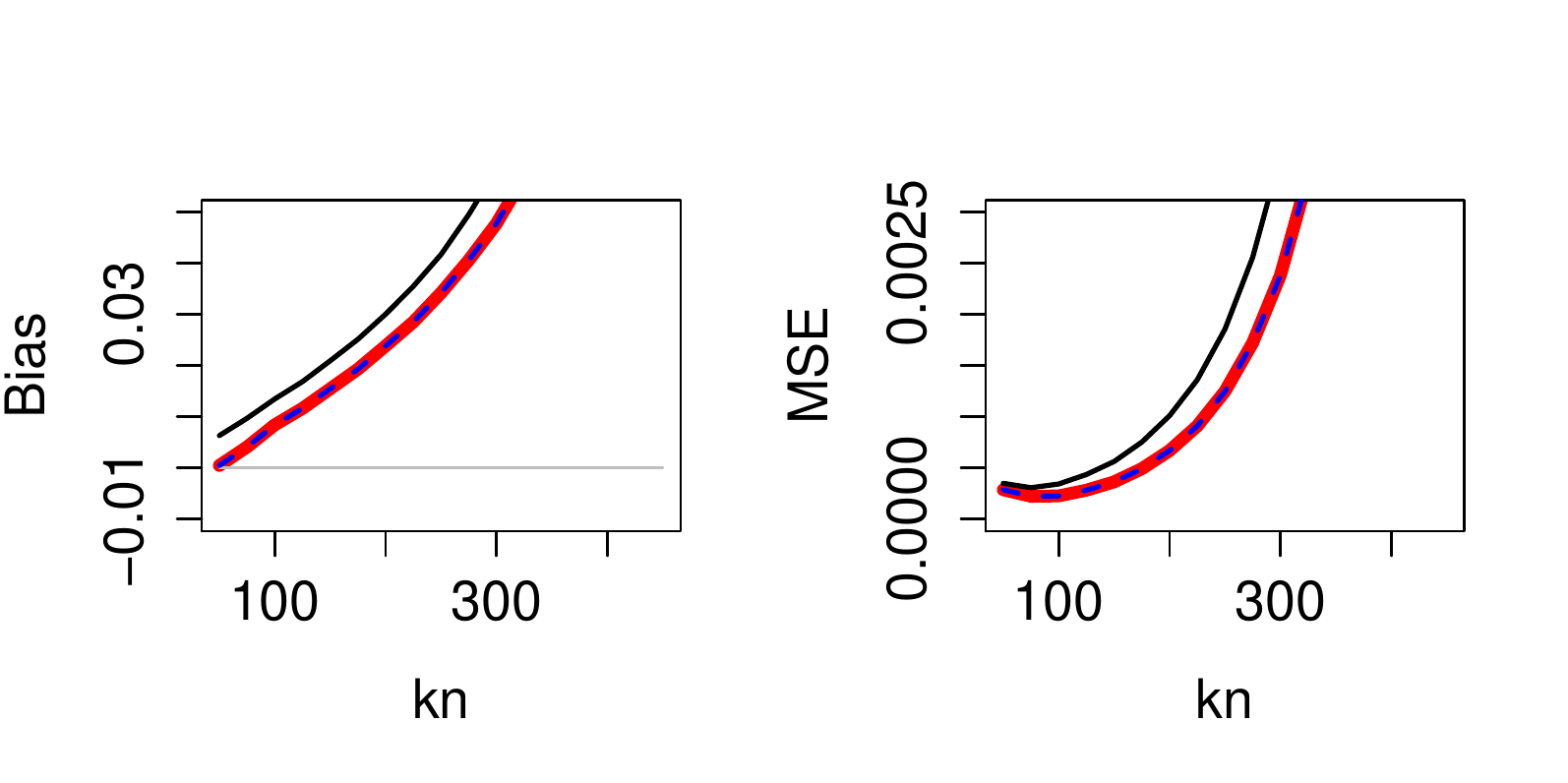}}
\hspace{.1in}
\subfigure[Case (a) but with $c_1=0.9$]{
\label{CadreFrechet1-09}
\includegraphics[height=3.5cm,width=.46\textwidth]{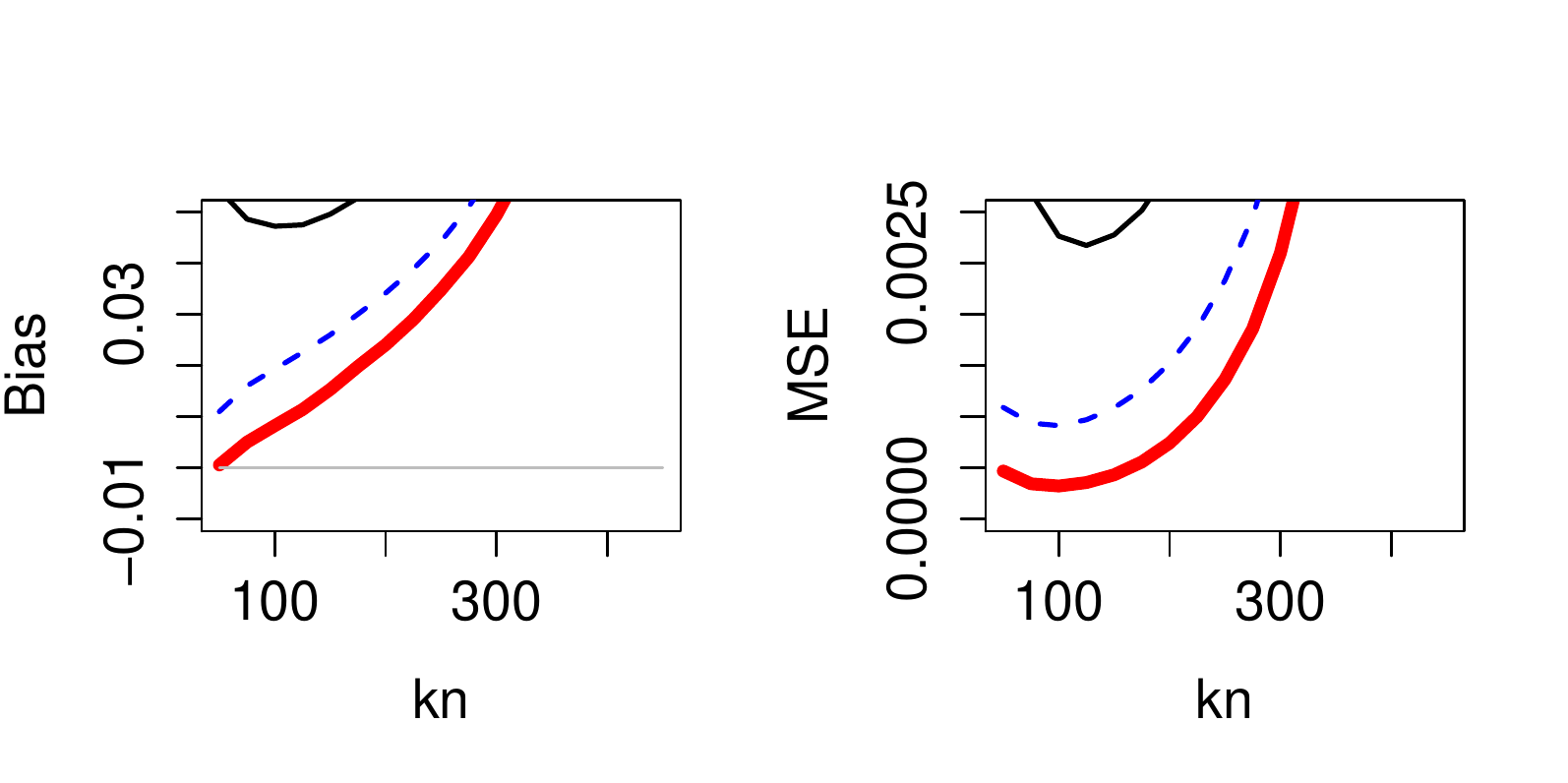}}
\hspace{.1in}
\subfigure[Fr\'echet case,\; $\gamma_1=0.1,\, \gamma_2=0.25,\, \gamma_C=0.2$, and  $c_1=1$]{
\label{CadreFrechet2-1}
\includegraphics[height=3.5cm,width=.46\textwidth]{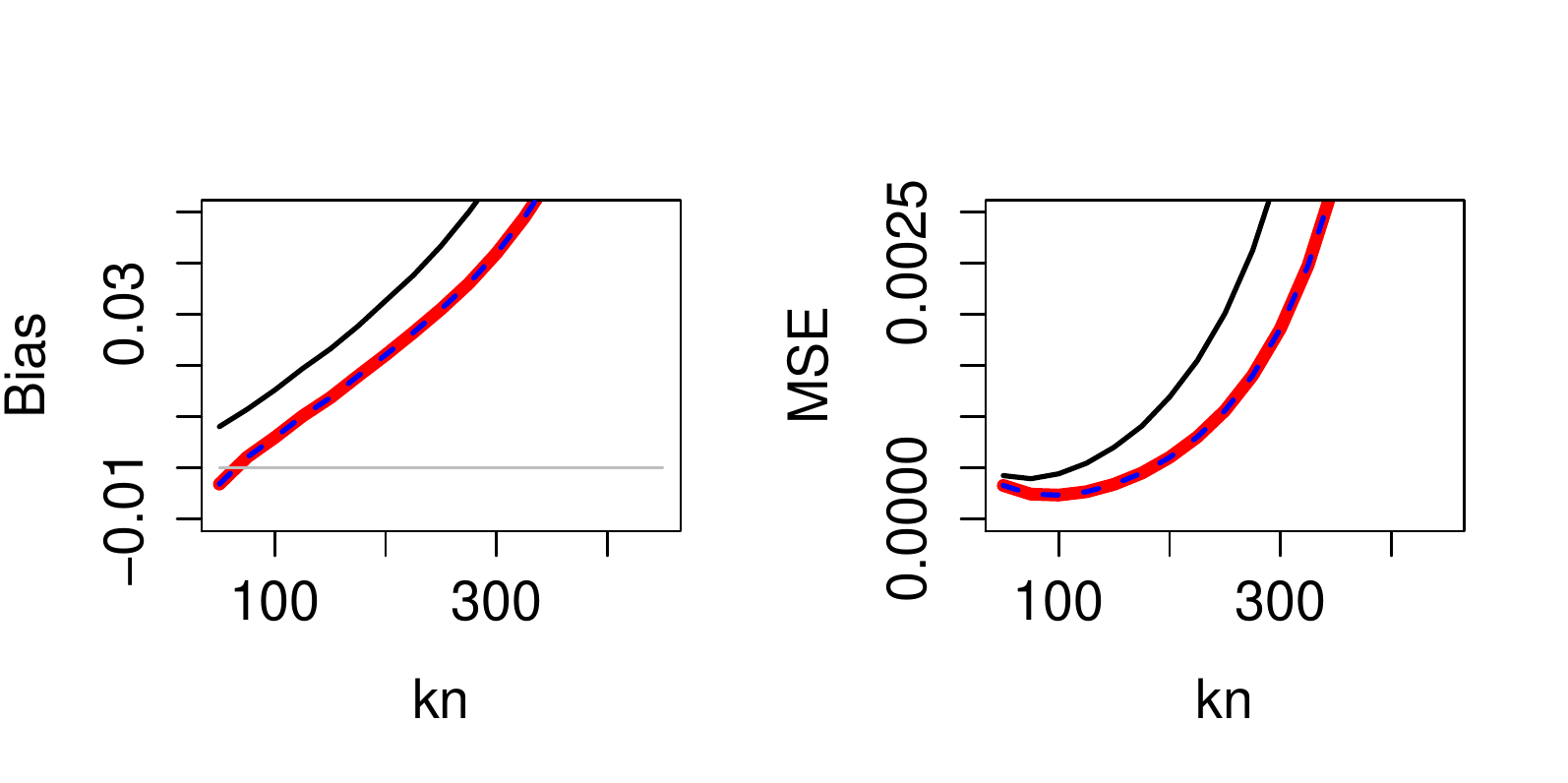}}
\hspace{.1in}
\subfigure[Case (c) but with $c_1=0.9$]{
\label{CadreFrechet2-09}
\includegraphics[height=3.5cm,width=.46\textwidth]{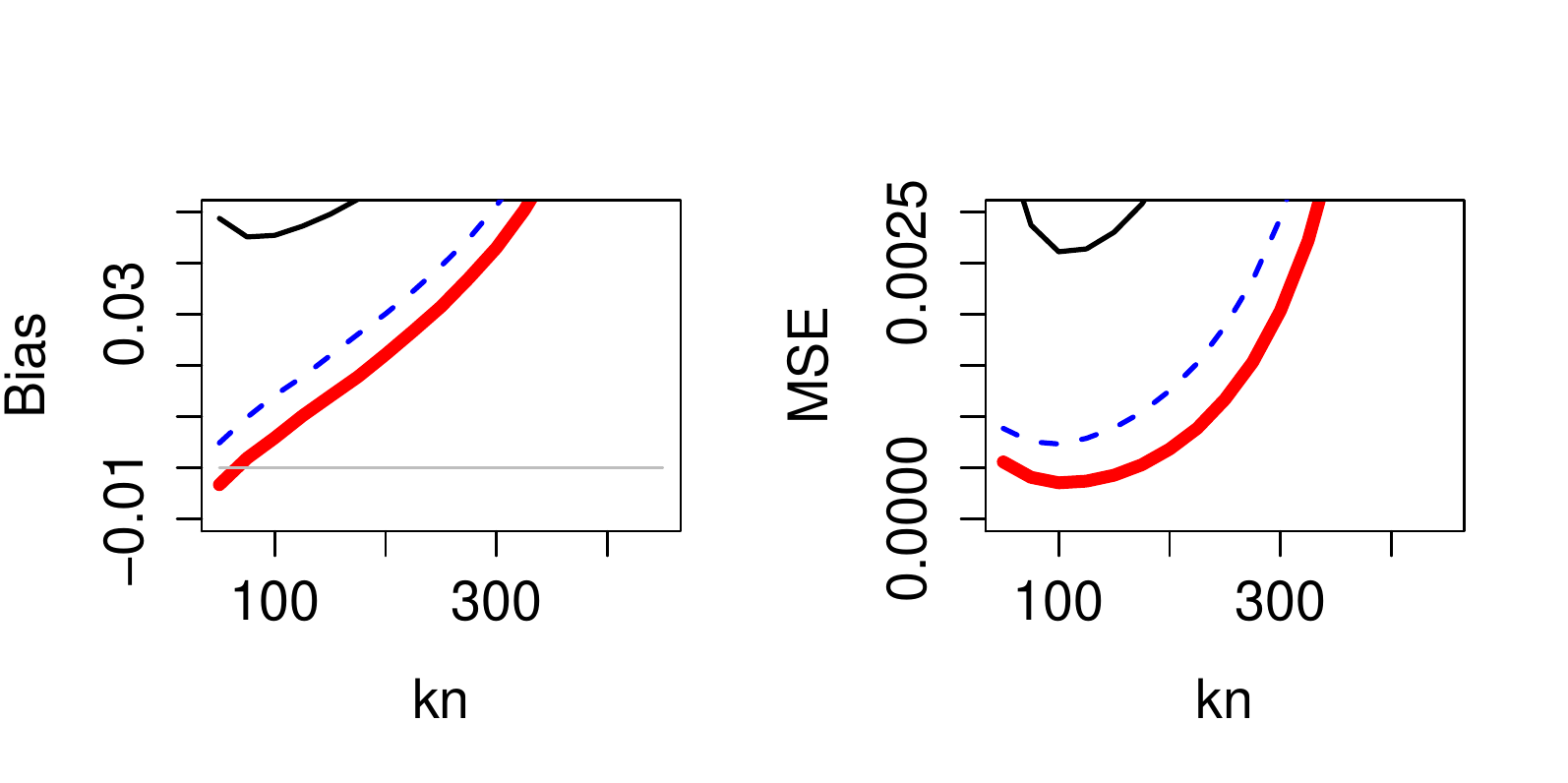}}
\hspace{.1in}
\subfigure[Fr\'echet case,\, $\gamma_1=0.1, \gamma_2=0.25, \gamma_C=0.45$, and $c_1=1$]{
\label{CadreFrechet3-1}
\includegraphics[height=3.5cm,width=.46\textwidth]{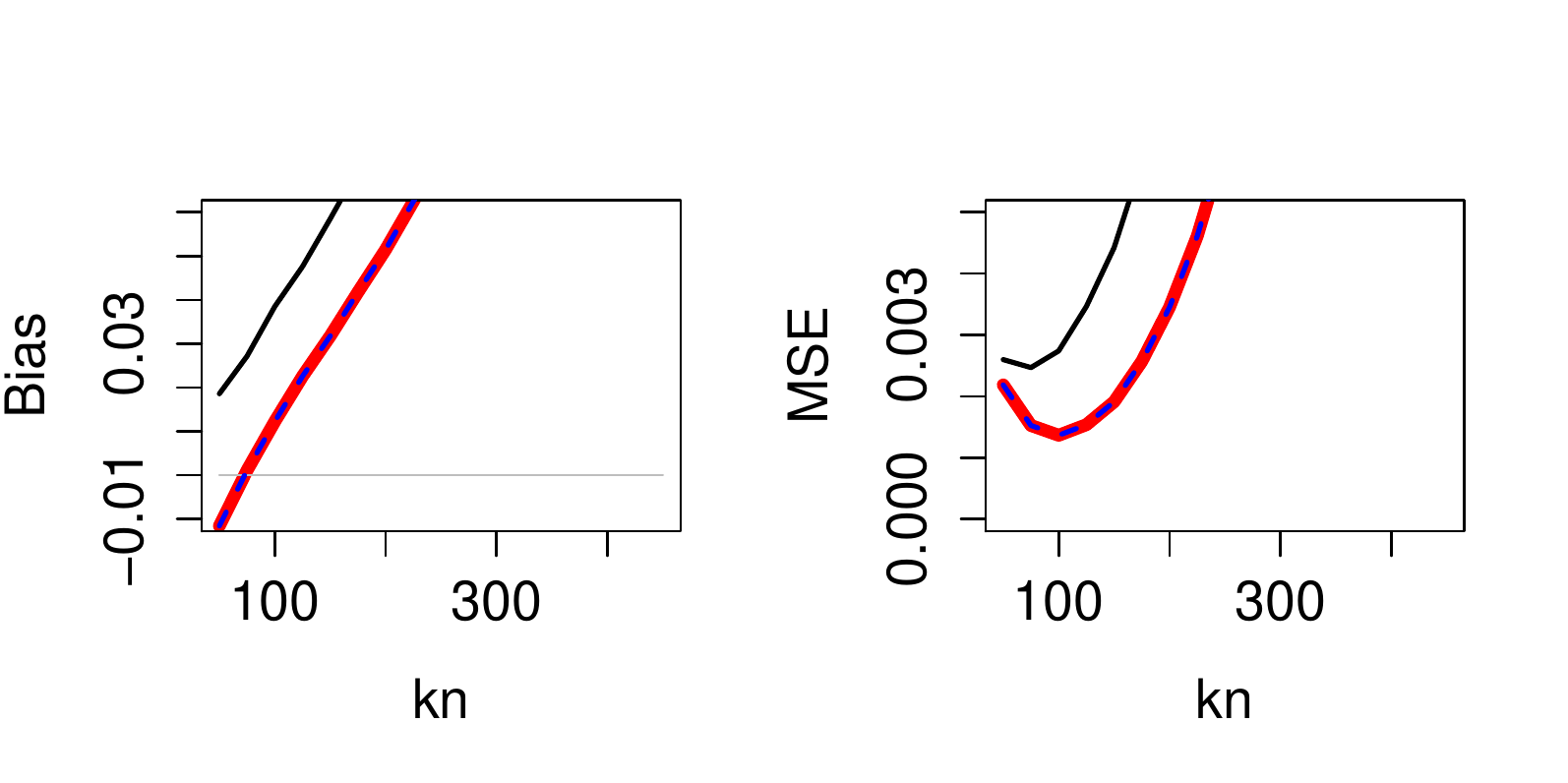}}
\hspace{.1in}
\subfigure[Case (e) but with $c_1=0.9$]{
\label{CadreFrechet3-09}
\includegraphics[height=3.5cm,width=.46\textwidth]{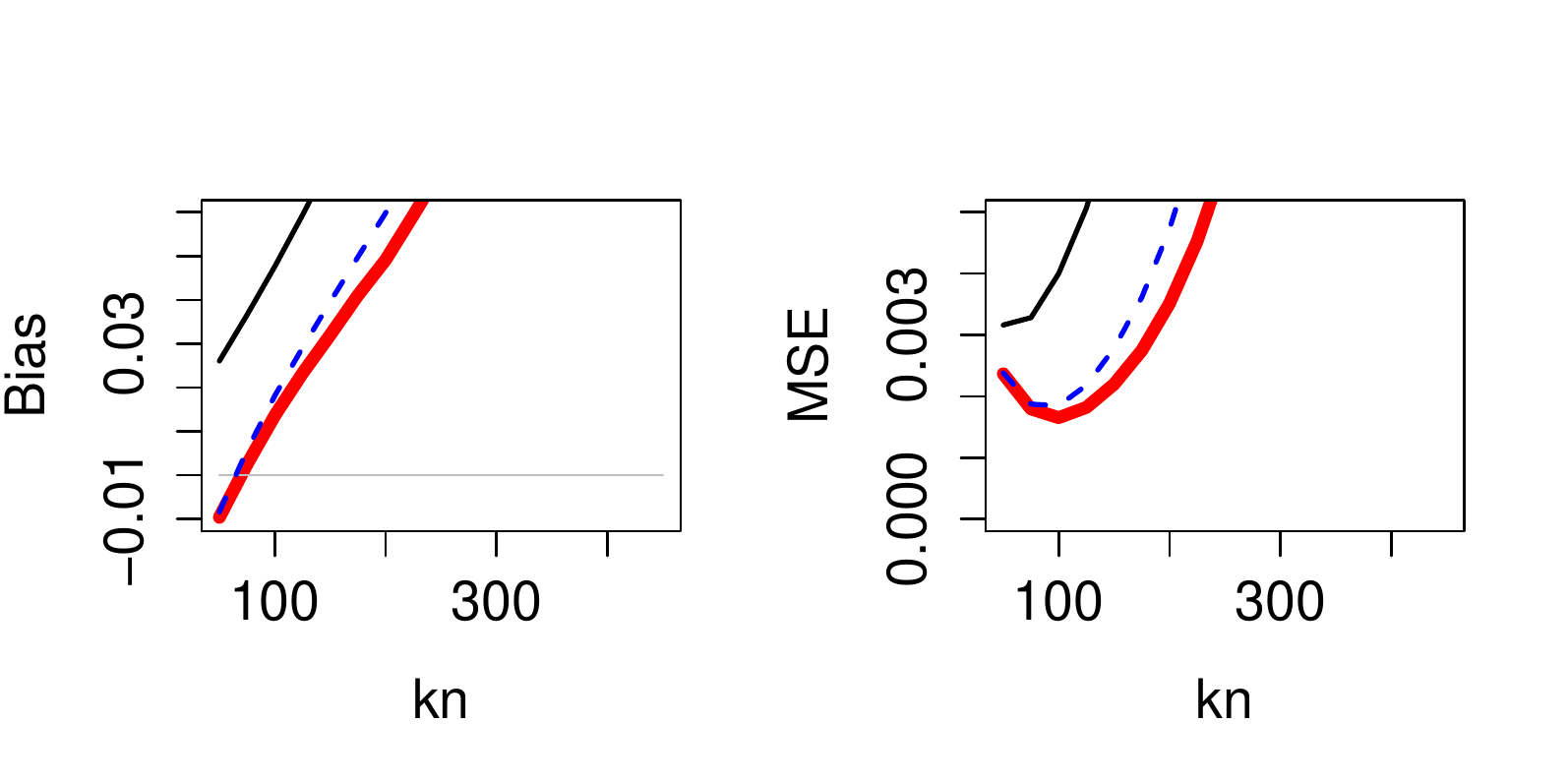}}
\hspace{.1in}
\caption{Comparison of bias  and MSE (respectively left and right in each subfigure) for $\gamchapnk$ (plain thick), $\hat\gamma_{1}^{(BDFG)}$ (plain thin) and $\hat\gamma_{1}^{(KM)}$ (dashed), for Fr\'echet distributed $X$ and $C$.}
\label{figgamchap-1}
\end{figure}

\begin{figure}[hbtp]
\centering
\subfigure[Burr case,\; $\gamma_1=0.1,\, \gamma_2=0.25,\, \gamma_C=0.3$, and $c_1=1$]{
\label{CadreBurr1-1}
\includegraphics[height=3.5cm,width=.46\textwidth]{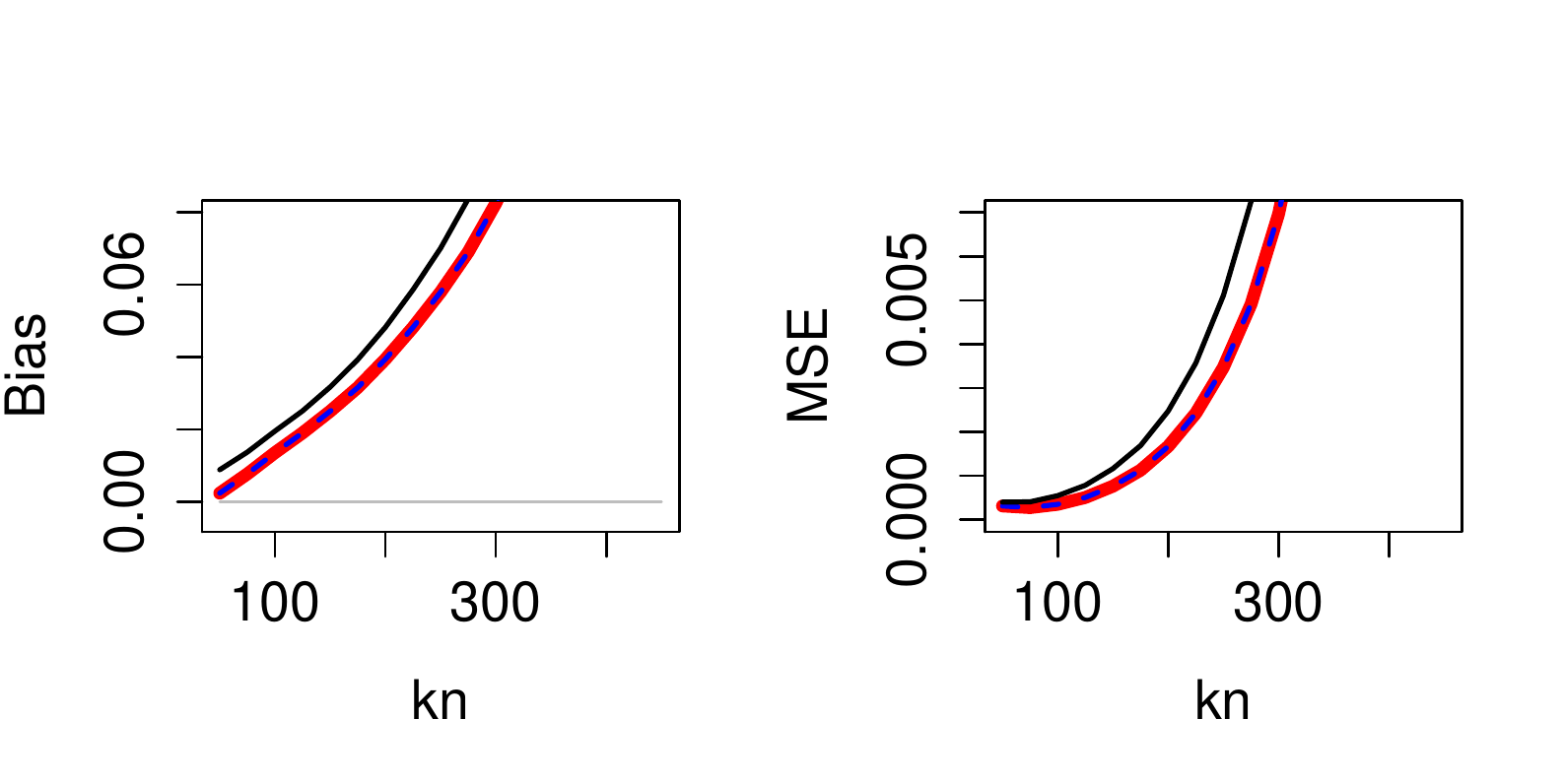}}
\hspace{.1in}
\subfigure[Case (a) but with $c_1=0.9$]{
\label{CadreBurr1-09}
\includegraphics[height=3.5cm,width=.46\textwidth]{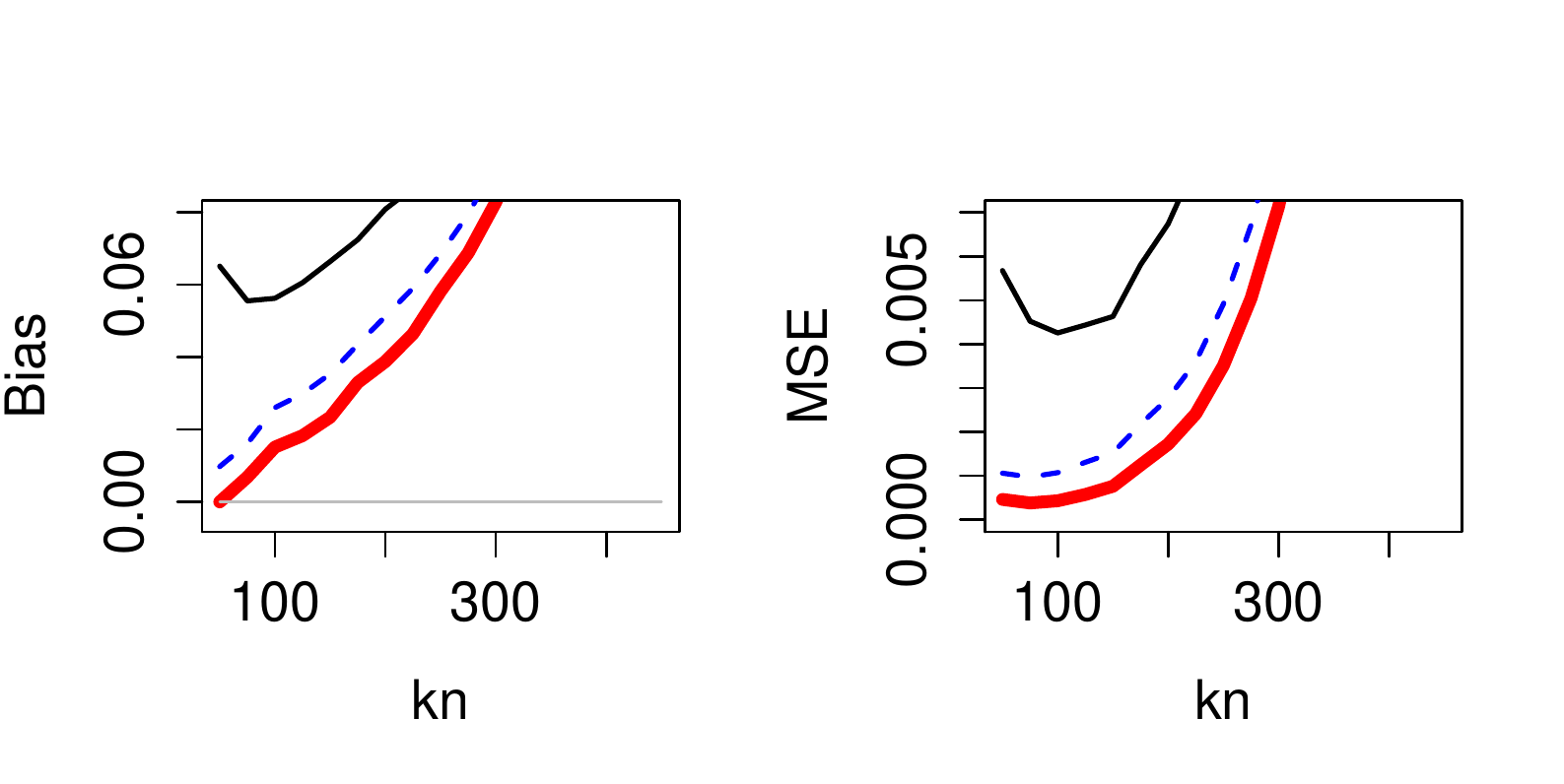}}
\hspace{.1in}
\subfigure[Burr case,\; $\gamma_1=0.1,\, \gamma_2=0.25,\, \gamma_C=0.2$, and  $c_1=1$]{
\label{CadreBurr2-1}
\includegraphics[height=3.5cm,width=.46\textwidth]{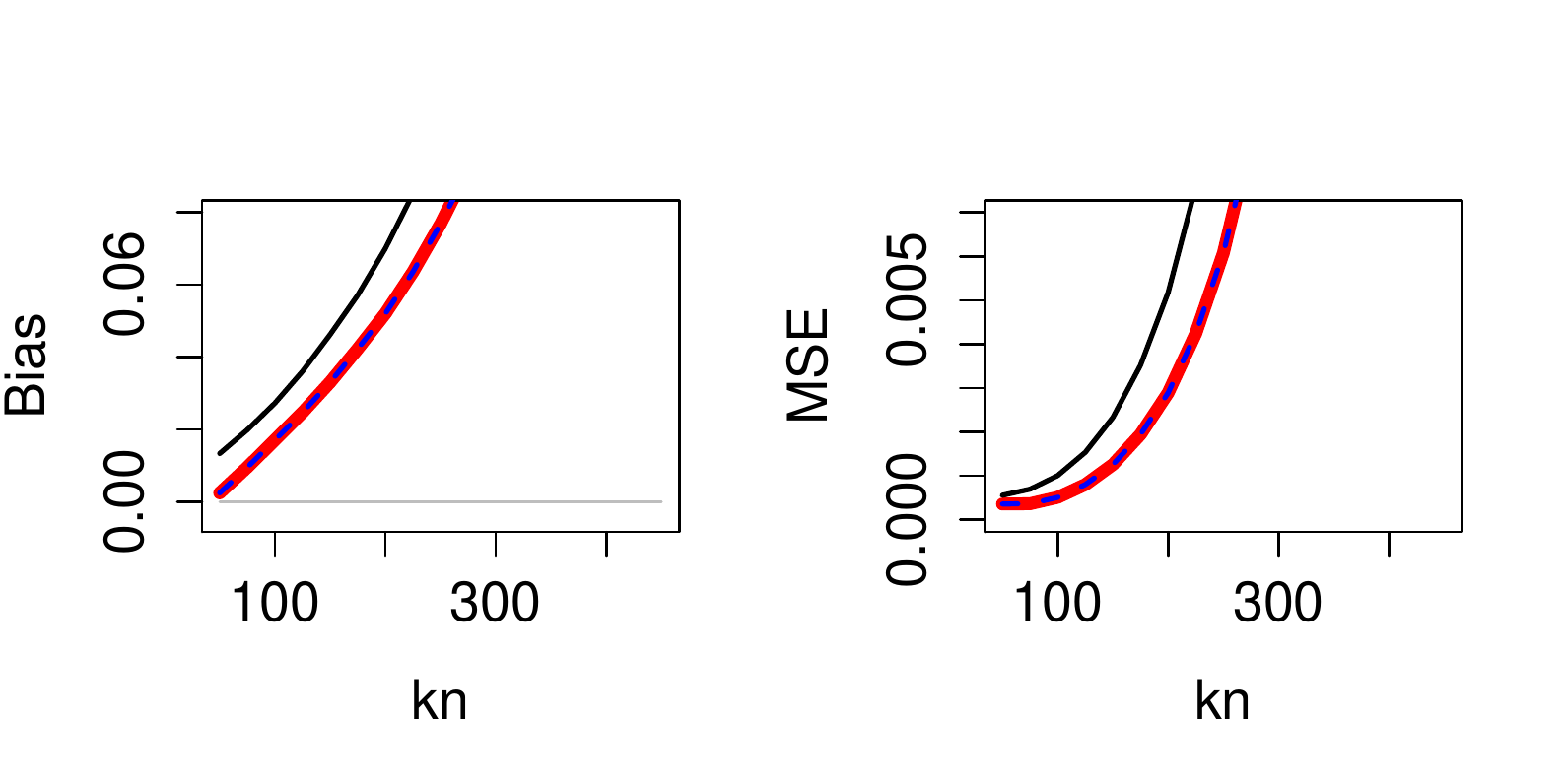}}
\hspace{.1in}
\subfigure[Case (c) but with $c_1=0.9$]{
\label{CadreBurr2-09}
\includegraphics[height=3.5cm,width=.46\textwidth]{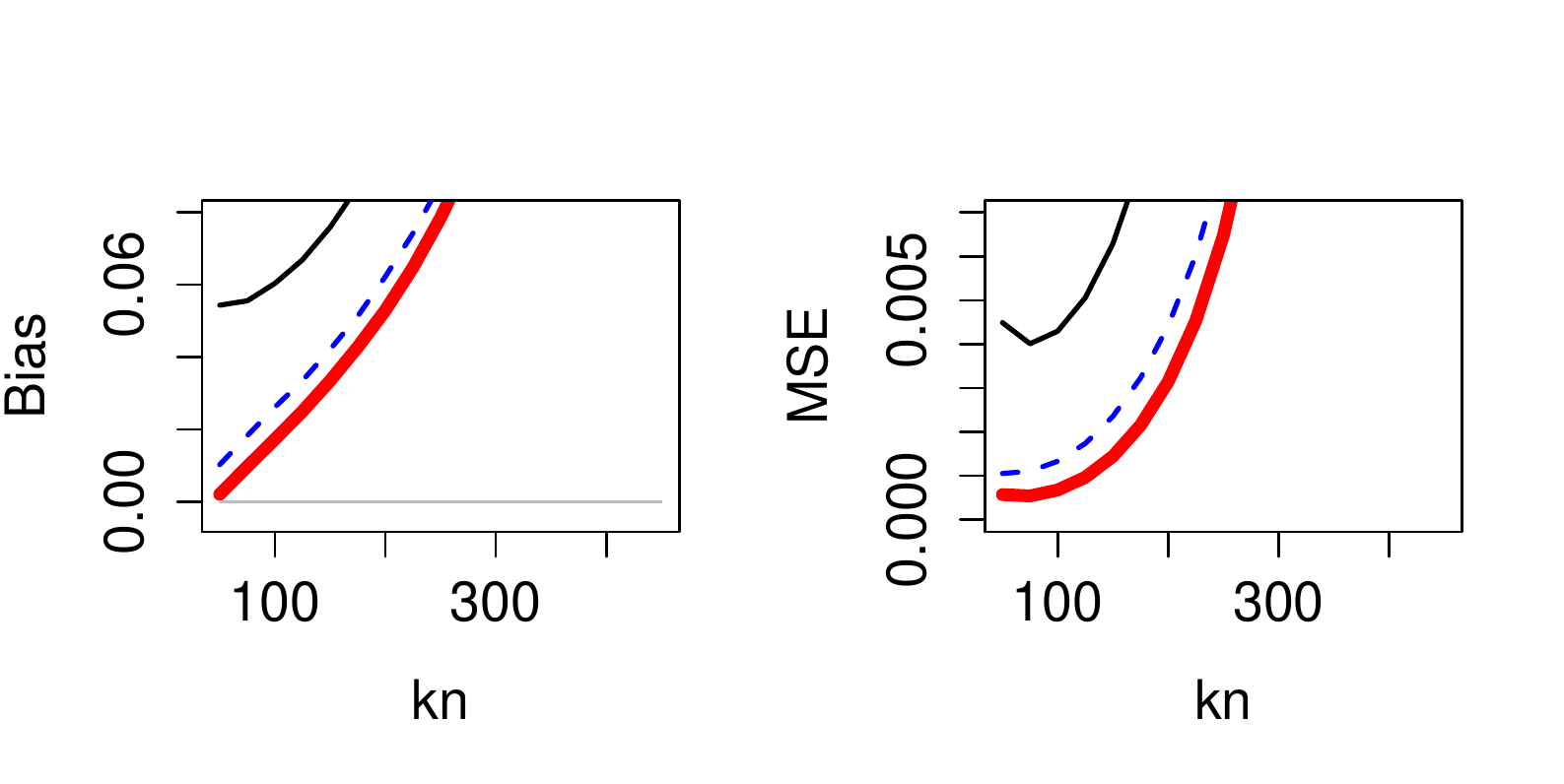}}
\hspace{.1in}
\subfigure[Burr case,\, $\gamma_1=0.1, \gamma_2=0.25, \gamma_C=0.45$, and $c_1=1$]{
\label{CadreBurr3-1}
\includegraphics[height=3.5cm,width=.46\textwidth]{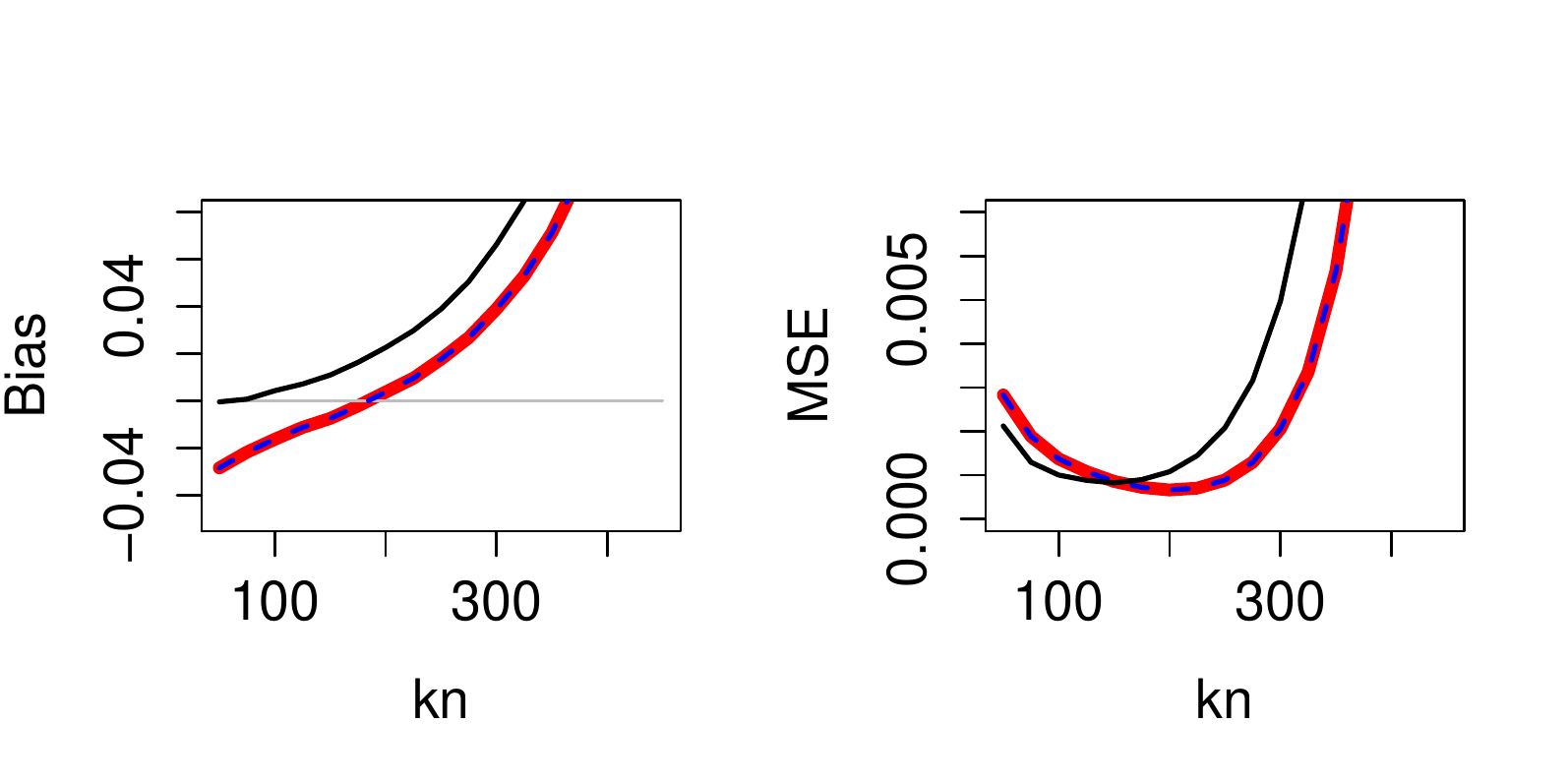}}
\hspace{.1in}
\subfigure[Case (e) but with $c_1=0.9$]{
\label{CadreBurr3-09}
\includegraphics[height=3.5cm,width=.46\textwidth]{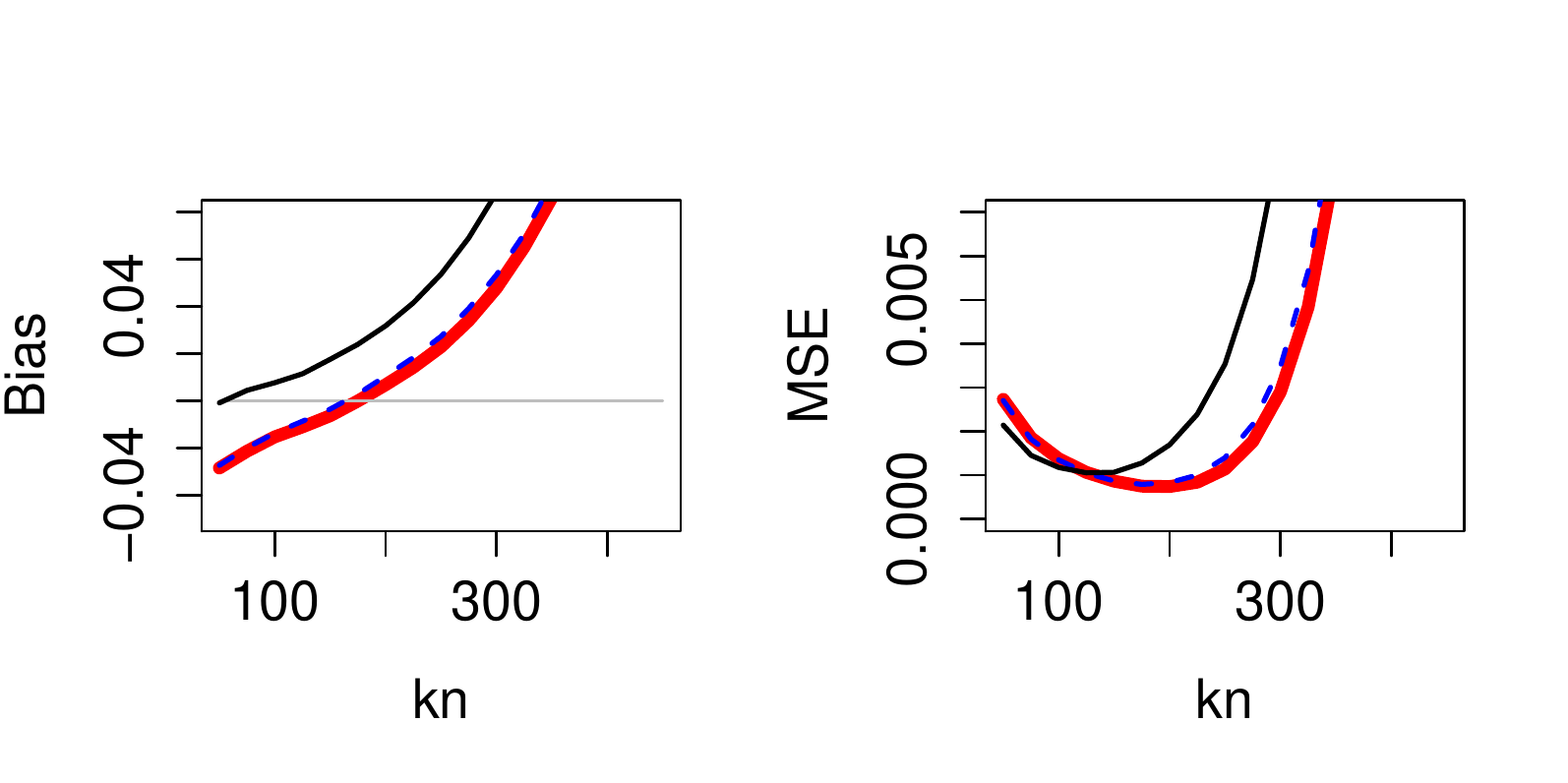}}
\hspace{.1in}
\caption{Comparison of bias  and MSE (respectively left and right in each subfigure) for $\gamchapnk$ (plain thick), $\hat\gamma_{1}^{(BDFG)}$ (plain thin) and $\hat\gamma_{1}^{(KM)}$ (dashed), for Burr distributed $X$ and $C$.  }
\label{figgamchap-2}
\end{figure}

Another aim is to illustrate the impact of dependency between the causes, when estimating the tail. The starting point is that,  if cause $2$ could be considered independent of cause $1$, then we could (and would)  include it in the censoring mechanism and we would be in the simple random censoring setting, without competing risk.  In this case, it would be possible to estimate $\gamma_1$ by one of the following  two estimators, the first  one being proposed in \cite{Beirlant07} (a Hill estimator weighted with a constant weight), and  the second one in \cite{WormsWorms14} (a Hill estimator weighted with varying Kaplan-Meier weights): 
\begin{eqnarray} \label{vieuxestim1}
  \hat\gamma_{1}^{(BDFG)} & = & 
  \frac{1}{k_n} \sum_{i=1}^{k_n} \frac{1}{\hat p_1 } \log(Z_{(n-i+1)}/Z_{(n-k_n)})  
  \zdeux\\
\hat\gamma_{1}^{(KM)} & = & 
  \frac 1 {n\Fbar_{n,b}(Z_{(n-k_n)})} \sum_{i=1}^{k_n} \frac{\delta_{(n-i+1)}\bI_{\cause_{(n-i+1)}=1}}{\Gbar_{n,b}(Z_{(n-i,n)})} \ \log(Z_{(n-i+1)}/Z_{(n-k_n)}) ,
  \label{vieuxestim2}
\end{eqnarray}  

\noindent where,  in  Equation $(\ref{vieuxestim1})$,  $\hat p_1=\frac{1}{k_n} \sum_{i=1}^{k_n} \delta_{(n-i+1)}\bI_{\cause_{(n-i+1)}=1}$, and in  Equation $(\ref{vieuxestim2})$,  the Kaplan Meier estimators $\Fbar_{n,b}$ and $\Gbar_{n,b}$ are based on the $\tilde\delta_i=\delta_i\bI_{\cause_i=1}$.  These two estimators consider the uncensored lifetimes associated to cause 2 as independent censoring times.  Comparing our new estimator with these latter two estimators, when $c_1<1$, will empirically prove that considering  cause $2$ as a competing risk independent of cause $1$ has a great (negative) impact on the estimation of $\gamma_1$. Note that when $c_1=1$, the new estimator  $\hat{\gamma}_1$ and  $\hat\gamma_{1}^{(KM)} $ are exactly the same (therefore the thick and dashed lines in sub-figures (a), (c) and (e) of Figures \ref{figgamchap-1} and \ref{figgamchap-2} are overlapping, identical).
\zdeux

We address these two aims for each set-up (Fr\'echet, or Burr), by generating $2000$ datasets of size $500$, with  three configurations of the triplet $(\gamma_1,\gamma_2,\gamma_C)$ :   
%
 $(0.1,0.25,0.3)$ ($\gamma_1<\gamma_2$, moderate censoring $\gamma_C>\gamma_F$),   $(0.1,0.25,0.2)$ ($\gamma_1<\gamma_2$, heavy censoring $\gamma_C<\gamma_F$), or $(0.25,0.1,0.45)$ ($\gamma_1>\gamma_2$, moderate censoring $\gamma_C<\gamma_F$).  Median bias and mean squared error (MSE) of the  different estimators are plotted against different values of $k_n$, the number of excesses used. When Burr distributions are simulated, the parameter $\beta$ is taken equal to $1$, and the parameters $(\tau_1,\tau_2,\tau_C)$ are taken equal to $(12,6,5)$ in configurations 1 and 2, and to $(6,12,5)$ in configuration 3. 
\zun

Figure \ref{figinfluencec}  illustrates  the behaviour  of our  estimator when $c_1$ varies. In terms of bias and MSE, we can see that the first configuration is a little better than the second one, which is itself much better than the third one. We observed this phenomenon  in many other cases, not reported here : our estimator behaves best when it is the smallest parameter $\gamma_k$ which is estimated, and when the censoring is not too strong. Our simulations also show that the quality of our estimator (especially in terms of the MSE) diminishes with $c_1$.\zun

Figures \ref{figgamchap-1} and \ref{figgamchap-2} present the comparison between  our new estimator and  the ones described in (\ref{vieuxestim1}) and (\ref{vieuxestim2}).  A general conclusion (confirmed by other simulations not reported here) is that $\hat\gamma_{1}^{(BDFG)}$ and $\hat\gamma_{1}^{(KM)}$ behave worse in most cases, even for a value of $c_1$ of $0.9$, which is only a slight modification of the situation without competing risk ($c_1=1$). Therefore, a contamination of the cause $1$ distribution by another cause rapidly yield inadequate estimations of $\gamma_1$ if dependency between causes is ignored ; this conclusion is true for both $\hat\gamma_{1}^{(BDFG)}$ and $\hat\gamma_{1}^{(KM)}$, but to a greater extent for $\hat\gamma_{1}^{(BDFG)}$ .    In the third configuration $(\gamma_1,\gamma_2,\gamma_C)=(0.25,0.1,0.45)$, the improvement provided by $\hat{\gamma}_1$ (with respect to $\hat\gamma_{1}^{(KM)}$)  becomes notable when $c_1$ drops below $0.7$.

\section{Conclusion} \label{conclusion}

In this paper,  we consider heavy tailed lifetime data subject to random censoring and competing risks, and use the Aalen-Johansen estimator of the cumulative incidence function to construct  an estimator for the extreme value index associated to the main cause of interest. To the best of our knowledge, this is the first estimator proposed in this context. Its asymptotic normality is proved and a small simulation study  exhibiting its finite-sample performance shows that accounting for the dependency of the different causes is important, but that the bias can be particularly high. Estimating second order  tail parameters would then be interesting  in order to reduce this bias. A first step  towards this aim  could be to study the following moments
\[
M_n^{(\alpha)} = \frac 1 {n \Fbarnktn} \sum_{i=1}^n \frac{\log^{\alpha}(Z_i/t_n)  }{ \GbarnZim }  \indksik \indZitn,
 \] 
which  asymptotic behaviour  can be derived  following the same lines as  in the proof of Theorem \ref{TLCGammachap}.

\section{Proofs}  \label{sec-proofs}

This section is essentially devoted to the proof of the main Theorem \ref{TLCGammachap}. Some  hints about the proof of the consistency result contained in Proposition \ref{consistance} are given in Subsection \ref{PreuveConsistance}, and Corollary \ref{coroquantiles} is proved in Subsection \ref{PreuveCoroQuantiles}. 
\zdeux

We  adopt  a  strategy developed by Stute in \cite{Stute95} in order to prove his Theorem 1.1, a well-known result which states that a Kaplan-Meier integral of the form $\int \phi \, dF_n$ can be approximated by a sum of independent terms. This idea is used in \cite{Suzukawa02} in the context of competing risks. 
 We  thus intend to approximate $\gamchapnk$ by the integral $\gamtildenk=\int \phi_n \, d\Fnk$ of some deterministic function $\phi_n$,  with respect to the Aalen-Johansen estimator, 
and approximate this integral by the mean $\gamtchenk$ of independent variables $U_{i,n}$ (defined a few lines below). The passage from $\gamchapnk$ to $\gamtildenk$ (which amounts to replacing $\Fbarnktn$ by $\Fbarktn$ in the denominator of $\gamchapnk$) will imply an additional  sum of independent variables $V_{i,n}$, which will participate to the asymptotic variance of our estimator. 
\zdeux

However, a major difference with  \cite{Stute95} or \cite{Suzukawa02}  is that the function we integrate here,  $\phi_n(u)= \frac{1}{\Fbark(t_n)} \log (u/t_n) \bI_{u>t_n}$, is  not only  an unbounded function, depending on $n$, but it also has a  "sliding" support $[t_n,+\infty[$,  which is therefore always close to the endpoint $+\infty$ of the distribution $H$. In \cite{Stute95},  a crucial point of the proof consists in temporarily considering  that the integrated function $\phi$ has a support which is bounded away from the endpoint  of $H$ (condition (2.3) there). Considering the kind of function $\phi_n$ we have to deal with here, we cannot follow the same strategy :
dealing with the remainder terms will thus be a particularly challenging part of our work. Finally note that, in order to deal with the ratio $\Fbarnktn/\Fbarktn$ (and somehow approximate $\gamchapnk$ by $\gamtildenk$) we will  have to consider simultaneously   integrals (with respect to $\Fnk$) of $\phi_n$ and of another function $g_n$, defined below, which basically shares the same flaws as $\phi_n$.
\zdeux

Let us first recall or define   the following objects :
\begin{eqnarray*}
\widehat\phi_n(u) & = & \unsurFbarnktn \log \left( \frac u {t_n} \right)  \bI_{u>t_n} \\
\phi_n(u) & = & \unsurFbarktn \log \left( \frac u {t_n} \right) \bI_{u>t_n} \\
\gammank & = & \int \phi_n(u) d\Fk(u)  \ \tqdninf \ \gamma_k \\
\gamtildenk & = & \int \phi_n(u) d\Fnk(u) \\
\gamchapnk & = & \int \widehat\phi_n(u) d\Fnk(u).
\end{eqnarray*}
We thus have $\gamchapnk = \Delta_n^{-1} \gamtildenk$, where
$$
 \Delta_n \; = \; \Fbarnktn / \Fbarktn \ = \ \int g_n(u) d\Fnk(u) \makebox[1.7cm][c]{and} 
g_n(u) \ = \ \unsurFbarktn \bI_{u>t_n},
$$
and we now introduce the following new quantities, related to the Stute-like decomposition of $\gamtildenk$ and $\Delta_n$  : 
\begin{eqnarray*}
\Uin{1} & = & \frac{\phi_n(Z_i)}{\Gbar(Z_i)} \delta_i \indcausek  \makebox[1.5cm][c]{and} 
\Vin{1} \ = \ \frac{g_n(Z_i)}{\Gbar(Z_i)} \delta_i \indcausek   \\ 
\Uin{2} & = & \frac{1-\delta_i}{\Hbar(Z_i)} \psi(\phi_n,Z_i)  \makebox[1.5cm][c]{and} 
\Vin{2} \ = \ \frac{1-\delta_i}{\Hbar(Z_i)} \psi(g_n,Z_i) \\
\Uin{3} & = & \int_0^{Z_i} \psi(\phi_n,u)\, dC(u) \makebox[1.5cm][c]{and} 
\Vin{3} \ = \ \int_0^{Z_i} \psi(g_n,u)\, dC(u)  \\
U_{i,n} & = & \Uin{1} + \Uin{2} - \Uin{3} \makebox[1.5cm][c]{and} V_{i,n} \ = \ \Vin{1} + \Vin{2} - \Vin{3} \end{eqnarray*}
where, for any function $f:\bR_+\ra\bR$, we note (for any given $z\geq 0$)
\[ 
\psi(f,z) = \int_z^{+\infty} f(t) d\Fk(t) \makebox[1.6cm][c]{ and } C(z)= \int_0^z \frac{dG(t)}{\Hbar(t) \Gbar(t)}.
\]  
This enables us to finally define the important objects 
\begeq{defgamtchetcheDeltachap}
\gamtchenk \ = \ \frac 1 n \sum_{i=1}^n U_{i,n}  \makebox[1.5cm][c]{and} \widehat\Delta_n \ = \ \frac 1 n \sum_{i=1}^n V_{i,n}
\fineq
which are the triangular sums of independent terms which will respectively approximate $\gamtildenk$ and $\Delta_n$. At the beginning of section \ref{sec-TLCpourZn}, it will be proved that $\bE(\Uin{1})=\gammank$ and $\bE(\Vin{1})=1$, while $\bE(\Uin{2})=\bE(\Uin{3})$ and $\bE(\Vin{2})=\bE(\Vin{3})$, yielding $\bE(\gamtchenk)=\gammank$ and $\bE(\widehat\Delta_n)=1$ ; the terms $\Uin{2}$, $\Uin{3}$, $\Vin{2}$ and $\Vin{3}$ only participate to the variance component of the estimator.   The relation between all these quantities is made clearer in the following Lemma : 

\begin{lem} \label{lemmedecompgammachap}
We have  
\begeq{decompgamchap}
\sqrt{v_n} (\gamchapnk - \gamma_k) \ = \ \Delta_n^{-1} \left( \ Z_n \, + \, \sqrt{v_n} R_n \, + \, \sqrt{v_n} (\gammank - \gamma_k) \    \right) 
\fineq
where 
\[
  Z_n \ = \ \sqrt{v_n}\left( \, (\gamtchenk-\gammank) - \gamma_k(\widehat\Delta_n - 1) \, \right)
  \makebox[1.5cm][c]{and}
  R_n \ = \ (\gamtildenk-\gamtchenk) - \gamma_k(\Delta_n - \widehat\Delta_n)
\]
\end{lem}

The proof of Lemma \ref{lemmedecompgammachap}  is simple :  
\begin{eqnarray*}
\sqrt{v_n} (\gamchapnk - \gamma_k)  & = &  \Delta_n^{-1} \sqrt{v_n}  \left( \gamtildenk - \Delta_n\gamma_k \right) \\
& = & \Delta_n^{-1} \sqrt{v_n}  \left\{  (\gamtchenk - \gammank) + \gamma_k(1-\Delta_n) + (\gamtildenk - \gamtchenk) + (\gammank-\gamma_k) \right\}
\end{eqnarray*}
which leads to the desired relation (\ref{decompgamchap}).  
\ztrois
 
The main theorem thus becomes an immediate consequence of the following four results, the second one being the most difficult to establish.  

\begin{prop}
\label{TLCpourZn}
Under condition  $(\ref{Ordre1})$  and assuming that 
\begeq{condvntn1}
 v_n \tqdninf + \infty ,
\fineq
if $\gamma_k < \gamma_C$, then 
 $$
  Z_n \ \stackrel{d}{\longrightarrow} \  {\cal N}(0,\sigma^2) \hs{0.5cm} \mbox{as $n\tinf$}
 $$ 
 where $\sigma^2$ is defined in the statement of Theorem \ref{TLCGammachap}. 
\end{prop}

\begin{prop}   
\label{PropGlobaleReste}
 Under conditions  $(\ref{Ordre1})$  and $(\ref{condvntn})$ , if $\gamma_k < \gamma_C$,
then 
\begeq{decomposition}
  \gamtildenk = \gamtchenk + R_{n,\phi}  \makebox[1.5cm][c]{and}   \Delta_n = \widehat\Delta_n + R_{n,g} 
\fineq
where $R_{n,\phi}$, $R_{n,g}$ (and consequently $R_n$ too) are $o_{\bP}(v_n^{-1/2})$. 
\end{prop}

\begin{cor}\label{coroDeltan}
 Under the conditions of Proposition \ref{PropGlobaleReste},  $\racinevn(\Delta_n-1)$ converges in distribution to ${\cal N}(0,1/(1-r))$ where $r=\gamma_k/\gamma_C\in ]0,1[$.  
\end{cor}

\begin{lem} \label{lemmebiais}
Under conditions $(\ref{Ordre1})$, $(\ref{Ordre2})$  and $\racinevn g(t_n)\ra \lambda\geq 0$, the bias term  $\racinevn (\gammank - \gamma_k)$  in  (\ref{decompgamchap})  converges to $\lambda  m$  as $n\tinf$,
where $m$ is defined in Theorem \ref{TLCGammachap}.
\end{lem}

Propositions \ref{TLCpourZn} and \ref{PropGlobaleReste} will be proved in Sections \ref{sec-TLCpourZn} and \ref{sec-Reste} respectively, sometimes with the help of other results stated and established in the Appendix. The proofs of  Corollary \ref{coroDeltan} and Lemma \ref{lemmebiais}   are short, we state them below. 
\zdeux 

Concerning Corollary \ref{coroDeltan}, once the proof of Proposition \ref{TLCpourZn} has been gone through, it will become clear to the reader that $\racinevn(\hat{\Delta}_n-1)$ converges in distribution to the centred gaussian distribution of variance $1/(1-r)$, because $\racinevn(\hat{\Delta}_n-1)=\sum_{i=1}^n \tilde W_{i,n}$ where $\tilde W_{i,n}=\frac{\racinevn}{n} \tildeVin = \frac{\racinevn}{n} (V_{i,n} -1)$ are centred, and $\Var(\tilde W_{i,n})=\frac 1 n \frac 1 {1-r} + o(1/n)$ (this is proved similarly as (\ref{expressionVarW1nInit}) and (\ref{expressionVarW1n})). Since Proposition \ref{PropGlobaleReste}   states that $\Delta_n=\hat{\Delta}_n +o_{\bP}(v_n^{1/2})$, the same central limit theorem holds for $\Delta_n$ and the corollary is proved. 
\zdeux 

Concerning now  Lemma \ref{lemmebiais},  remind that $\gammank = \int\phi_n(u)d\Fnk(u)$. An integration by parts and the fact that $\Fbark(x)=x^{-1/\gamma_k}l_k(x)$ yield 
$$
 \racinevn (\gammank-\gamma_k) = \racinevn\int_1^{+\infty} y^{-1/\gamma_k-1}\left( \frac{l_k(yt_n)}{l_k(t_n)} - 1\right) dy,
$$
and, using assumption $(\ref{Ordre2})$   and Proposition 3.1 in \cite{HaanFerreira06}, we can write 
$$
 \int_1^{+\infty} y^{-1/\gamma_k-1}\left( \frac{l_k(yt_n)}{l_k(t_n)} - 1\right) dy 
 = g(t_n) \int_1^{+\infty} y^{-1/\gamma_k-1} h_{\rho_k}(y) dy + o(g(t_n)).
$$
The result then follows from assumption $\racinevn g(t_n)\ra \lambda\geq 0$  and the fact that $\int_1^{+\infty} y^{-1/\gamma_k-1} h_{\rho_k}(y) dy = m$. 
 \ztrois

In the rest of the paper, we will very often handle the well-known sub-distributions functions $\HO$ and $\Hunk$ defined, for  all $t \geq 0$, by 
\[
\HO(t)= \bP(Z \leq t, \delta=0) \makebox[1.6cm][c]{ and }  \Hunk(t)= \bP(Z \leq t, \xi=k).
\]
Note that we have
\[
d\HO= \Fbar dG  \makebox[1.6cm][c]{ and }  d\Hunk= \Gbar d\Fk.
\]

\subsection{Proof of Proposition \ref{TLCpourZn}} \label{sec-TLCpourZn}

We first write
$$
 Z_n \ = \ \sum_{i=1}^n W_{i,n}    \makebox[1.5cm][c]{where}  
 W_{i,n} \ = \ \frac{\racinevn}{n} \left(\tildeUin - \gamma_k \tildeVin \right) 
   \makebox[1.5cm][c]{and} \left\{ \begar{l} \tildeUin = U_{i,n} - \gammank \\ \tildeVin = V_{i,n} - 1 \finar\right. 
$$
where $W_{i,n}$, $\tildeUin$ and $\tildeVin$ are centred, because the random variables $U_{i,n}$ and $V_{i,n}$ have expectations respectively  equal to $\gammank$ and $1$. Indeed, we have
$$
 \mbox{$
  \bE\left(\Uin{1}\right) = 
  \bE\left( \frac{\phi_n(Z_i)}{\Gbar(Z_i)} \delta_i\indcausek \right) 
  = \int \frac{\phi_n(u)}{\Gbar(u)} d\Hunk(u) = \int \phi_n(u) d\Fk(u) =\gammank
$}
$$
and
$$ 
 \textstyle
  \bE\left(\Uin{2}\right) \; = \;  
  \bE\left( \frac{1-\delta_i}{\Hbar(Z_i)} \psi(\phi_n,Z_i) \right) \; = \;   
  \int \int \frac 1{\Hbar(u)} \bI_{t>u} \phi_n(t) \,d\Fk(t)\,d\HO(u) \; = \;   
  \int \int \frac 1{\Gbar(u)} \bI_{t>u} \phi_n(t) \,d\Fk(t)\,dG(u) 
$$
as well as
$$
 \textstyle
  \bE\left(\Uin{3}\right) \; = \;  
  \bE\left( \int_0^{Z_i} \psi(\phi_n,u) \, dC(u) \right) \; = \;   
  \int \int \int \bI_{z>u} \bI_{t>u} \phi_n(t) \,dH(z)\,d\Fk(t)\,\frac{dG(u)}{\Gbar(u)\Hbar(u)}
  \; = \;   
  \bE\left(\Uin{2}\right).
$$
The proof for $\bE(V_{i,n})=1$ is similar. 
\zun\\
We will now prove the asymptotic normality of $Z_n$ by using the Lyapunov criteria. 
\zdeux

\begin{lem} \label{lem-varZn}
Under the conditions  $(\ref{Ordre1})$ and $(\ref{condvntn1})$ , if $\gamma_k<\gamma_C$ :
\begenum 
\item[$(i)$] we have 
 \begeq{expressionVarW1nInit}
   \Var(\tildeUunn-\gamma_k\tildeVunn) \ = \ 
   \bE(\Uunn^2)+\gamma_k^2\bE(\Vunn^2)-2\gamma_k\bE(\Uunn\Vunn)+o(1) 
 \fineq
\item[$(ii)$] we have 
 \begin{eqnarray}
  \bE(\Uunn^2) & = & \int \frac{\phi_n^2(u)}{\Gbar(u)} d\Fk(u) \ - \ \int (\psi(\phi_n,u))^2\, dC(u) 
  \label{formuleU1ncarre}  \\
  \bE(\Vunn^2) & = & \int \frac{g_n^2(u)}{\Gbar(u)} d\Fk(u) \ - \ \int (\psi(g_n,u))^2\, dC(u) 
  \label{formuleV1ncarre} \\
  \bE(\Uunn\Vunn) & = & \int \frac{\phi_n(u)g_n(u)}{\Gbar(u)} d\Fk(u) \ - \ 
        \int \psi(\phi_n,u)\psi(g_n,u)\, dC(u) 
  \label{formuleU1nV1n} 
 \end{eqnarray}
\item[$(iii)$] we have, noting $r=\gamma_k/\gamma_C$ (which belongs to $]0,1[$ under our conditions) as well as  $p=\gamma/\gamma_F=\gamma_C/(\gamma_F+\gamma_C)\in ]0,1[$,
 \begeq{expressionVarW1n}
   \Var(\tildeUunn-\gamma_k\tildeVunn) \ = \ 
   \frac 1 {\Fbarnktn\Gbar(t_n)} \left( \frac{ \gamma_k^2(1+r^2)} {(1-r)^3} + o(1) \right)
   \; - \; 
   \frac{1-p}{\Hbar(t_n)}  \left( \frac{2\gamma_k^3}{\gamma(2-\gamma_k/\gamma)^3} + o(1) \right) 
   + o(1)
 \fineq
 \finenum
\end{lem}

\begin{lem} \label{lem-Lyapunov}
Under the conditions  $(\ref{Ordre1})$ and $(\ref{condvntn1})$ , if $\gamma_k<\gamma_C$, then 
\begct
$\sum_{i=1}^n \bE \left|W_{i,n} \right|^{2+ \delta}  \rightarrow 0$, as $n$ tends to infinity,  for some $\delta >0$.
\finct
\end{lem} 
\zdeux

We can then immediately prove Proposition \ref{TLCpourZn}. Indeed, since $Z_n=\sum_{i=1}^n W_{i,n}$,  Lemma \ref{lem-varZn} yields
\[
\Var(Z_n) = n\Var(W_{1,n}) = \frac {v_n} n \Var(\tildeUunn-\gamma_k\tildeVunn) .
\]
which, since  $v_n=n\Fbarnktn\Gbar(t_n)$, becomes 
\[
\Var(Z_n) \; = \; \gamma_k^2(1+r^2)(1-r)^{-3}  \; - \left( 2(1-p)\gamma_k^3\gamma^{-1}(2-\gamma_k/\gamma)^{-3} \right) \left(\Fbarktn/\Fbar(t_n)\right) \; + \; o(1) .
\]
Therefore, depending on the limit $c$ of the ratio $\Fbarktn/\Fbar(t_n)$ when $n\tinf$ (for instance, it converges to $0$ when  $\gamma_k<\gamma_F$),  it is simple to check that the variance of $Z_n$ converges to the value $\sigma^2$ described in the statement of Theorem \ref{TLCGammachap}. Thanks to Lemma \ref{lem-Lyapunov}, the Lyapunov CLT applies and Proposition \ref{TLCpourZn} is proved.  
\ztrois

  The two subsections \ref{subsec-varZn} and \ref{subsec-Lyapunov} are  now respectively devoted to the proofs of Lemmas \ref{lem-varZn} and \ref{lem-Lyapunov}.

\subsubsection{Proof of Lemma \ref{lem-varZn}} \label{subsec-varZn} \ztrois

Part $(i)$ of the lemma is straightforward : since $\tildeUunn$ and $\tildeVunn$ are centred, we have indeed
\begin{eqnarray*}
 \textstyle 
  \Var(\tildeUunn-\gamma_k\tildeVunn) & = &  \textstyle 
  \bE((\tildeUunn-\gamma_k\tildeVunn)^2) \ = \ \; 
  \bE(\tildeUunn^2) +\gamma_k^2\bE(\tildeVunn^2) - 2\gamma_k\bE(\tildeUunn\tildeVunn) 
    \\ & = & \textstyle 
 (\bE(\Uunn^2)-\gammank^2) \, + \, \gamma_k^2(\bE(\Vunn^2) - 1) \, - \,  2\gamma_k(\bE(\Uunn\Vunn)-\gammank) 
\end{eqnarray*}
and the result comes by using the fact that $\gammank$ converges to $\gamma_k$ as $n\tinf$.
\zun

Now we proceed to the proof of part $(ii)$, and will only prove (\ref{formuleU1ncarre}) because, by definition of $\phi_n$ and $g_n$, the proofs for (\ref{formuleV1ncarre}) and (\ref{formuleU1nV1n}) will be completely similar. First of all, we obviously have
\begeq{Uincarre}
 (\Uunn)^2= (\Uunnb{1})^2+(\Uunnb{2})^2+(\Uunnb{3})^2
                 +2\Uunnb{1}\Uunnb{2}-2\Uunnb{1}\Uunnb{3}-2\Uunnb{2}\Uunnb{3}
\fineq
The first term in the right-hand side of (\ref{formuleU1ncarre}) is equal to $\bE((\Uunnb{1})^2)$, and the second one (without the minus sign) is equal to $\bE((\Uunnb{2})^2)$ and to $\bE(\Uunnb{1}\Uunnb{3})$  because
\[
 \bE((\Uunnb{2})^2) = \bE\left( \frac{1-\delta_1}{\Hbar^2(Z_1)} (\psi(\phi_n,Z_1))^2 \right)
 = \int \frac{1}{\Hbar^2(z)} (\psi(\phi_n,z))^2 \,d\HO(z)  =  \int  (\psi(\phi_n,z))^2 \,dC(z)
\]
and
\begin{eqnarray*}
 \bE(\Uunnb{1}\Uunnb{3}) & = & \int  \frac{\phi_n(z)}{\Gbar(z)} \left( \int_0^z \psi(\phi_n,u) dC(u) \right) d\Hunk(z)
 \\ & = &  \int \psi(\phi_n,u) \left(\int_u^{\infty} \phi_n(z) d\Fk(z)\right) \,dC(u)  \ = \ \,  \int  (\psi(\phi_n,z))^2 \,dC(z)
\end{eqnarray*}
The expectation $\bE(\Uunnb{1}\Uunnb{2})$ equals $0$ because $\delta_1(1-\delta_1)$ is constantly $0$, and we are now going to prove that $\bE((\Uunnb{3})^2)=2\bE(\Uunnb{2}\Uunnb{3})$, which ends the proof of (\ref{formuleU1ncarre}) in view of (\ref{Uincarre}). Indeed, noting $h(z)=\int_0^z \psi(\phi_n,u)dC(u)$ and using the simple fact that $h(z)=h(y)+\int_y^z \psi(\phi_n,u)dC(u)$ for every $y<z$, we have
\begin{eqnarray*}
 \bE((\Uunnb{3})^2) & = & 
 \int_0^\infty \left( \int_0^z \psi(\phi_n,y) h(z) dC(y) \right) dH(z)  
 \\ & = & 
 \int_0^\infty \left( \int_0^z \psi(\phi_n,y) h(y) dC(y) \right) dH(z)  
 \; + \; 
 \int_0^\infty \left\{ \int_0^z \left( \int_y^z \psi(\phi_n,u) dC(u) \right) \psi(\phi_n,y) dC(y) \right\} dH(z)  
 \\  & = &  
 \int_0^\infty \left( \int_0^z \psi(\phi_n,y) h(y) dC(y) \right) dH(z)  
 \; + \; 
 \int_0^\infty \left\{ \int_0^z \left( \int_0^u \psi(\phi_n,y) dC(y) \right) \psi(\phi_n,u) dC(u) \right\} dH(z)  
 \\ & = &
 2 \int_0^\infty \left( \int_0^z \psi(\phi_n,y) h(y) dC(y) \right) dH(z)  
\end{eqnarray*}
and
\begin{eqnarray*}
 \bE(\Uunnb{2}\Uunnb{3}) & = & 
 \bE\left(\, \frac{1-\delta_1}{\Hbar(Z_1)} \psi(\phi_n,Z_1) \, h(Z_1)\, \right)
 \; = \; 
 \int \Hbar(y) h(y) \psi(\phi_n,y) \frac{dG(y)}{\Hbar(y)\Gbar(y)} 
 \\ & = &
 \int_0^\infty \left( \int_y^{\infty} dH(z) \right)   h(y) \psi(\phi_n,y) \, dC(y) 
 \\ & = & 
 \int_0^\infty \left( \int_0^z \psi(\phi_n,y) h(z) dC(y) \right) dH(z)  
 \; = \; \frac 1 2  \bE((\Uunnb{3})^2) 
\end{eqnarray*}
as announced. 
\zun

We can now start proving part $(iii)$ of the lemma, in which the exact nature of the function $\phi_n$ matters. First remind that functions $\Fbark$, $\Gbar$ and $C$ are regularly varying of  respective orders $-1/\gamma_k$, $-1/\gamma_C$ and $1/\gamma$ (for $C$, this is proved in Lemma \ref{lem-Cdelta} with $\delta=1$).  Let us define the constants $c_j$ and $d_j$ ($j=0,1,2$) by
\[
 c_j = \frac{j!\gamma_k^j} {(1-q)^{j+1}}   
 \makebox[1.2cm][c]{and}  
 d_j = \frac{j!\gamma_k^{j+1}} {\gamma(2-\gamma_k/\gamma)^{j+1}}. 
\]
Since $\gamma_k<\gamma_C$ was assumed, then according to Lemma \ref{lemtechInteg}  part $(ii)$ (applied first with $a+b=1/\gamma_C - 1/\gamma_k <0$  for $c_j$, and then with $a+b=-2/\gamma_k + 1/\gamma = (1/\gamma_C-1/\gamma_k) + (1/\gamma_F-1/\gamma_k) <0$ for $d_j$) , we have  
\begeq{cjdj}
\int_{t_n}^{\infty} \log^j\left(\frac u {t_n}\right) \frac{d\Fk(u)}{\Gbar(u)}  
 \; \sim \; c_j \frac{\Fbark(t_n)}{\Gbar(t_n)} 
  \makebox[1.2cm][c]{and}   \int_{t_n}^{\infty} \log^j\left(\frac u {t_n}\right) \left(\fracFbark{u}{t_n}\right)^2 \frac{dC(u)}{C(t_n)} \ \tqdninf \ d_j . 
\fineq
 Hence, by definition of $\phi_n$, $g_n$, the  first terms of  $\bE(\Uunn^2)$, $\bE(\Vunn^2)$ and $\bE(\Uunn\Vunn)$ in relations (\ref{formuleU1ncarre}), (\ref{formuleV1ncarre}) and (\ref{formuleU1nV1n}) are respectively equivalent (as $n\tinf$) to $c_2D(t_n)$, $c_0D(t_n)$ and $c_1D(t_n)$ where $D(t_n)$ denotes
\[
 D(t_n) = \frac 1 {\Fbarktn\Gbar(t_n)}.
\]
Since $c_2+\gamma_k^2 c_0-2\gamma_k c_1$ is found to be equal to $\gamma_k^2(1+r^2)/(1-r)^3$, then in view of (\ref{expressionVarW1nInit}) this proves the first term in relation (\ref{expressionVarW1n}). We now need to obtain equivalent expressions for the quantities $\int (\psi(\phi_n,u))^2\, dC(u)$,  $\int (\psi(g_n,u))^2\, dC(u) $ and $\int \psi(\phi_n,u)\psi(g_n,u)\, dC(u) $ in order to prove the second part of relation (\ref{expressionVarW1n}) and therefore finish the proof of Lemma \ref{lem-varZn}.  
\zun

For saving space, we will use temporarily the following notations :
\[
   l_n(u) = \log\left(u/t_n\right)  \mkb[0.6cm][c]{,}  R_n(u) = \Fbark(u)/\Fbark(t_n)
\]
According to the technical Lemma \ref{lemFonctionPsi} of the Appendix and, after splitting the integral into $\int_0^{+\infty}$ and $\int_{t_n}^{+\infty}$, we can write
\begin{eqnarray}
& & \int (\psi(\phi_n,u))^2\, dC(u)  \; = \;  \gamma^2_{n,k} \int_0^{t_n} \! dC(u) \nonumber 
\\ 
& & \null\hspace*{1.2cm} + \int_{t_n}^{+\infty} \! \left( l^2_n(u) R^2_n(u) + \! \gamma^2_k \left(u/t_n \right)^{-2/\gamma_k} + \! 2 \gamma_k l_n(u) R_n(u) \left(u/t_n \right)^{-1/\gamma_k} \right) \! dC(u)  + o(C(t_n)), \label{JU}
\end{eqnarray}
where $o(C(t_n))$ in $(\ref{JU})$ is due to part $(ii)$ of Lemma \ref{lemtechInteg}  and to the fact that $\epsilon_n(u)$ in Lemma \ref{lemFonctionPsi} converges  to $0$ uniformly in $u$. According to  the second part of relation $(\ref{cjdj})$,  we thus have
\begin{eqnarray}
\int (\psi(\phi_n,u))^2\, dC(u) & = & (\gamma^2_k + d_2+ \gamma^2_k d_0 +2 \gamma_k d_1) C(t_n) + o(C(t_n)) \label{equivintpsiphi}
\end{eqnarray}
The other terms are treated similarly (using the fact that $\psi(g_n,u)=1$ when $u\leq t_n$, and $=\Fbark(u)/\Fbarktn$ when $u>t_n$) and we obtain
\begin{eqnarray}
\int (\psi(g_n,u))^2\, dC(u) & = & \left( 1+ d_0 \right) C(t_n) + o(C(t_n)), \label{equivintpsig} \\
\int \psi(\phi_n,u)\psi(g_n,u)\, dC(u) & =& \left( \gamma_k + d_1 + \gamma_k d_0 \right) C(t_n) + o(C(t_n)). \label{equivintpsiphig}
\end{eqnarray}
In view of (\ref{expressionVarW1nInit}), combining $(\ref{equivintpsiphi})$, $(\ref{equivintpsig})$ and $(\ref{equivintpsiphig})$ and using Remark \ref{equivC} (following Lemma \ref{lem-Cdelta}) to write that $C(t_n)\sim (1-p)/\Hbar(t_n)$ (as $n\tinf$), this  proves the second term in relation (\ref{expressionVarW1n}).

\subsubsection{Proof of Lemma \ref{lem-Lyapunov}} \label{subsec-Lyapunov} \zdeux 

We have to prove that, for some $\delta >0$ small enough, $n \bE \left| W_{1,n}  \right|^{2 +\delta}$ tends to $0$, as $n\tinf$.  In the sequel, $cst$ denotes an unspecified absolute positive constant. According to the definition of  $W_{1,n} $, it is clear that
\begin{eqnarray*}
 n\left| W_{1,n}  \right|^{2 +\delta}  & \leq  &  cst \ \frac{v_n^{1+\delta/2}}{n^{1+\delta}} \left( \sum_{j=1}^3 \left| \Uunnb{j}  \right|^{2+\delta}   +  \gamma_k^{2+\delta}  \sum_{j=1}^3 \left| \Vunnb{j}  \right|^{2+\delta}   +   \left| \gammank -   \gamma_k \right|^{2+\delta} \right)  
\end{eqnarray*}
First, we clearly have $n^{-1-\delta} \  v_n^{1+\delta /2}   \left| \gammank -   \gamma_k \right|^{2+\delta} \lra 0$ as $n\tinf$.  Secondly, since $\Vunnb{j}$ has the same form as $\Uunnb{j}$, with $g_n$ instead of $\phi_n$ ({\it i.e.} without the log factor), we will only prove that there exists some $\delta>0$ such that, as $n\tinf$, 
\begeq{U1nLyap}
n^{-1-\delta} \  v_n^{1+\delta /2} \bE \left| \Uunnb{j}  \right|^{2+\delta}  = n^{-\delta /2} \left( \Fbarktn  \bar{G}(t_n) \right)^{1+\delta /2}   \bE \left| \Uunnb{j}  \right|^{2+\delta} \ \tqdninf \ 0 \ \mbox{\hset{0.3cm} for  } j \in \{ 1,2,3\}
\fineq

For $j=1$, we have   
\begin{eqnarray*}
 \bE \left| \Uunnb{1}  \right|^{2+\delta} & = &  \int_0^{+\infty}   \left|  \frac{\phi_n(z)}{\bar{G}(z)} \right|^{2+\delta} d\Hunk(z) 
 \; = \;  ( \Fbark(t_n) )^{-2-\delta} \int_{t_n}^{+\infty} (\Gbar(z))^{-2-\delta} \left( \log (z / t_n) \right)^{2+\delta}\bar{G}(z) d \Fk(z) 
 \\
 & = &  \left( \Fbark(t_n) \bar{G}(t_n) \right)^{-1-\delta} \int_{t_n}^{+\infty} \left( \log (z / t_n) \right)^{2+\delta} \left( \frac{\bar{G}(t_n)}{\bar{G}(z)}\right)^{1+\delta} \ \frac{d \Fk(z)}{\Fbarktn} .
\end{eqnarray*}
Applying part $(ii)$ of Lemma \ref{lemtechInteg} for $\alpha=2+\delta$, $a=(1+\delta)/\gamma_C$ and $b=-1/\gamma_k$ (with $\delta$ sufficiently small so that $a+b=1/\gamma_C-1/\gamma_k+\delta/\gamma_C$ is kept $<0$), and using the fact that $v_n=n\Fbarktn\Gbar(t_n)  \tinf$, this ends the proof of (\ref{U1nLyap}) for $j=1$.  
\zdeux 

For $j=2$, we have  
\[
\bE \left| \Uunnb{2}  \right|^{2+\delta} =  \int_0^{+\infty}    \left|  \frac{\psi(\phi_n,z)}{\bar{H}(z)} \right|^{2+\delta} \Fbar(z) dG(z).
\]
By definition of $\psi$, $\phi_n$, and $\gammank$, we have $\psi(\phi_n,z)=\gammank$ when $z\leq t_n$. Therefore,  splitting the integral above into two integrals $\int_0^{t_n}$ and $\int_{t_n}^{+\infty}$ we obtain 
\[
\bE \left| \Uunnb{2}  \right|^{2+\delta} = I_1(t_n) + I_2(t_n), 
 \]
 where, on one hand,
 \[
 I_1(t_n)  
 \; =  \; (\gammank)^{2+\delta} \int_0^{t_n} \frac{\bar{F}(z) dG(y)}{(\bar{H}(z))^{2+\delta}}  
 \; = \;  (\gammank)^{2+\delta}  \int_0^{t_n} \frac{dC(z)}{(\bar{H}(z))^{\delta}} 
 \; \leq \;  (\gammank)^{2+\delta}  \frac{C(t_n)}{(\bar{H}(t_n))^{\delta}}
\]
and, on the other hand, using the technical Lemma \ref{lemFonctionPsi}, for some $\delta'>0$,  \begin{eqnarray*}
 I_2(t_n)  
 & \leq  & 
 \int_{t_n}^{+\infty} \left|  \log\left( \frac{z}{t_n} \right) \fracFbark{z}{t_n} + \gamma_k \left( \frac{z}{t_n} \right)^{-1/ \gamma_k} + \epsilon_n(u)  \left( \frac{z}{t_n} \right)^{-1/ \gamma_k + \delta'} \right|^{2+\delta}  \  \frac{\bar{F}(z) dG(z)}{(\bar{H}(z))^{2+\delta}} \\
 & \leq &  
  cst \left\{  \int_{t_n}^{+\infty}   \log^{2+\delta}\left( \frac{z}{t_n} \right) \left( \fracFbark{z}{t_n} \right)^{2+\delta}  \  \frac{dC(z)}{(\bar{H}(z))^{\delta}} \right. \\
 & &  
 \hset{0.5cm} + \; \left.   \gamma_k^{2+\delta} \int_{t_n}^{+\infty}    \left( \frac{z}{t_n} \right)^{-(2+\delta)/ \gamma_k} \  \frac{dC(z)}{(\bar{H}(z))^{\delta}} 
 \; + \;  
 \sup_{u>t_n}|\epsilon_n(u)|^{2+\delta}   \int_{t_n}^{+\infty}  \left( \frac{z}{t_n} \right)^{(-1/ \gamma_k + \delta')(2+\delta)}  \  \frac{dC(z)}{(\bar{H}(z))^{\delta}}  \right\}  \\
 & = &  
  cst \frac{C(t_n)}{(\bar{H}(t_n))^{\delta}} \left\{  \int_{t_n}^{+\infty}   \log^{2+\delta} \left( \frac{z}{t_n} \right) \left( \frac{\Fbark(z)}{\Fbarktn} \right)^{2+\delta}   \left(\frac{\bar{H}(t_n)}{\bar{H}(z)} \right)^{\delta}  \  \frac{dC(z)}{C(t_n)}  \right.  \\
 & & 
  +  \left.    \gamma_k^{2+\delta} \int_{t_n}^{+\infty} \!    \left( \frac{z}{t_n} \right)^{-(2+\delta)/ \gamma_k}  \!\left(\frac{\bar{H}(t_n)}{\bar{H}(z)} \right)^{\delta}  \!  \frac{dC(z)}{C(t_n)}  
  + 
 o(1)   \int_{t_n}^{+\infty} \! \left( \frac{z}{t_n} \right)^{(-1/ \gamma_k + \delta')(2+\delta)} \! \left(\frac{\bar{H}(t_n)}{\bar{H}(z)} \right)^{\delta} \! \frac{dC(z)}{C(t_n)}  \right\}  
\end{eqnarray*}
Applying Lemma \ref{lem-Cdelta} to $\delta=1$, we have $C(t_n) = O \left(  1/\Hbar(t_n) \right)$, therefore $ I_1(t_n) = O  \left(  (\bar{H}(t_n))^{-1-\delta} \right)$. It is then easy to check that  $n^{-\delta /2} \left( \Fbarktn  \bar{G}(t_n) \right)^{1+\delta /2}   \   I_1(t_n) $ tends to 0, because $\Fbark\leq\Fbar$ and $n\Hbar(t_n)\tinf$, since $\Hbar(t_n) \geq \Fbark(t_n) \Gbar(t_n)$. 
\zdeux 

\noindent For $ I_2(t_n)$, since by Lemma \ref{lem-Cdelta} the function $C$ is regularly varying with index $1/\gamma$, the application of  part $(ii)$ of Lemma \ref{lemtechInteg} to $\alpha=0$ or $2+\delta$ and to various couples of values of $a$ and $b$ finally  yields $ I_2(t_n) = O  \left(  (\bar{H}(t_n))^{-1-\delta} \right)$, and consequently $n^{-\delta /2} \left( \Fbarktn  \bar{G}(t_n) \right)^{1+\delta /2}   \   I_2(t_n) $ tends to 0. 
\zdeux

We now come to the study of relation (\ref{U1nLyap}) for $j=3$. We have
\[
\bE \left| \Uunnb{3}  \right|^{2+\delta}  =  \int_0^{+\infty}  \left(  \int_0^z \psi(\phi_n,u) dC(u) \right)^{2+\delta} dH(z). 
\]
Proceeding as above by splitting the integral into two integrals $\int_0^{t_n}$ and $\int_{t_n}^{+\infty}$, we obtain 
\[
\bE \left| \Uunnb{3}  \right|^{2+\delta} = J_1(t_n) + J_2(t_n), 
 \]
where 
\[
 J_1(t_n) \; = \; (\gammank)^{2+\delta} \int_0^{t_n}  (C(z))^{2+\delta} dH(z)  
 \; = \;  - (\gammank)^{2+\delta} \bar{H}(t_n) \left( C(t_n) \right)^{2+\delta}  \int_0^{t_n}  \left( \frac{C(z)}{C(t_n)}  \right)^{2+\delta} \frac{d\bar{H}(z)}{\bar{H}(t_n)}
 \]
and 
\[
 J_2(t_n)  \; \leq \;  cst  (J^{(1)}_2(t_n)  + J^{(2)}_2(t_n)),
\]
where 
\[
  J^{(1)}_2(t_n) 
  \; = \;  \int_{t_n}^{+\infty} \!\! \left( \int_0^{t_n} \! \psi(\phi_n,u) dC(u) \right)^{2+\delta} \!\!\!\! dH(z)
  \; = \;  \gammank^{2+\delta}  \int_{t_n}^{+\infty} \!\! \left( C(t_n) \right)^{2+\delta} \! dH(z)  
  \; = \;  \gammank^{2+\delta}  \left( C(t_n) \right)^{2+\delta}   \bar{H}(t_n)
\]
and, using the technical Lemma \ref{lemFonctionPsi} as we did some lines above, 
  \begin{eqnarray}
  J^{(2)}_2(t_n) & =&   \int_{t_n}^{+\infty}   \left( \int_{t_n}^z   \psi(\phi_n,u)  dC(u) \right)^{2+\delta} dH(z) \nonumber\\
      &  \leq  &  cst \left(  \int_{t_n}^{+\infty}    \!\! \left(  \int_{t_n}^z  \log \!\left( \frac{u}{t_n} \right) \frac{\Fbark(u)}{\Fbarktn}  \! dC(u) \right)^{2+\delta} \!\!\!\!  dH(z)  + \gamma_k  \int_{t_n}^{+\infty}  \!\left(  \int_{t_n}^z  \left( \frac{u}{t_n} \right)^{-1/\gamma_k} dC(u) \right)^{2+\delta}\!\!\!\!  dH(z)   \right. \nonumber \\
      &&  \left. +  \sup_{u>t_n}|\epsilon_n(u)|^{2+\delta}      \int_{t_n}^{+\infty}    \left(  \int_{t_n}^z  \left( \frac{u}{t_n} \right)^{-1/\gamma_k+\delta'} dC(u) \right)^{2+\delta}  dH(z)  \right). \label{majJ22}
 \end{eqnarray}
Using Lemma \ref{lem-Cdelta} and part $(iii)$  of Lemma \ref{lemtechInteg}, we find that both $ J_1(t_n)$  and $J^{(1)}_2(t_n)$ are $O \left( (\bar{H}(t_n))^{-1-\delta} \right) $ and, though the term $J^{(2)}_2(t_n)$ is more involved,  we are also going to prove below that the same property holds for $J^{(2)}_2(t_n)$ :  this will finish the proof of Lemma \ref{lem-Lyapunov} because 
$ n^{-\delta /2} \left( \Fbarktn  \bar{G}(t_n) \right)^{1+\delta /2} (\bar{H}(t_n))^{-1-\delta}$ tends to $0$,
as already seen in the proof for $j=2$. 
 \zun

We only treat the first integral in the right-hand side of $(\ref{majJ22})$, since the two others are very similar, {\it i.e.} we need to prove that 
\begin{eqnarray}
 \int_{t_n}^{+\infty}     \left(  \int_{t_n}^z  \log \left( \frac{u}{t_n} \right) \frac{\Fbark(u)}{\Fbarktn}  dC(u) \right)^{2+\delta}  dH(z) = O \left( \frac{1}{(\bar{H}(t_n))^{1+\delta}} \right). \label{integFkLiap}
\end{eqnarray} 
Now, 
\[
 \int_{t_n}^{+\infty}     \left(  \int_{t_n}^z  \log \left( \frac{u}{t_n} \right) \frac{\Fbark(u)}{\Fbarktn}  dC(u) \right)^{2+\delta} \!\!\!\!   dH(z) \; = \;   (C(t_n))^{2+\delta}   \int_{t_n}^{+\infty}  \left(  \int_{1}^{z/t_n}  \log  (y)  \frac{\Fbark(y t_n)}{\Fbarktn}  \frac{dC(y t_n)}{C(t_n)} \right)^{2+\delta} \!\!\!\!   dH(z) 
\]
 Using Potter-bounds $(\ref{BornesPotter})$ for  $\Fbark \in RV_{-1/ \gamma_k}$,  integration by parts and then Potter-bounds $(\ref{BornesPotter})$ for  $C \in RV_{1/ \gamma}$, it is easy to see that for $n$ sufficiently large and $\epsilon >0$, there exists some positive constants $c$, $c'$, $c''$ such that
\begin{eqnarray*}
   \left(  \int_{1}^{z/t_n}  \log  (y)  \frac{\Fbark(y t_n)}{\Fbarktn}  \frac{dC(y t_n)}{C(t_n)} \right)^{2+\delta}     & \leq & c \log^{2+\delta} \left(  \frac{z}{t_n} \right)   \ \left(  \frac{z}{t_n} \right)^a  + c'   \left(  \frac{z}{t_n} \right)^a  + c''. 
\end{eqnarray*}
where $a=(2+\delta)(\frac{1}{\gamma} -\frac{1}{\gamma_k} + 2 \epsilon)$. Consequently 
\begin{eqnarray*}
   & & \int_{t_n}^{+\infty}     \left(  \int_{t_n}^z  \log \left( \frac{u}{t_n} \right) \frac{\Fbark(u)}{\Fbarktn}  dC(u) \right)^{2+\delta} dH(z)  
   \\ 
   & \leq  &  (C(t_n))^{2+\delta}    \bar{H}(t_n) \left(  - c    \int_{t_n}^{+\infty}   \log^{2+\delta} \left(  \frac{z}{t_n} \right)   \left( \frac{z}{t_n} \right)^a  \frac{d \bar{H}(z)}{\bar{H}(t_n)}  
  - \; c' \int_{t_n}^{+\infty}  \left(  \frac{z}{t_n} \right)^a   \frac{d \bar{H}(z)}{\bar{H}(t_n)}  \, - \, c'' \right).
\end{eqnarray*}
This yields $(\ref{integFkLiap})$, by using part $(ii)$ of Lemma \ref{lemtechInteg} to this value of $a$, to $b=-1/\gamma$ (and to $\alpha=2+\delta$ or $\alpha=0$), as well as Lemma \ref{lem-Cdelta}.

\subsection{Proof of Proposition \ref{PropGlobaleReste}} \label{sec-Reste}

Let us start with an important note. In Proposition \ref{PropGlobaleReste}, the main result is that the remainder terms $R_{n,\phi}$ and $R_{n,g}$ are $o_{\bP}(v_n^{-1/2})$. Proving this will be conducted in a similar way as proving that $R_n$ is $o_{\bP}(n^{-1/2})$ in Theorem 1.1 of \cite{Stute95}. But, recall that in our situation, the function   that we integrate here is $\phi_n$, which  is depending on $n$, with  a  "sliding" support $[t_n,+\infty[$.  We will need to be particularly cautious with integrability issues, especially when dealing with U-statistics for the terms $R_{n,2}$ and $R_{n,3}$ in the remainder $R_{n,C}$,  defined below.   
\zdeux

Before we proceed with the proof, let us define the following   empirical (sub)-distribution functions : for $t \geq 0$, 
\[
H_n(t)= \frac 1 n \sum_{i=1}^n\bI_{Z_i \leq t} ,  \ \  H_n^{(0)}(t)= \frac 1 n \sum_{i=1}^n\bI_{Z_i \leq t} \  \delta_i ,  \ \   H_n^{(1,k)}(t)= \sum_{i=1}^n\bI_{Z_i \leq t} \  \delta_i \ \bI_{\cause_i=k}.
\]

First note that, since $g_n$ is the function $\phi_n$ without the log factor, it should be clear to the reader that proving that $\Delta_n=\hat\Delta_n+R_{n,g}$ and $\racinevn R_{n,g}=o_{\bP}(1)$ will be simpler than proving that $\gamtildenk=\gamtchenk+R_{n,\phi}$ and $\racinevn R_{n,\phi}=o_{\bP}(1)$. We will thus only prove the latter two relations. 
\zdeux

Let us start with the first one, in other words let us define the remainder term $R_{n,\phi}$.  Remind that the definitions of $\gamtildenk$ and $\gamtchenk$ are $\gamtildenk  =  \int \phi_n(u) d\Fnk(u)$ and  $\gamtchenk=\Ubarn{1}+\Ubarn{2}-\Ubarn{3}$, where $\Ubarn{j}$ denotes the mean of the $n$ variables $\Uin{j}$. We need to decompose the integral of $\phi_n$ with respect to $\Fnk$, which is a stepwise subdistribution function which jumps at the (ordered) observations $\Zi$ are equal to $\indksipk / (n\GbarnZimun)$. But it is known that (see Lemma 2.1 in \cite{Stute95})  
\[
 \frac 1 {\GbarnZimun} = \exp \left\{ n \int_0^{\Zi^-} \log(1+(n\Hbarn(x))^{-1})\, d\HnO(x)\right\}
\]
Therefore, using the fact that $\Gbar(z) = \exp\left( -\int_0^z \Hbar^{-1}d\HO \right)$, we have 
\begin{eqnarray*}
  \gamtildenk & = &  
    \frac 1 n \sum_{i=1}^n \frac{\phi_n(\Zi)}{\Gbar(\Zi^-)} \indksipk \frac{\Gbar(\Zi^-)} {\GbarnZimun}  \\
    & = & 
    \frac 1 n \sum_{i=1}^n \frac{\phi_n(Z_i)}{\Gbar(\Zim)} \indksik 
    \exp\left(  n \int_0^{\Zim} \log(1+(n\Hbarn(x))^{-1})\, d\HnO(x)   \; - \;  \int_0^{\Zim} \frac{d\HO}{\Hbar}\right).
\end{eqnarray*}
Consequently, using the mean value theorem for $\exp$, and introducing the important notations
\begin{eqnarray*}
\Bin & = & n \int_0^{\Zim} \log(1+(n\Hbarn(x))^{-1})\, d\HnO(x)   \; - \;  \int_0^{\Zim} \frac{d\HnO}{\Hbarn} \\
\Cin & = & \int_0^{\Zim} \frac{d\HnO}{\Hbarn} \; - \; \int_0^{\Zim} \frac{d\HO}{\Hbar},
\end{eqnarray*}
it is easy to see that 
\begin{eqnarray}
  \gamtildenk & = &    
  \frac 1 n \sum_{i=1}^n \frac{\phi_n(Z_i)}{\Gbar(\Zim)} \indksik  
  \; + \;
  \frac 1 n \sum_{i=1}^n \frac{\phi_n(Z_i)}{\Gbar(\Zim)} \indksik  \Bin 
  \; + \;
  \frac 1 n \sum_{i=1}^n \frac{\phi_n(Z_i)}{\Gbar(\Zim)} \indksik  \Cin
  \nonumber \\ & & \mbox{\null ${}^{}$} \hset{6.cm}  \; + \;
  \frac 1 {2n} \sum_{i=1}^n \phi_n(Z_i) \indksik (\Bin+\Cin)^2   e^{\Delta_{i,n}}
   \nonumber \\  
  & = & 
  \Ubarn{1} 
  \; + \; R_{n,B} 
  \; + \; \frac 1 n \sum_{i=1}^n \frac{\phi_n(Z_i)}{\Gbar(\Zim)} \indksik  \Cin
  \; + \; R_{n,\Delta}, \label{decompgammatilde}
\end{eqnarray}
where $\Ubarn{1}$ is the first term in the definition of $\gamtchenk$, and $\Delta_{i,n}$ is a random quantity lying between $ \int_0^{\Zim} \Hbar^{-1} d\HO$  and $n\int_0^{\Zim} \log(1+(n\Hbarn(x))^{-1})\, d\HnO(x)$. 
\zdeux\\
What we now need to do is to show that the term involving the quantity $\Cin$ in relation (\ref{decompgammatilde}) above can be written as $\Ubarn{2}-\Ubarn{3}$ plus a remainder term $R_{n,C}$, and therefore we have $\gamtildenk= \gamtchenk + R_{n,\phi}$, where 
\begeq{def-Rnphi}
 R_{n,\phi} \ = \ R_{n,B} + R_{n,C} + R_{n,\Delta}.
\fineq
 The rest of the proof will, afterwards, be devoted to showing that each term of $ R_{n,\phi}$ is $o_{\bP}(v_n^{-1/2})$.
\zdeux

Proceeding as in \cite{Stute95} or \cite{Suzukawa02}, and using the fact that for any given function $f$ we have $\int f \, d\Hnunk = \frac 1 n \sum_{i=1}^n f(Z_i)\indksik$, we can write
\begeq{decompCin}
  \frac 1 n \sum_{i=1}^n \frac{\phi_n(Z_i)}{\Gbar(\Zim)} \indksik  \Cin  \ = \ 
  -\CC{1} + 2\CC{2} - \CC{3} + \Rn{1}
\fineq
where 
\begin{eqnarray*}
\CC{1} & = &  \iiint \frac{\phi_n(z)}{\Gbar(z)\Hbar^2(v)} \bI_{z>v}\bI_{u>v} \, dH_n(u)d\HnO(v)d\Hnunk(z), \\
\CC{2} & = &  \iint \frac{\phi_n(z)}{\Gbar(z)\Hbar(v)} \bI_{z>v} \, d\HnO(v)d\Hnunk(z), \\
\CC{3} & = &  \iint \frac{\phi_n(z)}{\Gbar(z)\Hbar(v)} \bI_{z>v}  \, d\HO(v)d\Hnunk(z), \\
\Rn{1} & = &  \iint \frac{\phi_n(z)}{\Gbar(z)} \bI_{z>v}
 \frac{(\Hbarn-\Hbar)^2(v)}{\Hbar^2(v)\Hbarn(v)}  d\HnO(v)d\Hnunk(z) .
\end{eqnarray*}
Note that $\CC{1}$ and $\CC{2}$ are a kind of $U$-statistics, which need to be approximated by sums of independent variables called Hoeffding decompositions : more precisely, if we introduce the functions (important in the sequel)
\begeq{def-fonctionsh}
 h(v,w) \ = \ \frac{\phi_n(w)}{\Gbar(w^-)\Hbar(v)} \bI_{w>v} \bI_{v<\infty}\bI_{w<\infty} 
 \makebox[1.5cm][c]{and} \het(u,v,w) \ = \ h(v,w) \frac{\bI_{u>v}}{\Hbar(v)} 
\fineq
for $u\in\bR$, $v\in\bR\cup\{+\infty\}$ and $w\in\bR\cup\{+\infty\}$,  then these decompositions are defined by 
\begin{eqnarray} 
\Cchap{1} & = & 
\iiint \het(u,v,w) dH_n(u)d\HO(v)d\Hunk(w) + \nonumber \\
&  & \null \hset{0.7cm} 
\iiint \het(u,v,w) dH(u)d\HnO(v)d\Hunk(w) +  \nonumber \\
&  & \null \hset{1.4cm} 
\iiint \het(u,v,w) dH(u)d\HO(v)d\Hnunk(w) \nonumber \\
&  & \null \hset{2.cm} 
 - \; 2 \iiint \het(u,v,w) dH(u)d\HO(v)d\Hunk(w)  
\label{def-Cchap1} \\
\Cchap{2} & = & 
\iint h(v,w) d\HnO(v)d\Hunk(w) + \iint h(v,w) d\HO(v)d\Hnunk(w) \nonumber \\
&  & \null \hset{6.cm} 
- \iint h(v,w) d\HO(v)d\Hunk(w).  
\label{def-Cchap2}
\end{eqnarray}
Therefore, if we introduce the remainder terms 
\begeq{def-Rn2Rn3}
 \Rn{2} \ = \ \CC{1} - \Cchap{1}  
 \makebox[1.5cm][c]{and}
 \Rn{3} \ = \ \CC{2} - \Cchap{2}  
\fineq
then (\ref{decompCin}) becomes 
\[
  \frac 1 n \sum_{i=1}^n \frac{\phi_n(Z_i)}{\Gbar(\Zim)} \indksik  \Cin  \ = \ 
  -\Cchap{1} + 2\Cchap{2} - \CC{3} + R_{n,C}  
  \makebox[2.cm][c]{where} 
  R_{n,C} \; = \; \Rn{1}-\Rn{2}+2\Rn{3}.
\]
We are thus left  to prove that    $-\Cchap{1} + 2\Cchap{2} - \CC{3} = \ \Ubarn{2}-\Ubarn{3}$. This is indeed the case because, if we note 
\begeq{def-thetan}
 \mbox{$ \theta_n = \iint h(v,w) d\HO(v)d\Hunk(w) $},
\fineq
then, by definition of $\het$,  the last (fourth) term in $\Cchap{1}$ equals $-2\theta_n$, the third one equals $\CC{3}$, the second one is (because $d\Hunk(w)=\Gbar(w^-)d\Fk(w)$)
\begin{eqnarray*}
 \iiint \het(u,v,w) dH(u)d\HnO(v)d\Hunk(w)   
 & = & \iint h(v,w) d\HnO(v)d\Hunk(w)  \\
 & = & 
 \frac 1 n \sum_{i=1}^n  (1-\delta_i) \int h(Z_i,w)\Gbar(w^-)\, d\Fk(w)  \\
 & = & 
\frac 1 n \sum_{i=1}^n \frac{1-\delta_i}{\Hbar(Z_i)}  \int_{Z_i}^{\infty} \phi_n(w)\, d\Fk(w)  \\
 & = & 
 \frac 1 n \sum_{i=1}^n \frac{1-\delta_i}{\Hbar(Z_i)}  \psi(\phi_n,Z_i)   \ = \ \Ubarn{2} 
\end{eqnarray*} 
and the first one is (because $d\Hunk(w)=\Gbar(w^-)d\Fk(w)$ and $d\HO(v)=\Fbar(v)dG(v)$)
\begin{eqnarray*}
 \iiint \het(u,v,w) dH_n(u)d\HO(v)d\Hunk(w)   
 & = & 
 \frac 1 n \sum_{i=1}^n \iint \het(Z_i,v,w) d\HO(v)d\Hunk(w)  \\
 & = & 
\frac 1 n \sum_{i=1}^n \iint \frac{\phi_n(w)}{\Gbar(v)\Hbar(v)} \bI_{v<w} \bI_{v<Z_i} \, dG(v)\, d\Fk(w) \\
 & = &   \frac 1 n \sum_{i=1}^n \int_0^{Z_i} \psi(\phi_n,v) \, dC(v) 
 \ = \ \Ubarn{3} .
\end{eqnarray*} 
Likewise, the first term of $\Cchap{2}$ equals $\Ubarn{2}$, the second one equals $\CC{3}$, and the last one equals $-\theta_n$.  After straightforward simplifications, we obtain the desired equality $-\Cchap{1} + 2\Cchap{2} - \CC{3} = \ \Ubarn{2}-\Ubarn{3}$, and the proof of $\gamtildenk=\gamtchenk+R_{n,\phi}$ is over. 
\zdeux

The proof of Proposition \ref{PropGlobaleReste} is now based on the following two lemmas : Lemma \ref{lem-RnBRn1RnDelta} is proved in subsection \ref{subsec-RnB}, and Lemma \ref{lem-RnC} is the longest to establish, its proof will be split across subsections \ref{subsec-prelimsRn2Rn3} to \ref{subsec-CCtilde}.

\begin{lem} \label{lem-RnBRn1RnDelta}
If  conditions $(\ref{Ordre1})$ and  (\ref{condvntn}) hold with $\gamma_k<\gamma_C$,  then we have 
$$
\mbox{$\racinevn R_{n,B}=o_{\bP}(1)$, $\racinevn R_{n,1}=o_{\bP}(1)$, and $\racinevn R_{n,\Delta}=o_{\bP}(1)$.}
$$
\end{lem}

\begin{lem} \label{lem-RnC}
If  conditions $(\ref{Ordre1})$ and  (\ref{condvntn}) hold with $\gamma_k<\gamma_C$, then we have 
$$
\racinevn R_{n,j}=o_{\bP}(1) \mbox{\hset{0.5cm} for $j=2$ and for $j=3$} .
$$
\end{lem}

\subsubsection{Proof of Lemma \ref{lem-RnBRn1RnDelta}} \label{subsec-RnB}

$\bullet$  We start with the remainder term $R_{n,B} $, which is defined as 
\[
R_{n,B}  = \frac 1 n \sum_{i=1}^n \frac{\phi_n(Z_i)}{\Gbar(\Zim)} \indksik  \Bin, 
\]
where $ \Bin=n \int_0^{\Zim} \log(1+(n\Hbarn(x))^{-1})\, d\HnO(x)   \; - \;  \int_0^{\Zim} \frac{d\HnO}{\Hbarn}$.  Since, for all $x \geq 0$, \ $x-\frac{x^2}{2} \leq \log(1+x) \leq x$, we obtain
\[
-\frac{1}{2n}  \int_0^{\Zim} \frac{d\HnO(x)}{(\Hbarn(x))^2} \leq \Bin \leq 0
\] 
and  then 
\begin{eqnarray}
  |R_{n,B}|   
  & \leq  &  \frac 1 n \sum_{i=1}^n \left\{ \frac{\phi_n(Z_i)}{\Gbar(\Zim)} \indksik 
  \frac{1}{2n}  \left( \int_0^{\Zim} \frac{d\HnO(x)}{(\Hbarn(x))^2} \right) \right\} \nonumber  
  \\
  & \leq & \sup_{0\leq x < Z_{(n)} } \left( \frac{\Hbar(x)}{\Hbarn(x)} \right)^2 \ 
  \frac{1}{2n^2}  \sum_{i=1}^n \left\{\frac{\phi_n(Z_i)}{\Gbar(\Zim)} \indksik   
  \left( \int_0^{\Zim} \frac{d\HnO(x)}{(\Hbar(x))^2} \right) \right\} . 
  \label{inequRnB}
\end{eqnarray} 
But 
\[
  \int_0^{\Zim} \frac{d\HnO(x)}{(\Hbar(x))^2} =   \int_0^{\Zim}  \frac{d\HO(x)}{(\Hbar(x))^2} + \int_0^{\Zim}    \frac{d(\HnO-\HO)(x)}{(\Hbar(x))^2},
\]
so, if we define 
\[
T_n^{(1)} =  \frac 1 n \sum_{i=1}^n T_{i,n}^{(1)}  \mbox{ and } T_n^{(2)} =   \frac 1 n \sum_{i=1}^n T_{i,n}^{(2)}
\]
where 
\begin{eqnarray*} 
   T_{i,n}^{(1)}  & =  &   
   \frac 1 n  \frac{\phi_n(Z_i)}{\Gbar(\Zim)} \indksik   \int_0^{\Zim} \frac{d\HO(x)}{(\Hbar(x))^2} 
   \\
   T_{i,n}^{(2)}  & = &    
   \frac 1 n    \frac{\phi_n(Z_i)}{\Gbar(\Zim)} \indksik   \int_0^{\Zim} \frac{d(\HnO-\HO)(x)}{(\Hbar(x))^2},
\end{eqnarray*}
then it remains to prove (thanks to part $(i)$ of Lemma \ref{lemmeConsiderationsProcEmp}) that  $\racinevn T_n^{(1)}  = o_{\bP} (1)$ and $\racinevn T_n^{(2)}  = o_{\bP} (1)$. 
\zdeux

Concerning $T_{n}^{(1)}$, since $\Hbar \geq \HbarO$  implies that $T_{i,n}^{(1)} \leq  \frac 1 n  \frac{\phi_n(Z_i)}{\Gbar(\Zim)\HbarO(\Zim)} \indksik  $, then  $\racinevn T_n^{(1)}  = o_{\bP} (1)$ is a consequence  of Lemma \ref{lemGenRestes}, used with $\alpha=0$  and $d=1$. 
\zdeux

Concerning $T_{n}^{(2)}$, an integration by parts yields 
\begin{eqnarray*}
  \left|   \int_0^{\Zim} \frac{d(\HnO-\HO)(x)}{(\Hbar(x))^2} \right|  
  & \leq  &  \frac{|\HbarnO-\HbarO|(\Zim)}{\Hbar^2(\Zim)}  + 2 \int_0^{\Zim} \frac{|\HbarnO-\HbarO|(x)}{\Hbar^3(x)} \ dH(x) \, + \, |\HbarnO(0)-\HbarO(0)|
  \\
  &  \leq  & 
  \sup_{0\leq x < Z_{(n)} } \sqrt{n}  \frac{|\HbarnO-\HbarO|(x)}{(\HbarO(x))^{\!\frac{1}{2}- \alpha}}  \left(  \frac{\HbarO(x)}{\Hbar(x)} \right)^{\frac{1}{2}- \alpha}  \times
  \\ 
        & &    \left( \frac{1}{\sqrt{n} \left( \Hbar(\Zim) \right)^{\frac 3 2 + \alpha}}   + \frac{2}{\sqrt{n}}   \int_0^{\Zim}   \frac{dH(x)}{\left(\Hbar(x)\right)^{\frac 5 2 +\alpha}}  \right)  \, + \, |\HbarnO(0)-\HbarO(0)|,  
\end{eqnarray*}
for any given $0 < \alpha  < \frac 1 2 $.  Lemma \ref{lemmeConsiderationsProcEmp} (applied with $a=1/2-\alpha<1/2$) and the fact $\HbarO\leq\Hbar$ thus imply that 
\[
\left|   \int_0^{\Zim} \frac{d(\HnO-\HO)(x)}{(\Hbar(x))^2} \right|  \leq O_{\bP}(1) \frac{1}{\sqrt{n} \left( \Hbar(\Zim) \right)^{\frac 3 2 + \alpha}} \, + \, |\HbarnO(0)-\HbarO(0)|,
\]
 so that, by definition of $T_{i,n}^{(2)}$, the desired statement  $\racinevn T_n^{(2)}  = o_{\bP} (1)$ is a consequence   of Lemma \ref{lemGenRestes}, applied with $\alpha >0$ sufficiently small  and $d=\frac 3 2 $, and of
$$
 \frac {\racinevn}{n} \sum_{i=1}^n \frac 1 n    \frac{\phi_n(Z_i)}{\Gbar(\Zim)} \indksik |\HbarnO(0)-\HbarO(0)| 
 = \sqrt{n}|\HbarnO(0)-\HbarO(0)| \times \frac{\sqrt{\Fbarktn\Gbar(t_n)}}{n} \times \Ubarn{1} = o_{\bP}(1).
$$ 
Indeed $\Ubarn{1}$ converges to $\gamma_k$ and $\HbarnO(0)-\HbarO(0)$ equals $\frac 1 n \sum_{i=1}^n \bI_{\delta_i=0} - \bP(\delta=0)$, which is $O_{\bP}(n^{-1/2})$ by the standard central limit theorem.
\zdeux

$\bullet$ Let us now turn to the remainder term $R_{n,1}$, which is defined as  
\[
 R_{n,1} = \iint \frac{\phi_n(z)}{\Gbar(z)} \bI_{z>v}   \frac{(\Hbarn-\Hbar)^2(v)}{\Hbar^2(v)\Hbarn(v)}  d\HnO(v)d\Hnunk(z) .
\]
A simple calculation leads to 
\begin{eqnarray*}
 R_{n,1} & \leq  &  
 \sup_{0\leq x < Z_{(n)} }  \left(  \sqrt{n}  \frac{|H_n-H|(x)}{(\Hbar(x))^{\frac{1}{2}- \alpha}} \right)^2 \  
 \sup_{0\leq x < Z_{(n)} }  \frac{\Hbar(x)}{\Hbarn(x)}  
 \times  
 \frac{1}{n^2}  \sum_{i=1}^n \left\{ \frac{\phi_n(Z_i)}{\Gbar(\Zim)} \indksik  \left( \int_0^{\Zim} \frac{d\HnO(v)}{(\Hbar(v))^{2+2\alpha}} \right)\right\}
\end{eqnarray*}
for any $0< \alpha < \frac 1 2$. Taking $\alpha$ sufficiently small, the rest of the proof is very similar to the one for $R_{n,B}$ (compare to $(\ref{inequRnB})$)  and relies on Lemma \ref{lemmeConsiderationsProcEmp}  and Lemma \ref{lemGenRestes}. 
  \zdeux
   
$\bullet$ We can finally deal with the last remainder term  $R_{n,\Delta}$, defined as 
\[
R_{n,\Delta} =      \frac{1}{2n}   \sum_{i=1}^n \phi_n(Z_i) \indksik (\Bin+\Cin)^2   e^{\Delta_{i,n}},
 \]
where $\Delta_{i,n}$ is a random quantity lying between $a_n:=\int_0^{\Zim} \Hbar^{-1} d\HO$  and $b_n:=\int_0^{\Zim} \log(1+(n\Hbarn(x))^{-1})\, d\HnO(x)$. Since $e^a=1/\bar{G}(\Zim)$,  we have
\[
R_{n,\Delta} =   \frac{1}{2n}   \sum_{i=1}^n  \frac{\phi_n(Z_i)}{\bar{G}(\Zim)} \indksik   (\Bin+\Cin)^2   e^{\Delta_{i,n} -a_n}.
\]
 Since $b_n-a_n= \Bin + \Cin$, where $\Bin <0$ and $\Cin = \int_0^{\Zim} \frac{d\HnO}{\Hbarn} \; - \; \int_0^{\Zim} \frac{d\HO}{\Hbar}$, we clearly have
\[
e^{(\Delta_{i,n} -a)} \leq \max(1, e^{\Cin}) \leq e^{|\Cin|}.
\]
But $\Cin= \hat{\Lambda}_{n,G} (Z_i) - \Lambda_G (Z_i)$, where $\Lambda_G $ is the cumulative hazard function associated to $G$, and $\hat{\Lambda}_{n,G}$ its Nelson-Alen estimator. Relying on \cite{Zhou91} Theorem 2.1, we can deduce that $\sup_{1 \leq i \leq n} |\Cin| = O_{\bP}(1)$. Hence, $e^{(\Delta_{i,n} -a)}  = O_{\bP}(1)$. 

\noindent Now, 
\[
\Cin - \frac{1}{2n}  \int_0^{\Zim} \frac{d\HnO}{(\Hbarn)^2}   \leq \Bin + \Cin \leq \Cin. 
\]
By writing 
\[
\Cin = \int_0^{\Zim}   \frac{d(\HnO -\HO)}{\Hbarn}   + \int_0^{\Zim}  \left(  \frac{1}{\Hbarn} -   \frac{1}{\Hbar}  \right)  d\HO,
\]
we prove  (using Lemma \ref{lemmeConsiderationsProcEmp} and simple integrations as for the previous treatment of $T_n^{(2)}$ above) that  $|\Cin| \leq O_{\bP}(1)  \ 1/(\sqrt{n} (\Hbar(Z_i))^{1/2 +\alpha}) + |\HbarnO(0)-\HbarO(0)|$ for $0< \alpha < 1/2$. 
\zdeux

Hence, on one hand  $(\Cin)^2 \leq O_{\bP}(1)  (n (\Hbar(Z_i))^{1 +2 \alpha}) ^{-1} + O_{\bP}(n^{-1}) \leq  O_{\bP}(1)  (n (\HbarO(Z_i))^{1 +2 \alpha}) ^{-1} + O_{\bP}(n^{-1})$, and on the other hand 
\[
\left( \Cin - \frac{1}{2n}  \int_0^{\Zim} \frac{d\HnO}{(\Hbarn)^2}  \right)^2 \leq O_{\bP}(1) \left(   \frac{1}{n (\HbarO(Z_i))^{1 +2 \alpha}} +  \frac{1}{n^2 (\HbarO(Z_i))^2} \right) + O_{\bP}(n^{-1}), 
\]
for any given $0< \alpha < 1/2$ (where the $O_{\bP}(n^{-1})$ comes from $|\HbarnO(0)-\HbarO(0)|^2$, which does not depend on $i$). Therefore, it is  sufficient to prove that 
\[
\frac{1}{n}    \sum_{i=1}^n   \frac{\racinevn \log (Z_i / t_n)  \indksik  \indZitn }{\Fbarktn \bar{G}(\Zim) n (\HbarO(Z_i))^{1 +2 \alpha}  }  
\makebox[1.3cm][c]{and } \; \frac{1}{n}    \sum_{i=1}^n  \ \frac{\racinevn  \log (Z_i / t_n)  \indksik  \indZitn }{\Fbarktn \bar{G}(\Zim) n^2 (\HbarO(Z_i))^2  }    
\]
 are $o_{\bP}(1)$, and that
\[
n^{-1/2} \left( \Fbarktn\Gbar(t_n) \right)^{1/2} \times \frac{1}{n}    \sum_{i=1}^n   \frac{\log (Z_i / t_n)  \indksik  \indZitn }{\Fbarktn \bar{G}(\Zim)   }   
\]
is $o_{\bP}(1)$ as well. But the first two statements are consequences of Lemma $\ref{lemGenRestes}$ with  $\alpha >0$ sufficiently close to $0$ and, respectively, $d=1$ and $d=2$. And for the third statement, the expectation of the expression  turns out (thanks to Lemma \ref{lemtechInteg} part $(ii)$) to be equivalent to a constant times $n^{-1/2}( \Fbarktn\Gbar(t_n) )^{1/2}$, which tends to $0$.

\subsubsection{Preliminaries to the proof of Lemma \ref{lem-RnC}}
\label{subsec-prelimsRn2Rn3}

We start this section by introducing important objects, issued from an idea appearing (to the best of our knowledge) in \cite{Stute94}. We define the improper variables $(V_i)_{1\leq i\leq n}$ and $(W_j)_{1\leq j\leq n}$ by 
\[
 V_i = \left\{ \begar{ll} Z_i & \mbox{ if $\delta_i=0$} \\ +\infty & \mbox{ if $\delta_i=1$} \finar \right. 
 \makebox[1.5cm][c]{and}
 W_j = \left\{ \begar{ll} +\infty & \mbox{ if $\delta_j=0$ or $\cause_j\neq k$} \\ Z_j & \mbox{ if $\delta_j=1$ and $\cause_j=k$} \finar \right. 
\]
which have $\HO$ and $\Hunk$ for respective subdistribution functions. We thus have $1-\delta_i=\bI_{V_i<\infty}$ and $\bI_{\cause_j=k} = \bI_{W_j<\infty}$, which, according to the definitions of $\CC{1}$ and $\CC{2}$ on one hand, and of functions $h$ and $\het$ (in (\ref{def-fonctionsh})) on the other hand, leads to
\[
 \CC{2} \; = \;  \frac 1 {n^2} \sum_{i=1}^n\sum_{j=1}^n \frac{\phi_n(Z_j)}{\Gbar(Z_j^-)\Hbar(Z_i)} \bI_{Z_j>Z_i} (1-\delta_i)\delta_j\bI_{\cause_j=k} 
 \; = \; \frac 1 {n^2} \sum\sum_{\hset{-0.4cm}i\neq j} h(V_i,W_j) 
\] 
and 
\[
 \CC{1} \; = \;  \frac 1 {n^3} \sum_{l=1}^n\sum_{i=1}^n\sum_{j=1}^n \frac{\phi_n(Z_j)}{\Gbar(Z_j^-)\Hbar^2(Z_i)} \bI_{Z_j>Z_i} \bI_{Z_l>Z_i} (1-\delta_i)\delta_j\bI_{\cause_j=k} 
 \; = \;  \frac 1 {n^3} \sum\sum\sum_{\hset{-0.8cm}i\neq j,i\neq l} \het(Z_l,V_i,W_j).
\]
Since the latter triple sum is not convenient, we also define
\[
 \CC{1} \; = \; \CCtilde{1} + \CCdoubletilde{1}  \makebox[1.8cm][c]{where} 
 \CCtilde{1} = \frac 1 {n^3} \sum\sum\!\!\sum_{\hset{-0.95cm}i,j,l \ distincts} \het(Z_l,V_i,W_j)
 \makebox[1.2cm][c]{and} 
 \CCdoubletilde{1} = \frac 1 {n^3} \sum\sum_{\hset{-0.4cm}i\neq j}  h(V_i,W_j)/\Hbar(V_i),
\]
where $\CCtilde{1}$ will be the quantity  approximated by $\Cchap{1}$, and $\CCdoubletilde{1}$ will be a remainder. We can indeed rewrite (\ref{def-Rn2Rn3}) as 
\begin{eqnarray}
\Rn{3} & = & {\textstyle \frac{n(n-1)}{n^2}} \left( {\textstyle\frac{n^2}{n(n-1)}}\CC{2} - \Cchap{2} \right) - \Cchap{2}/n  ,
 \label{decompRn3}
 \\
 \Rn{2} & = & {\textstyle \frac{n(n-1)(n-2)}{n^3}} \left( {\textstyle \frac{n^3}{n(n-1)(n-2)}} \CCtilde{1} - \Cchap{1} \right)  - {\textstyle \frac{3n-2}{n^2}} \Cchap{1} + \CCdoubletilde{1} . 
 \label{decompRn2}
\end{eqnarray}  
The terms in parentheses in (\ref{decompRn3}) and (\ref{decompRn2}) turn out to be genuine U-statistics of 2 and 3 variables, denoted by
\begeq{defcalUncalVn}
  \calUn = {\textstyle \frac{1}{n(n-1)}} \sum\sum_{\hset{-0.4cm}i\neq j} \gdH(V_i,W_j) 
  \makebox[1.3cm][c]{and}  
  \calVn = {\textstyle \frac{1}{n(n-1)(n-2)}} \sum\sum\sum_{\hset{-0.95cm}i,j,l \ distincts} \gdHet(Z_l,V_i,W_j)
\fineq 
where functions $\gdH$ and $\gdHet$ will be defined in a few lines (relation (\ref{def-gdHgdHet})) after some preliminaries, certainly well-known in  the U-statistics literature, but which we include here to make our proof self-contained (and since we are dealing with improper variables).
\zdeux

If $V$ and $W$ denote {\it independent} improper random variables with subdistribution functions $\HO$ and $\Hunk$ ({\it i.e.} $V=Z\bI_{\delta=0} + \infty\bI_{\delta=1}$ and $W=Z' \delta'\bI_{\cause'=k} + \infty (1-\delta' + \bI_{\cause'\neq k})$ where $(Z,\delta,\cause)$ and $(Z',\delta',\cause')$ are independent copies of $(Z_1,\delta_1,\cause_1)$), we introduce the following notations : for any function $g:[0,\infty]\times[0,\infty]\rightarrow \bR$,
\[
 g_{1\bullet} (v) = \bE(g(v,W)) \makebox[1.5cm][c]{and}  g_{{\bullet}1} (w) = \bE(g(V,w))\ ,
\]
as well as, for any function $g:[0,\infty[\times[0,\infty]\times[0,\infty]\rightarrow \bR$, with $Z$ (of distribution function $H$)  independent of $V$ and $W$, 
\[
 g_{1\bullet\bullet} (u) = \bE(g(u,V,W)) \makebox[1.cm][c]{,}  
 g_{{\bullet}1\bullet} (v) = \bE(g(Z,v,W)) \makebox[1.5cm][c]{and}  
 g_{{\bullet\bullet}1} (w) = \bE(g(Z,V,w)).
\]
Since $h(v,w)=0$ whenever $v$ or $w$ equals $\infty$, we then have (the proof is simple) 
$$
 \theta_n =\iint h(v,w)d\HO(v) d\Hunk(w) = \bE(h(V,W)) = \bE(\het(Z,V,W)).
$$
Therefore, setting (for  $z$ in $[0,\infty[$ and $v$ and $w$ in $[0,\infty]$)
\begeq{def-gdHgdHet}
\begar{rcl} 
\gdH(v,w) & = & h(v,w) - \hunpt(v) - \hptun(w) +\theta_n \zundemi\\
\gdHet(z,v,w) & = & h(z,v,w) - \hetunptpt(z) - \hetptunpt(v) - \hetptptun(w) + 2\theta_n 
\finar
\fineq
it is then not difficult  to check (using (\ref{def-Cchap1}) and (\ref{def-Cchap2})) that $\calUn$ and $\calVn$  in relation (\ref{defcalUncalVn}) are indeed equal to the differences  in parentheses in relations (\ref{decompRn3}) and (\ref{decompRn2}), respectively.  Lemma \ref{lem-RnC}  thus becomes a consequence of the following facts : $\racinevn\calUn=o_{\bP}(1)$, $\racinevn\calVn=o_{\bP}(1)$, and
\begeq{rel-melimeloCCtilde}
 \mbox{ the three sequences \ $\CCdoubletilde{1}$ \ ,  \ $\Cchap{1}/n$  \ and \ $\Cchap{2}/n$   \ are  \ $o_{\bP}(v_n^{-1/2})$. }
\fineq  
We will prove these statements in the next 3 subsections.

\subsubsection{Proof of $\racinevn \calUn = o_{\bP}(1)$} \label{subsec-preuveUn}
\zun

We note ${\cal I}=\{I=(i,j) \, ; \, 1\leq i<j\leq n\}$, $\gdH_I=\gdH(V_i,W_j)$ when $I=(i,j)\in{\cal I}$, and $N=n(n-1)/2$. It is clear that it suffices to prove that 
$$
 S_N = o_{\bP}(1)  
 \makebox[1.8cm][c]{where} S_N = \sum_{I\in{\cal I}} {\textstyle\frac{\racinevn}N} \gdH_I . 
$$
The good point  is that $S_N$ turns out to be a sum of identically distributed centred and uncorrelated random variables $\gdH_I$, but unfortunately these variables  $\gdH_I$ are not square-integrable and potentially only have a moment of order slightly larger than $4/3$ when $\gamma_k<\gamma_C$. In order to deal with this difficulty, since we cannot handle directly the $L^p$ norm of $S_N$ of order $p=4/3$, we will follow a strategy similar to that found in \cite{Csorgo2008}, based on truncation. We set
\begeq{def-gdHs-gdHun-Mn}
 \gdHs(v,w) \; = \; \gdH(v,w)\bI_{|\gdH(v,w)|\leq M_n} \, - \, \bE(\gdHun\bI_{|\gdHun|\leq M_n}) 
 \makebox[1.4cm][c]{where}
 \left\{ \begar{lcl} \gdHun & = & \gdH(V_1,W_2) \\
 M_n & = & n^2/\racinevn \finar\right.
\fineq
The variables $\gdHs_I \, = \, \gdHs(V_i,W_j)$ ($I\in{\cal I}$) are centred and bounded, but they lose the non-correlation property of the variables $\gdH_I$. This is why we define now
\[
  \gdHss_I \; = \; \gdHss(V_i,W_j) 
  \makebox[1.4cm][c]{where}
  \gdHss(v,w) = \gdHs(v,w) - \gdHsunpt(v) - \gdHsptun(w)
\]
which are centred and bounded but are also uncorrelated  (see part $(i)$ of Lemma \ref{diversPreuveUstats}), and we write
\begeq{decompSN}
\textstyle
 S_N \;  = \;  S_N^{(1)} + S_N^{(2)} 
\; = \; \frac{\racinevn}{N} \sumI \gdHss_I \, + \, \frac{\racinevn}{N} \sumI (\gdH_I - \gdHss_I).
\fineq
We thus need to prove that $S_N^{(1)}$ and $S_N^{(2)}$ both converge to $0$ in probability. 
\zdeux

Concerning $S_N^{(1)}$, since the $\gdHss_I$ are centred and uncorrelated, we have
\[ 
 \textstyle
 \bE\left((S_N^{(1)})^2\right) 
 \, = \, (v_n/N^2) \bE\left( \left(\sumI \gdHss_I\right)^2\right) 
 \, = \, (v_n/N) \bE\left( (\gdHss(V_1,W_2))^2\right)
 \, \leq \, cst (v_n/n^2) \bE(\gdHun^2\bI_{|\gdHun|\leq M_n})
\]
where $\gdHun$ was defined in (\ref{def-gdHs-gdHun-Mn}) (the justification of the last inequality is postponed to part $(ii)$ of Lemma \ref{diversPreuveUstats}). Remind that $\gdHun$ is not square-integrable and $M_n=n^2/\racinevn = n^{3/2}/(\Fbarktn\Gbar(t_n))^{1/2}$, and introduce $m_n = n^{3/2}/(\Fbarktn\Gbar(t_n))^{1/2-\epsilon}=o(M_n)$ for some given small $\epsilon>0$. We then write 
\[
 \bE\left((S_N^{(1)})^2\right) 
 \, \leq \, cst \frac{v_n}{n^2} m_n^{2/3} \bE\left(|\gdHun|^{4/3}\right) 
 \; + \;  cst \frac{v_n}{n^2} M_n^{2/3} \bE\left(|\gdHun|^{4/3}\bI_{|\gdHun|>m_n}\right)
 \; = \; cst(A_n + B_n).
\]
Thanks to Lemma \ref{lem-integhetH} (parts $(i)$ and $(ii)$) and to the definition of $m_n$, the term $A_n$ is bounded by a quantity which is equivalent (as $n\tinf$) to $\frac{v_n}{n^2} m_n^{2/3} \left(\frac{v_n} n \right)^{-2/3} = (\Fbarktn\Gbar(t_n))^{2\epsilon/3}=o(1)$. 
We now rely on H\"older's inequality for dealing with the term $B_n$ . Let $p>1$ and $q>1$ such that $1/p+1/q =1$.  Since $\theta_n=\bE(\gdHun)$, again thanks to Lemma \ref{lem-integhetH} ($(i)$, $(ii)$ and $(v)$), for $p$ sufficiently close to $1$ so that $4p/3 < 1+(1+2\gamma_k/\gamma_C)^{-1}$,  we have
\begin{eqnarray*}
 B_n 
 & \leq & \left(\frac{v_n} n \right)^{2/3} \! \left(\bE |\gdHun|^{4p/3}\right)^{1/p} \left( \bP(|\gdHun|>m_n) \right)^{1/q} \\
 & \leq &  \left(\frac{v_n} n \right)^{2/3} \! \left( O\left( \left(\frac {v_n}n\right)^{2(1-4p/3)} \right) \right)^{1/p} \, m_n^{-1/q} \, (  4\theta_n )^{1/q}  \\
 & \leq & O(1) \left(\frac{v_n} n \right)^{2/3 + 2(1-1/q) - 8/3} M_n^{-1/q}  (\Fbarktn\Gbar(t_n))^{-\epsilon/q} (-\log\Gbar(t_n))^{1/q} \\
 & \leq &   O(1) \, v_n^{-3/2q}   (\Fbarktn\Gbar(t_n))^{-2\epsilon/q} \left( (\Gbar(t_n))^{\epsilon}(-\log\Gbar(t_n))\right)^{1/q} \; = \; o(1) \left( n(\Fbarktn\Gbar(t_n))^{1+4\epsilon/3}\right)^{-3/2q} \end{eqnarray*}
which converges to $0$ thanks to assumption (\ref{condvntn}), for $\epsilon>0$ small enough (we used part $(v)$ of Lemma \ref{lem-integhetH} in the third upper bound).  
\zun

We are thus left to prove that $S_N^{(2)}$ also converges to $0$, but this time in $L^1$. We start by writing that 
\begin{eqnarray*}
 \textstyle
 \bE\left(|S_N^{(2)}|\right) & \leq & \textstyle
 \frac{\racinevn} N \sumI \bE\left( |\gdH_I-\gdHss_I|\right) 
 \, = \, \racinevn \, \bE\left(\, |\gdHun-\gdH^*(V_1,W_2)+\gdHsunpt(V_1)+\gdHsptun(W_2)| \,\right) 
 \\ & \leq & 4\racinevn \,  \bE\left(|\gdHun|\bI_{|\gdHun|>M_n}\right),
\end{eqnarray*}
the last inequality being proved in the appendix (part $(iii)$ of Lemma \ref{diversPreuveUstats}). 
The follow-up is a bit similar to the treatment of $B_n$ above, relying on Lemma \ref{lem-integhetH} (parts $(i)$, $(ii)$ and $(v)$)  and on H\"older's inequality : for $p>1$ close to $1$ and a large $q$ such that $1/p+1/q=1$, we can write 
\begin{eqnarray*}
 \textstyle
 \bE\left(|S_N^{(2)}|\right) 
 & \leq & 4 \racinevn M_n^{-1/3} \left(\bE |\gdHun|^{4p/3}\right)^{1/p} \left( \bP(|\gdHun|>M_n) \right)^{1/q} 
 \\ 
 & \leq & O(1) \, v_n^{-3/2q} (-\log\Gbar(t_n))^{1/q}
 \; \leq \;  O(1) \left\{ \left( n(\Fbarktn\Gbar(t_n))^{1+\epsilon} \right)^{-3/2} (\Gbar(t_n))^{3\epsilon/2}   (-\log\Gbar(t_n)) \right\}^{1/q} 
\end{eqnarray*}
which, for $\epsilon>0$ small enough, is $o(1)$ thanks to assumption (\ref{condvntn}).

\subsubsection{Proof of $\racinevn \calVn = o_{\bP}(1)$} \label{subsec-preuveVn}
\zun

The proof is very similar to the one contained in the previous subsection. We nonetheless provide a few details to convince the reader of the validity of the result.  We note now ${\cal I}=\{I=(i,j,l) \, ; \, 1\leq i<j<l\leq n\}$ and $\gdHet_I=\gdHet(Z_l,V_i,W_j)$ when $I=(i,j,l)\in{\cal I}$, with $N=n(n-1)(n-2)/6$ denoting the cardinal of the index set ${\cal I}$. Since the observations $(Z_i)_{i\leq n}$ are i.i.d., it should be clear to the reader that  it suffices to prove that 
$$
 S_N = o_{\bP}(1)  
 \makebox[1.8cm][c]{where} S_N = \sum_{I\in{\cal I}} {\textstyle\frac{\racinevn}N} \gdHet_I.
$$
As previously, the problem lies with the moments of the centred and uncorrelated variables $\gdHet_I$, and now we only have a guaranteed moment of order slightly more than $6/5$ instead of $4/3$ in the previous situation. Fortunately, the cardinal $N$ is now of order $n^3$, which turns out to be the right compensation.  
\zun

We thus define, for $(u,v,w)\in [0,\infty[\times[0,\infty]\times[0,\infty]$,
\[
 \gdHets(u,v,w) \; = \; \gdHet(u,v,w)\bI_{|\gdHet(u,v,w)|\leq M_n} \, - \, \bE(\gdHetun\bI_{|\gdHetun|\leq M_n}) 
 \makebox[1.4cm][c]{where}
 \left\{ \begar{lcl} \gdHetun & = & \gdHet(Z_3,V_1,W_2) \\
 M_n & = & n^3/\racinevn \finar\right.
\]
as well as 
\[
  \gdHetss_I \; = \; \gdHetss(Z_l,V_i,W_j) 
  \makebox[1.4cm][c]{where}
  \gdHetss(u,v,w) = \gdHets(u,v,w) - \gdHetsunptpt(u) - \gdHetsptunpt(v) - \gdHetsptptun(w)
\]
which are centred and bounded but are also uncorrelated  (see part $(i)$ of Lemma \ref{diversPreuveUstats} in the Appendix), and we write
\[
\textstyle
 S_N \;  = \;  S_N^{(1)} + S_N^{(2)} 
\; = \; \frac{\racinevn}{N} \sumI \gdHetss_I \, + \, \frac{\racinevn}{N} \sumI (\gdHet_I - \gdHetss_I).
\]
Introducing $m_n=M_n(\Fbarktn\Gbar(t_n))^{\epsilon}$ and skipping details, we assess that
\[
 \bE\left((S_N^{(1)})^2\right) 
 \, \leq \, cst \frac{v_n}{n^3} m_n^{4/5} \bE\left(|\gdHetun|^{6/5}\right) 
 \; + \;  cst \frac{v_n}{n^3} M_n^{4/5} \bE\left(|\gdHetun|^{6/5}\bI_{|\gdHun|>m_n}\right)
\]
and that this quantity converges to $0$, as $n\tinf$, thanks to parts $(i)$ and $(iii)$ of Lemma \ref{lem-integhetH}. The same argument is used to prove that $\bE(|S_N^{(2)}|)\tqdninf 0$. 
\zdeux

\subsubsection{Proof of relation (\ref{rel-melimeloCCtilde}) }  \label{subsec-CCtilde}

Let us first prove that, for some $d\in ]4/5,1[$, $\bE(|\racinevn \CCdoubletilde{1}|^d)$ tends to 0, as $n$ tends to infinity. Recall that $\CCdoubletilde{1} = \frac 1 {n^3} \sum\sum_{i\neq j} h(V_i,W_j)/\Hbar(V_i)$. Since $d<1$, we have
\begin{eqnarray*}
\bE(|\racinevn \CCdoubletilde{1} |^d) & \leq & n^{d/2-3d} \ \left( \Gbar (t_n) \Fbark(t_n) \right)^{d/2} n(n-1) \bE\left(|h(V_1,W_2)/\Hbar(V_1) |^d\right).
\end{eqnarray*}
According to part $(iv)$ of Lemma \ref{lem-integhetH}, the right-hand side of the  inequality above is $O(1)\; v_n^{2-5d/2}$, which tends to $0$, since $d> 4/5$, and so we are done. 
\zun\\

Let us now prove that $\bE(|\racinevn \Cchap{1}/n|)$ tends to 0, as $n$ tends to infinity. $\Cchap{1}$ is defined in $(\ref{def-Cchap1})$, where the expectation of each of the four integrals is $\theta_n$ : therefore, we only need to prove that $\frac{\racinevn}{n} \theta_n$ tends to $0$. This is straightforward using part  $(v)$ of Lemma \ref{lem-integhetH}.
\zun\\
 We can prove in a very similar way that $\bE(|\racinevn \Cchap{2}/n|)$ tends to 0, as $n$ tends to infinity.

\subsection{Proof of Proposition \ref{consistance}}   \label{PreuveConsistance}

Using the same notations as in the begining of Section \ref{sec-proofs}, we have, 
\[
\gamchapnk - \gamma_k \ = \ \Delta_n^{-1} \left(\frac{Z_n}{\racinevn}  +  R_n  +  (\gammank - \gamma_k) \    \right) . 
\]
The fact that $\frac{Z_n}{\racinevn} \stackrel{\bP}{\rightarrow} 0$ is due to the application of a triangular weak law of large numbers (see \cite{ChowTeicher1997} for example)  to  $\frac 1 n \sum \tildeUin$ and to $\frac 1 n \sum \tildeVin$. By carrefully following the proof of proposition \ref{PropGlobaleReste}  in Section \ref{sec-Reste}, we can see that $R_n=o_{\bP}(1)$. The condition $\gamma_k < \gamma_C$ is  not used, neither in the treatment of $\frac{Z_n}{\racinevn}$ nor in that of $R_n$. Details are omited.

\subsection{Proof of Corollary \ref{coroquantiles}}   \label{PreuveCoroQuantiles}

The proof is very similar to the one of Theorem 2 in \cite{WormsWorms16}, with $\gamma_k$ and $\Fbark$ here replacing $\gamma_1$ and $\Fbar$ there. For completeness, we provide some details about it. Reminding the notations $d_n=\Fbark(t_n) / p_n \tinf$ and  $\Delta_n = \frac{\Fbarnk(t_n)}{\Fbarktn}$, we easily write 
\[
 \frac{\xchapeau}{\xpn} - 1   \; = \;  \frac{t_n}{\xpn} (\Delta_n d_n)^{\gamchapnk} -1   
 \; = \; \Delta_n^{\gamchapnk} \left(   \frac{t_n}{\xpn}\  d_n^{\gamma_k} T_n^1  +  T_n^2 + T_n^3 \right), 
\zun
\]
where $T_n^1 := d_n^{\gamchapnk-\gamma_k} - 1$,  \ $T_n^2 :=  \frac{t_n}{\xpn} d_n^{\gamma_k} -1$ \ and \ $T_n^3 := 1 - \Delta_n^{-\gamchapnk}$. 
We are going to prove that both $T_n^2$ and $T_n^3$ are $o_{\bP}(\log d_n/\sqrt{v_n})$, and that $\frac{\sqrt{v_n}}{\log d_n} T_n^1  \stackrel{\cal d}{\longrightarrow}   {\cal N} \left(\lambda m , \sigma^2 \right)$ : this will conclude the proof, since both $\Delta_n$  (Corollary \ref{coroDeltan}) and $ \frac{t_n}{\xpn}\  d_n^{\gamma_k}$ tend to $1$. 
\zun

Concerning $T_n^1$, the mean value theorem yields 
\[
\frac{\sqrt{v_n}}{\log d_n} T_n^1 = \frac{\sqrt{v_n}}{\log d_n} \left( e^{(\gamchapnk-\gamma_k)\log(d_n)} - 1 \right) =  \sqrt{v_n} (\gamchapnk - \gamma_k) \exp(E_n), 
\]
where $| E_n |  \leq |\gamchapnk - \gamma_k| \log d_n$ and therefore  $E_n$ tends to $0$ in probability thanks to Theorem \ref{TLCGammachap} and  assumption $(\ref{conditiondn})$. The desired result for $T_n^1$ is then implied by Theorem \ref{TLCGammachap} again. 
\zun 

Concerning the fact that $T_n^2=o_{\bP}(\log d_n/\sqrt{v_n})$, the proof is completely similar to the evoked one in \cite{WormsWorms16}, so we omit it here (basically, this is based on some uniform regular variation implied by the assumed negativity of the second order parameter $\rho_k$, and on the assumption that $\racinevn g(t_n)$ converges). 
\zdeux

Finally, concerning $T_n^3$ we use the mean value theorem to write 
\[
\frac{\racinevn}{\log d_n} T_n^3 = \frac{\gamchapnk D_n^{-\gamchapnk-1}}{\log d_n} . \racinevn (\Delta_n -1),
\]
where $D_n$ lies between $\Delta_n$ and $1$. But Corollary \ref{coroDeltan} (and  the consistency of $\gamchapnk$) implies that $\gamchapnk D_n^{-\gamchapnk-1}\stackrel{\bP}{\longrightarrow} \gamma_k$ on one hand, and $\racinevn (\Delta_n -1)=O_{\bP}(1)$ on the other hand ;  therefore,  $\frac{\racinevn}{\log d_n} T_n^3=O(1/\log(d_n))=o(1)$.

\section{Appendix}

This appendix contains various results : some of them are used repeatedly in the proof of the main result (in particular Proposition  $\ref{PropBornesPotter}$,  Lemmas \ref{lemtechInteg},  and \ref{lemmeConsiderationsProcEmp}, and to a lesser extent Lemmas \ref{lemFonctionPsi} and \ref{lem-Cdelta}), the other ones concern parts  of the main proof which are postponed to the appendix for better clarity of the main flow of the proof  (Lemmas \ref{lemGenRestes},  \ref{lem-integhetH} and \ref{diversPreuveUstats}). 
\ztrois

\begin{defi} \label{defRV}   
An ultimately positive function $f$ : $\bR^+  \rightarrow \bR$  is {\it regularly varying} (at infinity) with index  $\alpha \in \bR$, if
\[
\lim_{t \rightarrow + \infty} \frac{f(tx)}{f(t)} = x^{\alpha} \hspace{0.2cm}  (\forall x >0). 
\]
This is noted $f \in RV_{\alpha}$.  If $\alpha=0$, $f$ is said to be slowly varying. 
\end{defi}

\begin{prop} \label{PropBornesPotter} (See  \cite{HaanFerreira06}  Proposition B.1.9) \\
Suppose $f \in RV_{\alpha}$. 
If $x < 1$ and $\epsilon >0$, then there exists $t_0=t_0(\epsilon)$ such that for every $t\geq t_0$, 
\[
(1-\epsilon) x^{\alpha+\epsilon} < \frac{f(tx)}{f(t)} < (1+\epsilon) x^{\alpha-\epsilon} 
\]
and if $x \geq  1$ , 
\begeq{BornesPotter}
(1-\epsilon) x^{\alpha-\epsilon} < \frac{f(tx)}{f(t)} < (1+\epsilon) x^{\alpha+\epsilon} .
\fineq
\end{prop}

\begin{lem} \label{lemtechInteg}
Let $x \in \bR_+^* $ , $\alpha \in \bR_+$, $\beta >-1$, and for $a$ and $b$ real numbers, $f$ and $g$ are two  regular varying functions at infinity, with index, respectively, $a$ and $b$. Then, as $t \rightarrow + \infty$,   
\begitem
\item[$(i)$] $ \displaystyle J_{\beta}(x) = \int_1^{+\infty}  \log^{\beta} (y) \ y^{-x-1} dy = \frac{\Gamma(\beta +1)}{x^{\beta +1}}$.
\item[$(ii)$] $ \displaystyle I_{\alpha,a,b} = \int_{1}^{+\infty}   \log^{\alpha} (y) \ \frac{f(yt)}{f(t)}  \ \frac{dg(yt)}{g(t)} \rightarrow \frac{b \Gamma(\alpha +1)}{(-a-b)^{\alpha +1}}$, if  $a+b <0$
\item[$(iii)$] $  \displaystyle J_{a,b} = \int_0^{1}  \frac{f(yt)}{f(t)}  \ \frac{dg(yt)}{g(t)} \rightarrow   \frac{b}{a+b}$, if $a+b > 0$
\finit
\end{lem}
Proof : 
\begitem 
\item[$(i)$] A simple change of variable and the definition of the $\Gamma$ function yields the result. 
\item[$(ii)$] For the sake of simplicity, we are going to treat  the case $a<0$ and $b<0$. The only difference for the other cases is the sign in front of the $\epsilon$ or $\epsilon'$ appearing below (coming from the application of $(\ref{BornesPotter})$ several times), which can  depend on the sign of $a$, $b$ or another constant, but does not affect the result.  Using Potter-bounds $(\ref{BornesPotter})$ for $f$ yields, for $n$ sufficiently large and $\epsilon >0$ ,
\[
 (1 + \epsilon ) \int_1^{+\infty}  \log^{\alpha} (y) \ y^{a + \epsilon}  \ \frac{dg(yt_n)}{g(t_n)}   \leq  I_{\alpha,a,b}  \leq  (1 -\epsilon ) \int_1^{+\infty}  \log^{\alpha} (y) \ y^{a  - \epsilon}  \ \frac{dg(yt_n)}{g(t_n)} .
\]
Let us treat  only the upper  bound  and the case $\alpha \neq 0$  (the other cases  being similar).  By integration by parts, with $a +b <0$, we have
\[
\int_1^{+\infty}  \!\!\!\log^{\alpha} (y) \ y^{a - \epsilon}  \ \frac{dg(yt_n)}{g(t_n)}  = -\alpha \int_1^{+\infty}  \log^{\alpha -1} (y) \ y^{a  -1 - \epsilon}  \ \frac{g(yt_n)}{g(t_n)}  dy - (a - \epsilon) \int_1^{+\infty}   \!\!\! \log^{\alpha} (y) \ y^{a  -1 - \epsilon}  \ \frac{g(yt_n)}{g(t_n)}  dy .
\]
 Using Potter-bounds $(\ref{BornesPotter})$ for $g$ yields, for $n$ sufficiently large and $\epsilon' >0$ 
 \[
\int_1^{+\infty}  \log^{\alpha} (y) \ y^{a -\epsilon}  \ \frac{dg(yt_n)}{g(t_n)} \leq -\alpha (1 - \epsilon') J_{\alpha-1} (-a-b + \epsilon + \epsilon') -  (a - \epsilon) (1 + \epsilon') J_{\alpha} (-a-b +  \epsilon - \epsilon') .
\]
Doing the same with the lower bound and making $\epsilon$ and $\epsilon'$ tend to $0$, yields the result after simplifications. 
\zun
\item[$(iii)$] As in $(ii)$, using Potter-bounds $(\ref{BornesPotter})$ for $f$, integration by parts and then again  $(\ref{BornesPotter})$ for $g$ yields the result. 
\finit

\begin{lem} \label{lem-Cdelta}
For any $\delta>0$, let $C_{\delta}$ denote the function 
$$ 
 C_{\delta}(t)=\int_0^t \frac{dG(v)}{\Gbar(v)\Hbar^{\delta}(v)}.
$$
Under condition $(\ref{Ordre1})$, this function is regularly varying of order $\delta/\gamma$ and  we have $C_{\delta}(t)\sim (\gamma/\gamma_C)/(\delta\Hbar^{\delta}(t))$,  as $t\ra +\infty$.
\end{lem}
Proof : by writing $\Hbar^{\delta}(t)C_{\delta}(t) = -\int_0^1 \frac{\Hbar^{\delta}(t)}{\Hbar^{\delta}(tu)} \frac{\Gbar(t)}{\Gbar(tu)} \frac{d\Gbar(tu)}{\Gbar(t)}$, the lemma is an immediate consequence of part $(iii)$ of Lemma \ref{lemtechInteg}, with $a+b = (\delta/\gamma + 1/\gamma_C) + (-1/\gamma_C) = \delta/\gamma>0$ and $-b/(a+b) = (\gamma/\gamma_C)/\delta$.   

\begin{rmk} \label{equivC}
In the Lemma above, $C_1$ is the important function $C$ introduced at the beginning of Section \ref{sec-proofs}, and thus  $C(t)\sim (\gamma/\gamma_C)/\Hbar(t) =  (1-\gamma/\gamma_F)/\Hbar(t)$, as $t \rightarrow + \infty$.  Hence, $C$ is regularly varying at infinity with index $1/\gamma$, a property which proves useful several times in the main proofs. 
\end{rmk}

\begin{lem} \label{lemFonctionPsi}
Let $\psi(\phi_n,u)= \int_u^{+\infty} \phi_n(s) d F^{(k)}(x)$, for $u \geq 0$ and $\phi_n(u)= \frac{1}{\Fbark(t_n)} \log (u/t_n) \bI_{u>t_n}$. Under condition $(\ref{Ordre1})$, we have 
\begin{eqnarray*}
\psi(\phi_n,u) & =& \gamma_{n,k}, \mbox{ if } u \leq t_n \\
 &  = &  \log \left( \frac{u}{t_n} \right) \fracFbark{u}{t_n}  +  \gamma_k \left( \frac{u}{t_n} \right)^{-1/ \gamma_k}  + \epsilon_n(u)  \left( \frac{u}{t_n} \right)^{-1/ \gamma_k + \delta} \mbox{ if } u > t_n, 
\end{eqnarray*}
where  $\epsilon_n(u)$  is  a sequence  tending  to $0$ uniformly in $u$, as $n\tinf$, and $\delta$ a positive real number such that $-\frac{1}{\gamma_k} + \delta<0$.  
\end{lem}
Proof :  We only consider the second situation where $u > t_n$ (the first one is straightforward) : 
\[
 \int_u^{+\infty} \phi_n(s) d F^{(k)}(x) = -  \int_{\frac{u}{t_n}}^{+\infty}  \log (y) \frac{d\Fbark (yt_n)}{\Fbark(t_n)} 
 \]
An integration by part and the fact that $\Fbark$ is regularly  varying at infinity with index $-1/ \gamma_k$, yields  
\[
 \int_u^{+\infty} \phi_n(s) d F^{(k)}(x) =  \log \left( \frac{u}{t_n} \right) \fracFbark{u}{t_n}  +  \gamma_k \left( \frac{u}{t_n} \right)^{-1/ \gamma_k}  + \Delta_n(u),
 \]
where 
\[
 \Delta_n(u)=  \int_{\frac{u}{t_n}}^{+\infty} \left(\frac{\Fbark (yt_n)}{\Fbark(t_n)}  - y^{-1/\gamma_k}  \right) \frac{dy}{y}
\]
Let $\delta$ be a positive real number. Then 
\begin{eqnarray*}
| \Delta_n(u)| &  = & \left|  \int_{\frac{u}{t_n}}^{+\infty}  y^{-1/\gamma_k-1 +\delta}  \left( y^{1/\gamma_k - \delta} \fracFbark{yt_n}{t_n}  -y^{-\delta}  \right) dy \right|  \\ 
& \leq & \sup_{y \geq 1}  \left| y^{1/\gamma_k - \delta} \fracFbark{yt_n}{t_n}  -y^{-\delta}  \right|  \int_{\frac{u}{t_n}}^{+\infty}  y^{-1/\gamma_k-1 +\delta} dy, 
\end{eqnarray*}
where the  function $y \rightarrow y^{1/\gamma_k-\delta} \Fbark(y)$  is regularly varying with index $-\delta$.  Then since
\[ 
 \sup_{y \geq 1}  \left| y^{1/\gamma_k - \delta} \fracFbark{yt_n}{t_n}  -y^{-\delta}  \right|  \tqdninf 0 
\]
and, when $-\frac{1}{\gamma_k} + \delta<0$, we have $ \int_{\frac{u}{t_n}}^{+\infty}  y^{-1/\gamma_k-1 +\delta} dy = cst \left( u/t_n \right)^{-1/ \gamma_k + \delta}$, this concludes the proof. 

\begin{lem} \label{lemmeConsiderationsProcEmp}
Recalling that $H$ is a distribution function with infinite right endpoint, we have :
\begitem
\item[$(i)$]  $\sup_{0\leq x<Z^{(n)}} \Hbar(x)/\Hbarn(x) = O_{\bP}(1)$
\item[$(ii)$] for any $a<1/2$, 
$$
 \sqrt{n} \sup_{t \geq 0} \frac{ |\Hbarn(t)-\Hbar(t)|}{(\Hbar(t))^a} = O_{\bP}(1) 
 \makebox[1.4cm][c]{and}
 \sqrt{n} \sup_{t \geq 0 } \frac{ |\HbarnO(t)-\HbarO(t)|}{(\HbarO(t))^a} = O_{\bP}(1) .
$$
\finit
\end{lem}

\noindent
Proof :  part $(i)$ is well known (see for instance section 3 of chapter 10 of \cite{ShorackWellner}), while the two statements in $(ii)$ are proved by usual empirical processes techniques, showing that the family of functions $(f_t)_{t<\infty}$ defined in one case by $f_t(z) = \bI_{z>t}/(\Hbar(t))^a$, and in the other case by $f_t(\delta,z) = (1-\delta) \bI_{z>t}/(\HbarO(t))^a$ are Donsker whenever $a<1/2$ (using respective square integrable envelope functions $f^*(z) = 1/(\Hbar(z))^a$ and $f^*(\delta,z)=(1-\delta)/(\HbarO(z))^{a}$, which bound from above the functions $f_t$ uniformly in $t$) .

\begin{lem} \label{lemGenRestes}
Under conditions (\ref{Ordre1}) and (\ref{condvntn}), suppose that   $\alpha \geq 0$ and $d \geq 1$ are real numbers. If   $\gamma_k < \gamma_C$ and
$$
\Xin =  \frac{ \sqrt{v_n} }{ n^{1+d}} \frac{ \phi(Z_i) }{ \Gbar(Z_i) (\HbarO(Z_i))^{d+\alpha} } \indksik , 
$$
then we have $ \sum_{i=1}^n  \Xin  \stackrel{\bP}{\longrightarrow} 0$, as $n$ tends to infinity, if $\alpha$ is $0$  or sufficiently close to it.
\end{lem}
Proof : 
\zdeux\\
According to the  LLN for triangular arrays, we need to prove the  following three statements : 
\begin{eqnarray*}
 (i) &  \forall \epsilon >0, \ \ \sum_{i=1}^n \bP(|\Xin| > \epsilon) & \tqdninf \; 0  \\
 (ii) &  \sum_{i=1}^n \bE((\Xin)^2 \bI_{|\Xin| \leq 1}) & \tqdninf \; 0  \\
  (iii) &  \sum_{i=1}^n \bE(\Xin \bI_{|\Xin| \leq 1}) &  \tqdninf \;  0  
\end{eqnarray*}
But, $\Xin$ being positive, $(iii)$ clearly implies $(ii)$. We thus need to prove that $(i)$ and $(iii)$ hold. 
\zdeux

Let us start with assertion $(i)$. If $\epsilon >0$ is given, then 
\[
 \Xin \; = \; \frac{v_n^{1/2}}{n^{1+d}}   \frac{\log (Z_i/t_n)}{\Fbarktn \Gbar(t_n) (\HbarO(t_n))^{d+\alpha}}   \indZitn  \indksik   \frac{\Gbar(t_n)}{\Gbar(Z_i)}  \left( \frac{ \HbarO(t_n)}{ \HbarO(Z_i)}  \right)^{d+\alpha} .
\]
Now, put $a= \frac{1}{\gamma_C} + \frac{d+\alpha}{\gamma} $ ($>0$); since, for a given $\epsilon'>0$,  there exists $c>0$ such that $\forall x \geq 1$, $\log (x) \leq c x^{\epsilon'} $,  and  using Potter-bounds $(\ref{BornesPotter})$ for $\Gbar^{-1} (\HbarO)^{-(d+\alpha)} \in RV_{-a}$, we can write (using the definition of $v_n$) 
 \begin{eqnarray*}
\{\, |\Xin| > \epsilon\,\}  & \subset &   \left\{  v_n^{-1/2} n^{-d}  \left( \HbarO(t_n)  \right)^{-(d+\alpha)} c (1+ \epsilon') \left(  \frac{Z_i}{t_n} \right)^{a+2\epsilon'}  > \epsilon  \right\}   \cap  \left\{ \xi_i=k  \mbox{ and } Z_i > t_n \right\} \\
& \subset &   \left\{  \, Z_i > c(\epsilon, \epsilon') \, t_n w_n  \, \right\}   \cap  \left\{  \xi_i=k   \mbox{ and } Z_i > t_n  \right\} , 
\end{eqnarray*}
where $w_n=\left( v_n^{1/2} n^d \left(  \HbarO(t_n)\right)^{d+\alpha} \right)^{1/(a+2\epsilon')}$ and $c(\epsilon,\epsilon')$ is a constant depending on $\epsilon$ and $\epsilon'$ only. Consequently, if $w_n$ tends to infinity,  
\begin{eqnarray*}
\sum_{i=1}^n \bP(|\Xin| > \epsilon) & \leq & n   \bE(\bI_{Z_i > c(\epsilon, \epsilon') \, t_n  w_n} \indksik)  = n \int_{c(\epsilon, \epsilon')\, t_n w_n}^{\infty} \Gbar(x) d F^k (x) \\
  & \leq & v_n \frac{(\Fbark \Gbar)(c(\epsilon, \epsilon')\, t_n w_n)}{(\Fbark \Gbar)(t_n)} \\
    & \leq & cst \ v_n w_n^{-\beta},
\end{eqnarray*}
where $\beta=\frac{1}{\gamma_C}+\frac{1}{\gamma_k} - \epsilon'$ and  the last inequality is due to Potter-bounds $(\ref{BornesPotter})$ applied to $\Fbark \Gbar \in RV_{-\frac{1}{\gamma_C}-\frac{1}{\gamma_k}}$. Then, assertion $(i)$ above will be true as soon as we prove  that $w_n\tinf$ and $v_n w_n^{-\beta}\ra 0$, as $n\tinf$. 

\noindent Since $\HbarO(t)$ is equivalent to a positive constant times $\Hbar(t)$ when $t\ra +\infty$,  and $\Hbar(t_n) \geq v_v /n$, then $w_n^{a+2\epsilon'} \geq cst \ (n^{-\eta} v_n)^r$, for $r= \frac 1 2 + d+\alpha >0$ and $\eta=\frac{\alpha}{r}\geq 0$.  Assumption (\ref{condvntn}) finally yields that $w_n$ tends to $+\infty$, since $0\leq \eta  \leq \eta_0$ for $\alpha$ sufficiently close to $0$.  \zun

\noindent Now, proving that  $v_n w_n^{-\beta}$ tends to $0$ is equivalent to proving that 
$v_n^{-(a+2\epsilon')/\beta} v_n^{1/2} n^d \left(  \HbarO(t_n)\right)^{d+\alpha}$ 
tends  to $+ \infty$. 
The same arguments as in the previous paragraph yield that it is sufficient to prove that $v_n^A n^{-\alpha} = \left( n^{-\eta} v_n \right)^A$ tends to $+ \infty$, for $A=-(a+2\epsilon')/\beta + 1/2 + d+\alpha$ and $\eta=\frac{\alpha}{A}$.  This is a consequence of hypothesis $(\ref{condvntn})$, since $A >0$  and $\alpha \leq \eta_0 A$, for  $\alpha$ sufficiently close to $0$.  This ends the proof of $(i)$.
\zdeux 

Let us now start the proof of assertion $(iii)$. If $\epsilon >0$ is given, using Potter-Bounds (\ref{BornesPotter}) for $\Gbar^{-1} (\HbarO)^{-(d+\alpha)}$ which belongs to $RV_{-a}$, and introducing $h(x)= \log (x) x^{a-\epsilon}$, we find that  (for some positive constant $c$)
\[
\bI_{|\Xin| \leq 1} \bI_{Z_i > t_n} \leq \bI_{h (Z_i / t_n) \leq c w_n} \bI_{Z_i > t_n}
\]
where we set $w_n=v_n^{1/2} n^d \left( \HbarO(t_n) \right)^{d+\alpha}$. Hence, denoting by $h^{-1}$ the inverse function of $h$, 
\[
\bI_{|\Xin| \leq 1} \bI_{Z_i > t_n} \indksik \leq \bI_{t_n < Z_i < t_n h^{-1} (c w_n)} \indksik.
\]
Consequently, using once again Potter-Bounds $(\ref{BornesPotter})$ and bounding the log with a constant times a power of $z/t_n$, we get
\begin{eqnarray*}
n \bE(X_{1,n} \bI_{|X_{1,n}| \leq 1} ) & \leq & \frac{v_n^{1/2}}{n^d} \int_{t_n}^{t_n h^{-1} (c w_n)} \frac{\log (z/t_n)}{\Fbark(t_n)  \Gbar(z) (\HbarO(z))^{d+\alpha}} dH^{(1,k)} (z) \\
& \leq & cst \ \frac{v_n}{w_n} \int_1^{ h^{-1} (c w_n)}  s^{b+2\epsilon'} \ \dfracFk{st_n},
\end{eqnarray*}
where $b=  \frac{d+\alpha}{\gamma}$ and $\epsilon'>0$ is some given positive value (the inequality $\log (s)\leq cst\, s^{\epsilon'}, \ \forall s\geq 1$, was used). But, by integration by parts and $(\ref{BornesPotter})$ applied to $\Fbark$,  setting $h_n= h^{-1} (c w_n)$, we have
\[
\frac{v_n}{w_n}  \int_1^{ h^{-1} (c w_n)} s^{b+2\epsilon'} \ \dfracFk{st_n}  
\leq 
cst \frac{v_n}{w_n} \left( 1 + h_n^{b-1/\gamma_k+3\epsilon'} \right) .
\]
Proceeding similarly as in the previous paragraphs, we find that $w_n/v_n\tinf$ (and thus $w_n$ and $h_n$ as well) thanks to assumption (\ref{condvntn}), for  $\alpha$ close to $0$. We are thus left to prove that $(v_n/w_n)\times h_n^{b'}$ tends to $0$, where $b'=b-1/\gamma_k+3\epsilon'$. If $b-1/\gamma_k$ is negative, this is immediate. We thus suppose that $b-1/\gamma_k\geq 0$ and, after some simple computations,  we find out that  $(v_n/w_n)h_n^{b'}$ tends to $0$ if $v_n^{-a+\epsilon'} w_n^{a-b'-\epsilon'}$ tends to $\infty$, a property which holds true thanks to assumption (\ref{condvntn}), for $\alpha$ close to $0$ (we omit the details). 

\zdeux

\begin{lem} \label{lem-integhetH}
Suppose that $V_1$ and $W_2$ are independent improper random variables of respective subdistribution functions $\HO$ and $\Hunk$, and $Z_3$ is independent of $V_1$ and $W_2$ and has distribution $H$. Consider $h$, $\het$, $\gdH$ and $\gdHet$ the functions defined in (\ref{def-fonctionsh}) and (\ref{def-gdHgdHet}).
\begitem
\item[$(i)$] For any $d\geq 1$, there exist some positive constants $c$ and $c'$ such that 
\[
 \bE\,(\, |\gdH^d(V_1,W_2)|\,) \leq c \,  \bE\,(\, h^d(V_1,W_2)\,) 
  \makebox[1.4cm][c]{and}
 \bE\,(\, |\gdHet^d(Z_3,V_1,W_2)|\,) \leq c' \,  \bE\,(\, \het^d(Z_3,V_1,W_2) \,)  . 
\]
\item[$(ii)$] For any $d\in ]1,1+(1+2\gamma_k/\gamma_C)^{-1}[$, we have 
$$
  \textstyle \bE\,(\, h^d(V_1,W_2)\,)  = O\left( (\Fbarktn\Gbar(t_n))^{2(1-d)} \right).
$$
In particular, if $\gamma_k < \gamma_C$, then $\bE(h^{4/3}(V_1,W_2))$ is of the order of $(\Fbarktn\Gbar(t_n))^{-2/3}$ and $\bE(h^d(V_1,W_2)) $ is finite whenever $d$ is (greater than but) sufficiently close  to $4/3$. 
\zun
\item[$(iii)$] For any $d\in ]1,1+(1+3\gamma_k/\gamma_C)^{-1}[$, we have 
$$
  \textstyle \bE\,(\, \het^d(Z_3,V_1,W_2)\,)  = O\left( (\Fbarktn\Gbar(t_n))^{3(1-d)} \right).
$$
In particular, if $\gamma_k < \gamma_C$, then $\bE(\het^{6/5}(Z_3,V_1,W_2))$ is of the order of $(\Fbarktn\Gbar(t_n))^{-3/5}$ and $\bE(\het^d(Z_3,V_1,W_2)) $ is finite whenever $d$ is (greater than but) sufficiently close to $6/5$. 
\item[$(iv)$] For any $d\in ]1/2,(2\gamma_C^{-1}+\gamma_F^{-1}+\gamma_k^{-1})/(3\gamma_C^{-1}+2\gamma_F^{-1})[$, we have $\bE\left(h^d(V_1,W_2)/\Hbar^d(V_1)\right)=O\left((\Fbarktn\Gbar(t_n))^{2-3d}\right)$. In particular, if $\gamma_k<\gamma_C$ then taking $\delta$ (greater than but) sufficiently close to $4/5$ is permitted, otherwise it is $2/3$ instead of $4/5$. 
\item[$(v)$] The integral $\theta_n=\iint h(v,w)d\HO(v)d\Hunk(w)$ is equivalent, as $n\tinf$, to $\gamma_k(-\log\Gbar(t_n))$.
\finit
\end{lem}

Proof : 
\begitem
\item[$(i)$] Let $d\geq 1$, and remind that $h$ is a non-negative function. Using several times the inequality $|a+b|^d\leq 2^{d-1}(|a|^d+|b|^d)$, we can write 
$$
 \bE ( |\gdH^d(V_1,W_2)| ) \leq cst \left\{ \bE(h^d(V_1,W_2)) + \bE[(\hunpt(V_1))^d] + \bE[(\hptun(W_2))^d] + ( \bE(h(V_1,W_2)) )^d  \right\}.
$$ 
But using the fact that the $L^1$ norm is bounded by the $L^d$ norm whenever $d\geq 1$, we have  $(\bE(h(V_1,W_2)) )^d \leq \bE(h^d(V_1,W_2))$ and it is quite simple to prove (by independency of $V_1$ and $W_2$) that it is also the case of $\bE[(\hunpt(V_1))^d] = \bE[(\bE(h(V_1,W_2)|V_1))^d] \leq \bE[\bE(h^d(V_1,W_2)|V_1)] = \bE(h^d(V_1,W_2))$, as well as for $\bE[(\hptun(W_2))^d]$. The inequality is thus proved. The other one (concerning $\gdHet$ and $\het$) is proved similarly.\zun

\item[$(ii)$] Let $d>1$. Since $h(v,\infty)=h(\infty,w)=0$ ($\forall v,w$), we have
\begin{eqnarray*}
 \bE(h^d(V_1,W_2)) 
  & = & 
 (\Fbarktn)^{-d} \iint \log^d(w/t_n) (\Hbar(v)\Gbar(w))^{-d} \bI_{w>t_n} \bI_{w>v} d\HO(v)d\Hunk(w) \\
  & = &
 (\Fbarktn)^{1-d} \int_{t_n}^{\infty} \log^d(w/t_n) \left(\int_0^w \frac{dG(v)}{\Gbar(v)\Hbar^{d-1}(v)} \right) \Gbar^{1-d}(w) \dfracFk{w}  \\
 & = & 
 \frac{ C_{d-1}(t_n) }{ (\Fbarktn\Gbar(t_n))^{d-1} } 
 \int_{t_n}^{\infty} \log^d(w/t_n) \left( \fracGbar{t_n}{w} \right)^{d-1} \frac{ C_{d-1}(w) }{ C_{d-1}(t_n) }  \dfracFk{w}
\end{eqnarray*}
where the function $C_{d-1}$ was defined in the statement of Lemma \ref{lem-Cdelta}. This lemma and Lemma \ref{lemtechInteg}, applied with $\alpha=d$, $a=(d-1)/\gamma_C+(d-1)/\gamma$ and $b=-1/\gamma_k$ (the  constraint specified on $d$ certifies that $a+b<0$), imply that the integral in the previous line converges to a constant. And Lemma \ref{lem-Cdelta} also implies that the ratio in front of this integral is equivalent, as $n\tinf$, to a positive constant times $\left(\Hbar(t_n)\Fbarktn\Gbar(t_n)\right)^{1-d}$, which is itself lower than $\left(\Fbarktn\Gbar(t_n)\right)^{2(1-d)}$, as desired. 
\zun
\item[$(iii)$] Let $d>1$. By definition of $\het$ in (\ref{def-fonctionsh}), and proceeding as in the previous item, $\bE(\het^d(Z_3,V_1,W_2))$ equals 
\begin{eqnarray*}
  &  & 
 (\Fbarktn)^{-d} \iiint \log^d\left({\textstyle\frac{w}{t_n}}\right) (\Hbar(v))^{-2d}(\Gbar(w))^{-d} \bI_{w>t_n} \bI_{w>v} \bI_{u>v} dH(u)  d\HO(v)d\Hunk(w) \\
 & = & 
 \frac{ C_{2d-2}(t_n) }{ (\Fbarktn\Gbar(t_n))^{d-1} } 
 \int_{t_n}^{\infty} \log^d(w/t_n) \left( \fracGbar{t_n}{w} \right)^{d-1} \frac{ C_{2d-2}(w) }{ C_{2d-2}(t_n) }  \dfracFk{w},
\end{eqnarray*}
which is equivalent to $O\left( (\Hbar(t_n))^{2-2d} (\Fbarktn\Gbar(t_n))^{1-d} \right) =   O\left( (\Fbarktn\Gbar(t_n))^{3(1-d)} \right)$ as soon as, thanks to Lemma \ref{lemtechInteg}, the sum $\left( (d-1)/\gamma_C+(2d-2)/\gamma \right) - 1/\gamma_k$ is negative, which turns out to be true whenever $d<1+(2+3\gamma_k/\gamma_C)^{-1}$, as specified.  
\zun
\item[$(iv)$] The proof is very similar to the previous ones, starting from
$$
 \bE\left(h^d(V_1,W_2)/\Hbar^d(V_1)\right) = 
 (\Fbarktn)^{-d} \iint \log^d(w/t_n) (\Hbar(v))^{-2d}(\Gbar(w))^{-d} \bI_{w>t_n} \bI_{w>v} d\HO(v)d\Hunk(w) 
$$ 
so we omit the details. 
\item[$(v)$] Noting that $-\log\Gbar$ is slowly varying at infinity null at $0$, we have
$$
 \theta_n = \int_{t_n}^{\infty} \log(w/t_n) \left( \int_0^w dG(v)/\Gbar(v) \right) \dfracFk{w}  
 =  (-\log \Gbar(t_n)) \left( - \int_1^{\infty} \log(u) \frac{ -\log\Gbar(ut_n) }{ -\log\Gbar(t_n) } \dfracFbark{ut_n}\right) 
$$
which can be dealt with using part $(ii)$ of Lemma \ref{lemtechInteg} with $\alpha=1$, $a=0$ and $b=-1/\gamma_k$ : the obtained constant is indeed equal to $\gamma_k$. 
\finit

\begin{lem} \label{diversPreuveUstats}
In this Lemma, various notations defined in sections \ref{subsec-prelimsRn2Rn3} to \ref{subsec-preuveVn} are used.
\begitem
\item[$(i)$] The variables $\gdHss_I$ for $I\in \{(i,j) \, ; \, 1\leq i<j\leq n\}$ are centred and uncorrelated  . This is also true for  the variables $\gdHetss_I$ for $I\in \{(i,j,l) \, ; \, 1\leq i<j<l\leq n\}$.
\item[$(ii)$] We have $\bE\left[ (\gdHss(V_1,W_2))^2\right] \leq  48 \bE[\gdHun^2\bI_{|\gdHun|\leq M_n}]$.
\item[$(iii)$] We have $\bE\left(\, |\gdHun-\gdH^*(V_1,W_2)+\gdHsunpt(V_1)+\gdHsptun(W_2)| \,\right) 
 \leq  4\bE\left(|\gdHun|\bI_{|\gdHun|>M_n}\right)$
\finit
\end{lem}

\noindent
Proof : 
\begitem
\item[$(i)$] Let us consider the first situation, where ${\cal I}=\{(i,j) \, ; \, 1\leq i<j\leq n\}$. First, if  $I=(i,j)\in{\cal I}$, then $\bE(\gdHss_I) = 0 - \bE(\gdHsunpt(V_i))-\bE(\gdHsptun(W_j))$ ; but, by definition of $\gdHsunpt$ and independency of $V_i$ and $W_j$, we have $\bE(\gdHsunpt(V_i))=\bE(\gdHs(V_i,W_j))=0$,  and $\bE(\gdHsptun(W_j))=0$ is obtained similarly,  so we proved that $\bE(\gdHss_I)=0$. Note that we can prove (with similar arguments)  that $\gdHssunpt(v)=\gdHssptun(w)=0$ for every $v,w$ in $[0,\infty]$, a property which is repeatedly used below . Let us now deal with the non-correlation of $\gdHss_I$ and $\gdHss_{I'}$,  by considering the various cases where $I\neq I'$ with $I=(i,j)$ and $I'=(k,l)$ are in ${\cal I}$.\zun

If all four indices $i,j,k,l$ are distinct, then non-correlation of $\gdHss_I$ and $\gdHss_{I'}$ is immediate by mutual independence of the variables $Z_1,\ldots,Z_n$.

 \noindent  If $i=k$ but $j\neq l$, then $\bE(\gdHss_I\gdHss_{I'})=\bE(\psi(V_i))$ where $\psi(v)=\bE(\gdHss(v,W_j)\gdHss(v,W_l))=(\gdHssunpt(v))^2=0$, by independence of $V_i$ with $(W_j,W_l)$, and of $W_j$ and $W_l$.
 
  \noindent  The case $i\neq k$ and $j=l$ is similar using $\gdHssptun(\cdot)\equiv 0$.  
  
  \noindent If $i=l$ but $j\neq k$, then $\bE(\gdHss_I\gdHss_{I'})=\bE(\psi(V_i,W_i))$ where $\psi(v,w)=\bE(\gdHss(v,W_j)\gdHss(V_k,w))=\gdHssunpt(v)\gdHssptun(w)=0\times 0=0$ ; the case $j=k$ and $i\neq l$ is treated similarly. 
  
  \noindent Note that the case $i=l$ and $j=k$ ({\it i.e.} $\gdHss_I=\gdHss(V_i,W_j)$, $\gdH_{I'}=\gdHss(V_j,W_i)$) is not permitted (it would lead to dependency) since we cannot have simultaneously $i<j$ and $j<i$ ; this is the reason why, in the beginning of section \ref{subsec-preuveUn}, we restricted the study of the sum $\calUn$ to that of the sum $S_N$ having terms $\gdH(V_i,W_j)$ satisfying $i<j$.\zun

The second situation, for $\gdHetss_I$ and $\gdHetss_{I'}$ with $I\neq I'$ in ${\cal I}=\{I=(i,j,l) \, ; \, 1\leq i<j<l\leq n\}$, is a bit more tedious (with more cases to detail) but very similar, so we omit its proof.

\item[$(ii)$] We start by the trivial bound 
$$
 \bE[(\gdHss(V_1,W_2))^2] 
 \leq 4\left\{ \bE[(\gdHs(V_1,W_2))^2] + \bE[(\gdHsunpt(V_1))^2] + \bE[(\gdHsptun(W_2))^2] \right\}.
$$
Noting $\gdHun^-=\gdHun\bI_{|\gdHun|\leq M_n}$, we can write, on one hand, by definition of $\gdHs$,  $\bE[(\gdHs(V_1,W_2))^2] \leq 2\left\{\bE[(\gdHun^-)^2] + (\bE[\gdHun^-])^2 \right\} \leq 4\bE[(\gdHun^-)^2]$.  On the other hand, if $W$ is independent of $V_1$, we have $\bE[(\gdHsunpt(V_1))^2] = \bE[(\bE[\gdHs(V_1,W) | V_1])^2] \leq \bE[\bE[ (\gdHs(V_1,W))^2 | V_1]]=\bE[(\gdHs(V_1,W_2))^2]$, which is the same term as the first one, and is thus lower than $4\bE[(\gdHun^-)^2]$. The same is true of   $\bE[(\gdHsptun(W_2))^2]$, so the desired inequality is proved. 

\item[$(iii)$] First recall that $\gdH_1$ denotes $\gdH(V_1,W_2)$. Now, since $\gdH_1$ is centred and we trivially have $\gdH_1=\gdH_1\bI_{|\gdH_1|\leq M_n} + \gdH_1\bI_{|\gdH_1|> M_n}$, noting $\gdH_1^+ = \gdH_1\bI_{|\gdH_1|>M_n}$ yields 
$$\gdH_1-\gdHs(V_1,W_2)=\gdH_1^{+} - \bE(\gdH_1^+).$$
 Secondly, using the fact that $\gdHunpt(\cdot)\equiv 0$ (simple to prove), we can write 
$$
 \gdHsunpt(v)= \bE(\gdH(v,W)\bI_{|\gdH(v,W)|\leq M_n}) - \bE(\gdH_1\bI_{|\gdH1|\leq M_n}) 
  = - \gdHunpt^+(v) + \bE(\gdH_1^+) ,
$$
where $ \gdHunpt^+(v)$ denotes $\bE(\gdH(v,W)\bI_{|\gdH(v,W)|>M_n})$ and satisfies $\bE(\gdHunpt^+(V_1))=\bE(\gdHun^+)$,  and similarly 
$$\gdHsptun(w)= - \gdHptun^+ (w) + \bE(\gdH_1^+)$$
 with $\gdHptun^+(w)=\bE(\gdH(V,w)\bI_{|\gdH(V,w)|>M_n})$ and $\bE(\gdHptun^+(W_2))=\bE(\gdHun^+)$. Summing these three terms  finally leads to 
\begin{eqnarray*}
 \bE\left(\, |\gdHun-\gdH^*(V_1,W_2)+\gdHsunpt(V_1)+\gdHsptun(W_2)| \,\right) 
 & = & 
 \bE\left(\, |\gdH_1^+  -\gdHunpt^+(V_1) -\gdHptun^+(W_2) + \bE(\gdH_1^+) | \,\right)
\end{eqnarray*}
 which is lower than $4\bE(|\gdHun^+|)$, as announced.  
\finit

\vspace{1.cm}


\end{document}